\documentclass[times,final]{elsarticle}

\usepackage[margin=2.5cm]{geometry}
\usepackage{setspace}
\usepackage[tableposition=above]{caption}
\usepackage[]{booktabs}

\usepackage{graphicx} 
\usepackage[section]{placeins} 

\usepackage[]{xfrac} 
\usepackage[]{amsmath} 
\usepackage[bb = dsserif]{mathalpha} 

\usepackage[colorlinks]{hyperref}
\usepackage[capitalise]{cleveref} 
\biboptions{sort&compress}
\begin{document}

\begin{frontmatter}

\title{An adaptive wavelet method for nonlinear partial differential equations with applications to dynamic damage modeling}

\author[1]{Cale Harnish}
\author[2]{Luke Dalessandro}
\author[1]{Karel Matouš\corref{cor1}}
\cortext[cor1]{Corresponding author}
\ead{kmatous@nd.edu}
\author[3]{Daniel Livescu}

\address[1]{Department of Aerospace and Mechanical Engineering, University of Notre Dame, Notre Dame, Indiana 46556, USA}
\address[2]{Department of Intelligent Systems Engineering, Indiana University, Bloomington, Indiana 47405, USA}
\address[3]{Computer and Computational Sciences Division, Los Alamos National Laboratory, Los Alamos, New Mexico 87545, USA}

\begin{abstract}
Multiscale and multiphysics problems need novel numerical methods in order for them to be solved correctly and predictively. To that end, we develop a wavelet based technique to solve a coupled system of nonlinear partial differential equations (PDEs) while resolving features on a wide range of spatial and temporal scales. The algorithm exploits the multiresolution nature of wavelet basis functions to solve initial-boundary value problems on finite domains with a sparse multiresolution spatial discretization. By leveraging wavelet theory and embedding a predictor-corrector procedure within the time advancement loop, we dynamically adapt the computational grid and maintain accuracy of the solutions of the PDEs as they evolve. Consequently, our method provides high fidelity simulations with significant data compression. We present verification of the algorithm and demonstrate its capabilities by modeling high-strain rate damage nucleation and propagation in nonlinear solids using a novel Eulerian-Lagrangian continuum framework.
\end{abstract}


\end{frontmatter}
\section{Introduction}
This work details a novel numerical method capable of modeling multiphysics problems while resolving features on a wide range of spatial and temporal scales. Considering that many engineering applications must incorporate failure modes into material models, the potential of our algorithm is demonstrated by modeling the dynamic damage of a nonlinear solid. The behavior of damaged material has been the subject of extensive research \cite{historyIrwin,historyDugdale,historyWillis,historySpanoudakis,historyBazant2,historyJu,historyNeedleman,historyFish,historyBazant,historyBuehler,historyVoyiadjis,historyKitey,historyPark,historyHuynh,historyKim}. Early damage models were simple in nature and classified material failure based on comparing the state of stress or strain to empirically measured values \cite{Hankinsonphenom,Hillphenom,DPphenom,TWphenom,HBphenom,JHphenom,MCphenom}. Fracture mechanics provides a more theoretical approach and analytically describes the forces required to open and grow a crack \cite{Griffith,historyBarenblatt,Freund}. Alternatively, this work applies thermodynamically consistent computational damage mechanics to model the effects of damage \cite{historyKrajcinovic}.

Continuum damage mechanics has been used to model several failure modes such as: fatigue \cite{fatigue1971}, creep \cite{creep1974}, ductile damage \cite{ductile1985,ductile1985b,ductile1985c}, brittle damage \cite{brittle1981,brittle1983,brittle1985}, isotropic damage \cite{isotropic1984}, and anisotropic damage \cite{anisotropic1980,anisotropic1981}. This work leverages concepts from the reference map technique \cite{RMT} to extend the Lagrangian constitutive theory from \cite{matous2008,matous2010,matous2015} into the Eulerian frame and model rate-dependent isotropic damage. Our coupling of the Eulerian-Lagrangian frames combines the advantages of resolving large deformations on an Eulerian discretization with a constitutive damage model from a Lagrangian description and is one of the novelties of this work. The derivation of the energy based damage model is thermodynamically motivated and incorporates the irreversible nature of the damage process \cite{historySimoJu,historySimoJu2}.

In general, material failure does not occur only between two adjacent atomic layers \cite{historyBuehler}. Instead, cracks develop with a finite thickness of a characteristic length-scale, $l_{\omega}$, and have a process zone at the crack-tip \cite{zoneXue,zoneSteel,zoneConcrete,zoneRubber}. At the microscale one may find voids and material defects, microcracks can be observed at the mesoscale, and at the macroscale these imperfections coalesce into macrocracks \cite{damageScales}. These features have very small spatial and temporal scales which prove challenging for both experimental measurement \cite{expPearson,expDekkers,expHuang,expKinlock,expRadford,expNakamura,expSingh} and computational modeling \cite{matous2008,compChallenges2}. However, computational models of damage have become increasingly important to bridge the gap between theoretical predictions and experimental design.

The computational difficulty of damage modeling has so far been addressed through a variety of numerical techniques. For example, the integral equations from peridynamics are used in \cite{peridynamicDamage} to better model the discontinuous nature of cracks. A mesh-free smooth particle hydrodynamic (SPH) simulation was used in \cite{SPHdamage} to model damage under compressive loading. Finite element methods (FEM) have incorporated rate dependent damage models \cite{FEMdamageA} and cohesive zone models \cite{FEMdamageB} to simulate fracture in brittle materials (\emph{e.g.}, concrete, ceramics, PMMA). The extended finite element method (XFEM) has been used \cite{XFEMdamage} to model a crack with arbitrary discontinuities instead of remeshing. Alternatively, \cite{matous2003,matous2021} only resolves the damaged material on the microscale and assumes a separation of scales to model its effect on the macroscale using computational homogenization \cite{matous2008,matous2009,matous2010,matous2015}. In contrast, an adaptive space-time discontinuous Galerkin method is used in \cite{vspike2010} to provide high resolution simulations, revealing new solution features (\emph{e.g.}, a quasi-singular structure in the velocity response). 

In other works, the challenge of multiscale modeling has been met through adaptive mesh refinement (AMR) \cite{AMR_def, AMR_hesthaven}, multigrid methods \cite{multigrid_def, multigrid_hackbusch,matousDewen}, Chimera overset grids \cite{chimera}, or remeshing/refining FEM \cite{karniadakis, FEM_hp1, FEM_hp2, FEM_rh}. However, these approaches are computationally expensive in the context of damage modeling since the spatial and temporal location of interesting solution features are not known \emph{a priori}. This work develops an alternative adaptive discretization strategy to resolve all relevant scales from the characteristic length scale of damage to the size of the macroscopic material. In particular, we use wavelet basis functions to approximate solutions to coupled systems of nonlinear partial differential equations (PDEs) in a monolithic multiscale solver.

Wavelet based numerical methods are particularly well suited for modeling multiscale and multiphysics problems, because they can dynamically adapt the discretization to discover spatial scales \cite{Jawerth1994, Schneider2010}. Furthermore, current wavelet solvers have achieved several notable accomplishments, including: significant data compression \cite{Liandrat1990, Beylkin1997, Bertoluzza1996a}, bounded energy conservation \cite{Ueno2003, Qian1993}, modeling stochastic systems \cite{Kong2016}, multiscale model reduction \cite{Rody2019}, and solving coupled systems of nonlinear PDEs \cite{Paolucci2014PT1, Paolucci2014PT2, Nejadmalayeri2015, Dubos2013, Sakurai2017}. However, some implementations only solve PDEs in infinite or periodic domains (\emph{e.g.}, \cite{Frohlich1994, Goedecker, Iqbal2014}), some do not utilize the data compression ability of wavelets, resulting in a costly uniform grid (\emph{e.g.}, \cite{Qian1993, Le2015, Lin2001}), and some use finite difference methods to calculate the spatial derivatives, reducing the ability to control accuracy (\emph{e.g.}, \cite{Paolucci2014PT1, Paolucci2014PT2, Nejadmalayeri2015, Holmstrom1999}). 

We present a predictor-corrector algorithm designed to advance the state of wavelet-based methods by leveraging the successes of past solvers and developing strategies to mitigate their limitations. We use differentiable wavelet bases to numerically solve nonlinear initial–boundary value problems. We use second-generation wavelets near spatial boundaries to produce multiresolution spatial discretizations on finite domains. Furthermore, the lowest resolution is defined with the minimum number of collocation points required to support the basis functions, thereby improving the ability to compress data. Additionally, we compute spatial derivatives by operating directly on the wavelet bases. We use error estimates for the wavelet representation of fields, their spatial derivatives, and the aliasing errors associated with nonlinear terms to construct a sparse, dynamically adaptive computational grid for each unknown that provides the required accuracy and approaches the theoretical \emph{a priori} guarantee. Moreover, we employ an embedded Runge-Kutta integration to control the temporal accuracy as the solution of the PDE evolves. Therefore our adaptive multiresolution discretization is guided directly by the wavelet operators and is one of the novelties of this work.

An overview of wavelet theory and the wavelet operations required to solve PDEs is presented in \cref{sec:wavelets}. Then, our algorithm, the Multiresolution Wavelet Toolkit (MRWT) is described and verification examples are presented in \cref{sec:MRWT}. Next, the Eulerian-Lagrangian model of dynamic damage growth in nonlinear materials is derived in \cref{sec:model}. Finally, \cref{sec:damage_results} shows results from modeling high-strain rate damage nucleation and propagation, illustrating the benefits of applying this novel numerical method to multiscale and multiphysics problems.
\section{Multiresolution analysis}
\label{sec:wavelets}
The proposed numerical method exploits the properties of wavelet basis functions and builds the computational grid using the wavelet operators. Therefore it is essential to provide a brief summary of wavelet theory. In particular, the following outlines the construction of wavelet basis functions and defines the operations needed to build a sparse grid and solve nonlinear PDEs.

\subsection{Wavelet basis functions}
\label{sec:basis}
A multiresolution analysis (MRA) provides the formal mathematical definition of a wavelet family of basis functions \cite{Daubechies10}. An MRA of an $N$-dimensional domain $\Omega \subset \mathbb{R}^{N}$ consists of a sequence of nested approximation spaces $\boldsymbol{V}_{j}$ and their associated dual spaces $\widetilde{\boldsymbol{V}}_{j}$ such that the union of these spaces is the $L^{2}(\Omega)$ space \cite{Cohen2000_interval}. The wavelet spaces $\boldsymbol{W}_{j}$, and their associated dual spaces $\widetilde{\boldsymbol{W}}_{j}$ are then defined as the complements of the approximation spaces $\boldsymbol{V}_{j}$ in $\boldsymbol{V}_{j+1}$ \cite{Bacry1992,Qian1993}
\begin{align} 
\label{eqn:MRAv}
    \boldsymbol{V}_{j} &\subset \boldsymbol{V}_{j+1},
    &
    \overline{\bigcup_{j = 1}^{\infty} \boldsymbol{V}_{j}} &= L^{2}(\Omega),
    &
    \boldsymbol{V}_{j+1} &= \boldsymbol{V}_{j} \oplus \boldsymbol{W}_{j}.
\end{align}
Furthermore, $N$-dimensional representations can be constructed through $N$ tensor products of the one-dimensional ($1$D) spaces $V_{j}$ and $W_{j}$. For example, the approximation spaces in a three-dimensional ($3$D) setting become
\begin{align}
\label{eqn:MRA_3D}
    \boldsymbol{V}_{j} &= V_{j} \otimes V_{j} \otimes V_{j},\nonumber\\
    \boldsymbol{V}_{j} &= \left( V_{j-1} \oplus W_{j-1} \right) \otimes \left( V_{j-1} \oplus W_{j-1} \right) \otimes \left( V_{j-1} \oplus W_{j-1} \right).
\end{align}
In general, this produces $2^{N}$ types of basis from each possible combination of $V_{j}$ and $W_{j}$ spaces. The scaling functions $\phi_{k}^{j} (x)$ and dual scaling functions $\widetilde{\phi}_{k}^{j} (x)$ are the $1$D basis functions in the spaces $V_{j}$ and $\widetilde{V}_{j}$ respectively, whereas the wavelets $\psi_{k}^{j} (x)$ and dual wavelets $\widetilde{\psi}_{k}^{j} (x)$ are the $1$D basis functions in the spaces $W_{j}$ and $ \widetilde{W}_{j}$ respectively. Note that the multiresolution nature of $1$D wavelets requires the use of two types of indices. One to define the resolution level $j \in \mathbb{Z}$, and another to define the spatial locations $k \in \mathbb{Z}$ on a particular resolution level $j$. Furthermore, any basis with $j$ and $k$ designation is constructed through dilations and translations such that 
\begin{align}
\label{eqn:trans_dia}
    \phi_{2k}^{j + 1}(x) &= \phi_{k}^{j} (2x - 2k \Delta x_{j+1}),
    &
    \widetilde{\phi}_{2k}^{j + 1}(x) &= 2 \ \widetilde{\phi}_{k}^{j} (2x - 2k \Delta x_{j+1}), \nonumber
    \\
    \psi_{2k}^{j + 1}(x) &= \psi_{k}^{j} (2x - 2k \Delta x_{j+1}),
    &
    \widetilde{\psi}_{2k}^{j + 1}(x) &= 2 \ \widetilde{\psi}_{k}^{j} (2x - 2k \Delta x_{j+1}),
\end{align}
where $\Delta x_{j+1}$ is the grid spacing on level $j + 1$. An additional index $\lambda \in \mathbb{Z}$ is introduced for multidimensional wavelets to simplify notation and denote specific $V_{j}$ and $W_{j}$ combinations as shown in \cref{tab:basis}. The local index $k$ is reset for each $j$ and $\lambda$ combination.
\begin{table}[htb]
    \begin{center}
        \caption{Notation for multidimensional basis in $\mathbb{R}^{3}$.\label{tab:basis}}
        \begin{tabular}{ c | c }
            \toprule
            ${}^{\lambda}\Psi_{\vec{k}}^{j}(\vec{x})$ & Three-dimensional basis\\ 
            \midrule
            ${}^{0}\Psi_{\vec{k}}^{j}(\vec{x})$ & $\phi_{k_{1}}^{j}(x_{1}) \ \phi_{k_{2}}^{j}(x_{2}) \ \phi_{k_{3}}^{j}(x_{3})$
            \\
            ${}^{1}\Psi_{\vec{k}}^{j}(\vec{x})$ & $\psi_{k_{1}}^{j}(x_{1}) \ \phi_{k_{2}}^{j}(x_{2}) \ \phi_{k_{3}}^{j}(x_{3})$
            \\
            ${}^{2}\Psi_{\vec{k}}^{j}(\vec{x})$ & $\phi_{k_{1}}^{j}(x_{1}) \ \psi_{k_{2}}^{j}(x_{2}) \ \phi_{k_{3}}^{j}(x_{3})$
            \\
            ${}^{3}\Psi_{\vec{k}}^{j}(\vec{x})$ & $\psi_{k_{1}}^{j}(x_{1}) \ \psi_{k_{2}}^{j}(x_{2}) \ \phi_{k_{3}}^{j}(x_{3})$
            \\
            ${}^{4}\Psi_{\vec{k}}^{j}(\vec{x})$ & $\phi_{k_{1}}^{j}(x_{1}) \ \phi_{k_{2}}^{j}(x_{2}) \ \psi_{k_{3}}^{j}(x_{3})$
            \\
            ${}^{5}\Psi_{\vec{k}}^{j}(\vec{x})$ & $\psi_{k_{1}}^{j}(x_{1}) \ \phi_{k_{2}}^{j}(x_{2}) \ \psi_{k_{3}}^{j}(x_{3})$
            \\
            ${}^{6}\Psi_{\vec{k}}^{j}(\vec{x})$ & $\phi_{k_{1}}^{j}(x_{1}) \ \psi_{k_{2}}^{j}(x_{2}) \ \psi_{k_{3}}^{j}(x_{3})$
            \\
            ${}^{7}\Psi_{\vec{k}}^{j}(\vec{x})$ & $\psi_{k_{1}}^{j}(x_{1}) \ \psi_{k_{2}}^{j}(x_{2}) \ \psi_{k_{3}}^{j}(x_{3})$ \\ 
            \bottomrule
        \end{tabular}
    \end{center}
\end{table}

In many cases the wavelet basis functions do not have a closed-form expression, instead they are defined in terms of four types of filter coefficients (\emph{i.e.}, $h_{i}, \ \widetilde{h}_{i}, \ g_{i},$ and $\widetilde{g}_{i}$) \cite{Goedecker,interval} which satisfy orthogonality relations
\begin{align}
\label{eqn:ortho}
    h_{i} \ \widetilde{h}_{i} &= 1, & g_{i} \ \widetilde{g}_{i} &= 1, \nonumber\\
    h_{i} \ \widetilde{g}_{i} &= 0, & g_{i} \ \widetilde{h}_{i} &= 0,
\end{align}
symmetry relations
\begin{align}
\label{eqn:sym}
    g_{i+1} &= (-1)^{i+1} \ \widetilde{h}_{-i}, & \widetilde{g}_{i+1} &= (-1)^{i+1} \ h_{-i},
\end{align}
and refinement relations
\begin{align}
\label{eqn:ref}
    \phi_{k}^{j}(x) &= \sum_{i = -p}^{p} h_{i} \ \phi_{2 k + i}^{j + 1} (x),
    &
    \widetilde{\phi}_{k}^{j}(x) &= \sum_{i = -p}^{p} \widetilde{h}_{i} \ \widetilde{\phi}_{2 k + i}^{j + 1} (x), \nonumber
    \\
    \psi_{k}^{j}(x) &= \sum_{i = -p}^{p} g_{i} \ \phi_{2 k + i}^{j + 1} (x),
    &
    \widetilde{\psi}_{k}^{j}(x) &= \sum_{i = -p}^{p} \widetilde{g}_{i} \ \widetilde{\phi}_{2 k + i}^{j + 1} (x).
\end{align} 
The parameter $p$ is defined by the length of the filter coefficients and reflects the properties of the basis (\emph{e.g.}, the number of vanishing moments or the degree of continuity) \cite{Goedecker}. Typically the wavelet basis are defined for infinite or periodic domains and modifications are necessary for representation on finite domains \cite{Beylkin1992, Sweldens1998}.

There are many families of wavelet basis functions which satisfy \crefrange{eqn:MRAv}{eqn:ref} and this work in particular uses the Deslauriers-Dubuc wavelets \cite{wavelet_FEM_beam, maxwell}. These basis functions are sometimes referred to as the biorthogonal interpolating wavelet family \cite{donoho_interpolating}, or the auto-correlation of the Daubechies wavelets \cite{Bertoluzza1996a}. Additionally, second generation wavelets are used near spatial boundaries \cite{interval}.

\subsection{Wavelet transformations}
\label{sec:transforms}
Space is discretized by projecting each continuous function $f(\vec{x})$, defined on the $N-$dimensional finite domain $\Omega \subset \mathbb{R}^{N}$, onto the basis functions ${}^{\lambda}\Psi_{\vec{k}}^{j}(\vec{x})$. With an infinite number of resolution levels the representation is exact and, due to the nested nature of the approximation spaces, the $\lambda=0$ basis are needed only on the lowest resolution
\begin{align}
\label{eqn:fullTransform}
    f(\vec{x}) = \sum_{k_{i} \in [0, 2p]} {}^{0}\mathbb{d}_{\vec{k}}^{1} \ {}^{0}\Psi_{\vec{k}}^{1}(\vec{x}) + \sum_{j=1}^{\infty} \sum_{\lambda = 1}^{2^{N}-1} \sum_{\vec{k}} {}^{\lambda}\mathbb{d}_{\vec{k}}^{j} \ {}^{\lambda}\Psi_{\vec{k}}^{j}(\vec{x}).
\end{align}
Note that in \cref{eqn:fullTransform} the lowest resolution level (\emph{i.e.}, $\lambda = 0$ and $j = 1$) is defined by a uniform grid of $2p + 1$ points in each direction. This choice maximizes the compression ability of wavelets since it is the minimum number of points required to support one ${}^{0}\Psi_{\vec{k}}^{1}(\vec{x})$ basis at the center of the domain (see \cref{tab:basis}). All other points on the lowest resolution level (coarse grid) correspond to the second generation wavelets defined for a finite domain. The coefficients ${}^{\lambda}\mathbb{d}_{\vec{k}}^{j}$ are defined by integrating the field with the dual basis
\begin{align}
\label{eqn:coeff_def}
   {}^{\lambda}\mathbb{d}_{\vec{k}}^{j} ( f ) = \int f(\vec{x}) \ {}^{\lambda}\widetilde{\Psi}_{\vec{k}}^{j}(\vec{x}) \ \mathrm{d} \Omega.
\end{align}
Using the Deslauriers-Dubuc wavelet family, the integrals in \cref{eqn:coeff_def} are solved exactly and are replaced with the matrix operator $\mathbbb{F}$, defined in terms of the filter coefficients, as described in \cite{Harnish2018,Harnish2021}. Repeated application of this operator yields all of the wavelet coefficients on each resolution level. For continuous functions $f(\vec{x})$, the magnitude of coefficients ${}^{\lambda}\mathbb{d}_{\vec{k}}^{j}$ decreases as the resolution level increases. Therefore, the infinite sum in \cref{eqn:fullTransform} is truncated by limiting the number of resolution levels to $j_{\mathrm{max}} \in \mathbb{Z}$ such that the magnitude of all wavelet coefficients above $j_{\mathrm{max}}$ are below a prescribed tolerance $\varepsilon$. This truncated wavelet projection results in the discretization
\begin{align}
\label{eqn:fJ}
    f_{j}(\vec{x}) = \sum_{k_{i} \in [0, 2p]}  {}^{0}\mathbb{d}_{\vec{k}}^{1} \ {}^{0}\Psi_{\vec{k}}^{1}(\vec{x}) + \sum_{j=1}^{j_{\mathrm{max}}} \sum_{\lambda = 1}^{2^{N}-1} \sum_{\vec{k}} {}^{\lambda}\mathbb{d}_{\vec{k}}^{j} \ {}^{\lambda}\Psi_{\vec{k}}^{j}(\vec{x}),
\end{align}
defined on a dense grid with uniform spacing given by $j_{\mathrm{max}}$. However this discretization can be made sparse by leveraging the unique mapping between the wavelet indices and the $N$-dimensional spatial coordinate of a collocation point in the domain given by $\vartheta : (j, \lambda, \vec{k}) \rightarrow \vec{x} \in \mathbb{R}^N$
\begin{align}
\label{eqn:xMap}
    x_{i} &= a_{i} + (2 k_{i} + \mathrm{bin}(\lambda)_{i}) \frac{b_{i} - a_{i}}{p \ 2^{j+1}}.
\end{align}
In the above, $a_{i}$ and $b_{i}$ are respectively the lower and upper bounds of the $i^{\mathrm{th}}$ direction in the $N$-dimensional finite domain. Additionally, $\mathrm{bin}(\lambda)_{i}$ corresponds to the $i^{\mathrm{th}}$ bit of the binary representation of $\lambda$ (\emph{e.g.}, in $\mathbb{R}^3$, $\mathrm{bin}(\lambda = 6) = 011$). The inverse mapping, $\vartheta^{-1} : \vec{x} \in \mathbb{R}^N \rightarrow (j, \lambda, \vec{k})$, requires first specifying $j$, then $\lambda$, then $\vec{k}$
\begin{align}
    \label{eqn:invMap}
    j &= \max \left( j_{i} \right), \qquad \mathrm{with} \qquad j_{i} = \min \left( \left\{ j \in [1,j_{\mathrm{max}}] :  x_{i} - a_{i} \bmod \frac{b_{i} - a_{i}}{p 2^{j + 1}} = 0 \right\} \right),    
    \nonumber \\
    \mathrm{bin}(\lambda)_{i} &= \frac{\left( x_{i} - a_{i} \right) p 2^{j + 1} }{b_{i} - a_{i}} \bmod 2,
    \nonumber \\
    k_{i} &= \frac{\left( x_{i} - a_{i} \right) p 2^{j} }{b_{i} - a_{i}} - \frac{\lambda_{i}}{2}.
\end{align}

Therefore, the omission of small wavelet coefficients on resolution levels $j \in [1, j_{\mathrm{max}}]$ corresponds to the omission of collocation points in the computational domain and this procedure results in 
\begin{align}
\label{eqn:fEpsilon}
    f_{\varepsilon}(\vec{x}) = \sum_{k_{i} \in [0, 2p]}  {}^{0}\mathbb{d}_{\vec{k}}^{1} \ {}^{0}\Psi_{\vec{k}}^{1}(\vec{x}) + \sum_{j=1}^{j_{\mathrm{max}}} \sum_{\lambda = 1}^{2^{N}-1} \sum_{\substack{\vec{k} \\ |\mathbb{d}| \geq \varepsilon}} {}^{\lambda}\mathbb{d}_{\vec{k}}^{j} \ {}^{\lambda}\Psi_{\vec{k}}^{j}(\vec{x}),
\end{align}
defined on a sparse, multiresolution spatial discretization with \emph{a priori} knowledge of the spatial error \cite{Paolucci2014PT1, Nejadmalayeri2015, Holmstrom1999}
\begin{align}
\label{eqn:ferror}
    ||f(\vec{x}) - f_{\varepsilon}(\vec{x})||_{\infty} \leq \mathcal{O}(\varepsilon).
\end{align}

Using the Deslauriers-Dubuc wavelet family, the wavelet coefficients are mapped back to field values using the matrix operator $\mathbbb{B} = \mathbbb{F}^{-1}$, also defined in terms of the filter coefficients as described in \cite{Harnish2018,Harnish2021}. The structure of these matrices are similar to those used in \cite{Goedecker,nasa_matrices}. Though like \cite{Dahmen1999}, we modify these matrices near spatial boundaries to account for a finite domain.

\begin{figure}[!htb]
\centering
    \includegraphics[width=0.9\textwidth]{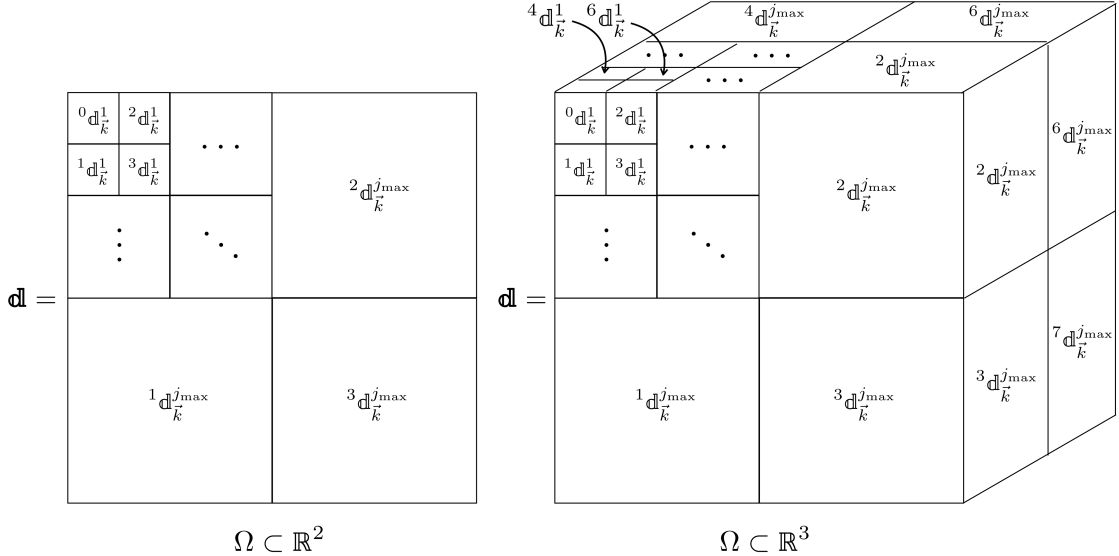}
\caption{Array of $\mathbbb{d}$ in $2$D and $3$D settings.\label{fig:dArrays}}
\end{figure}

The use of matrix operators presents the opportunity to replace the cumbersome notation of \cref{eqn:fEpsilon} with index notation and implied summation
\begin{align}
\label{eqn:matrixNotation}
	f_{\varepsilon}(\vec{x}) &= \mathbbb{d}_{m_{1} \ldots \ m_{N}} \ \mathbf{\Psi}(\vec{x})_{m_{1} \ldots \ m_{N}},
\end{align}
where $\mathbbb{d}$ is an $N$-dimensional array of wavelet coefficients and $\mathbf{\Psi}(\vec{x})$ is an $N$-dimensional array of wavelet basis functions. \Cref{fig:dArrays} illustrates $\mathbbb{d}$ in $2$D and $3$D settings and shows how the notation has been simplified in \cref{eqn:matrixNotation} as these arrays have a unique index $\vec{m}$ for every $j$, $\lambda$, and $\vec{k}$ through the mapping
\begin{align}
\label{eqn:indexMap}
    m_{i}(j, \lambda, \vec{k}) = k_{i} + \mathrm{bin}(\lambda)_{i} \times (p \ 2^{j} + 1).
\end{align}

With $\vec{x} \in \mathbb{R}^{3}$, the $3$D representation of $f(\vec{x})$ given by \cref{eqn:coeff_def} and its inverse become
\begin{align}
\label{eqn:FWT_BWT}
    \mathbb{d}_{m_{1} m_{2} m_{3}} = \mathbb{F}_{m_{1} n_{1}} \ \mathbb{F}_{m_{2} n_{2}} \ \mathbb{F}_{m_{3} n_{3}} \ f_{n_{1} n_{2} n_{3}} \qquad \mathrm{and} \qquad
	 f_{m_{1} m_{2} m_{3}} = \mathbb{B}_{m_{1} n_{1}} \ \mathbb{B}_{m_{2} n_{2}} \ \mathbb{B}_{m_{3} n_{3}} \ \mathbb{d}_{n_{1} n_{2} n_{3}},
\end{align}
on each level $j$, where $\boldsymbol{f}$ is a $3$D array of field values. The application of the $\mathbbb{F}$ operators is sometimes called a forward wavelet transform (FWT) or wavelet analysis, whereas application of the $\mathbbb{B}$ operators is sometimes called a backward wavelet transform (BWT) or wavelet synthesis \cite{Goedecker, Farge1992}. These matrices are sparse, banded, and constant in time. Due to these properties, the $\mathbbb{F}$ and $\mathbbb{B}$ matrices are never fully assembled, and only nonredundant nonzero entries are stored in memory. Therefore, the FWT and BWT operations have a matrix-free computational implementation.

\subsection{Wavelet derivatives}
\label{sec:derivatives}
Additionally, the evaluation of spatial derivatives must be defined in order to approximate the solution of PDEs. The Deslauriers-Dubuc wavelet family is continuous and differentiable, with the smoothness of the basis analyzed in \cite{Rioul1992} and summarized in \cite{Harnish2018,Harnish2021}. This continuity allows the spatial derivative operator to act directly on the basis functions
\begin{align}
\label{eqn:D_basis}
    \frac{\partial^{\alpha}}{\partial x_{i}^{\alpha}} f(\vec{x}) &\approx \frac{\partial^{\alpha}}{\partial x_{i}^{\alpha}}\bigg( \mathbbb{d}_{m_{1}  \ldots \  m_{N}} \ \mathbf{\Psi}(\vec{x})_{m_{1}  \ldots \  m_{N}} \bigg) = \mathbbb{d}_{m_{1}  \ldots \  m_{N}} \frac{\partial^{\alpha} \mathbf{\Psi}(\vec{x})_{m_{1}  \ldots \ m_{N}}}{\partial x_{i}^{\alpha}}.
\end{align}
As in \cite{Beylkin1997}, we project the spatial derivative of the basis back onto the wavelet basis and this combination of differentiation and projection transforms \cref{eqn:D_basis} into
\begin{align} 
\label{eqn:D_product}
    \frac{\partial^{\alpha}}{\partial x_{i}^{\alpha}} f(\vec{x}) &\approx {}^{(\alpha, i)}\mathbb{D}_{m_{1}  \ldots \  m_{N}}^{n_{1} \ldots \  n_{N}} \ \mathbbb{d}_{n_{1} \ldots \  n_{N}} \ \mathbf{\Psi}(\vec{x})_{m_{1} \ldots \  m_{N}},
\end{align}
where the operator ${}^{(\alpha, i)}\mathbbb{D}$ is defined in terms of ``connection coefficients'' that result from evaluating the integral
\begin{align}
\label{eqn:Gammas}
    {}^{(\alpha, i)}\mathbb{D}_{\lambda, j, \vec{k}}^{\beta, r, \vec{l}} &= \int \frac{\partial^{\alpha}}{\partial x_{i}^{\alpha}} \bigg( {}^{\beta}\Psi_{\vec{l}}^{r}(\vec{x}) \bigg) \ {}^{\lambda}\widetilde{\Psi}_{\vec{k}}^{j}(\vec{x}) \ \mathrm{d}\vec{x},
\end{align}
for all possible combinations of $\beta, \lambda, r, j \in \mathbb{Z}$ present in the multidimensional wavelet discretization. Much like the integrals in \cref{eqn:coeff_def}, leveraging the properties of wavelet basis functions allows the integrals in \cref{eqn:Gammas} to be solved exactly and expressed in terms of eigenvector solutions and the four types of filter coefficients (\emph{i.e.}, $h, \ \widetilde{h}, \ g,$ and $\widetilde{g}$). Application of the ${}^{(\alpha, i)}\mathbbb{D}$ operator results in a discrete approximation of the $\alpha^{\mathrm{th}}$ order derivative in the $i$-direction with the spatial error
\begin{align}
\label{eqn:D_error}
    \bigg| \bigg| \frac{\partial^{\alpha} f }{\partial x_{i}^{\alpha}}  &- \frac{\partial^{\alpha} f_{j} }{\partial x_{i}^{\alpha}} \bigg| \bigg|_{\infty} \leq \mathcal{O} \left( \varepsilon^{1-\frac{\alpha}{p}} \right).
\end{align}
Note that \cref{eqn:D_error} reflects the error of the dense $f_{j}$ discretization and not the sparse $f_{\varepsilon}$ discretization since the derivation from \cite{Harnish2018} assumes a dense grid up to level $j_{\mathrm{max}}$. Details on how to apply the ${}^{(\alpha, i)}\mathbbb{D}$ operator can be found in \ref{sec:fullDexpression}. Much like the $\mathbbb{F}$ and $\mathbbb{B}$ operators, the $\mathbbb{D}$ operator is also sparse and constant in time. Therefore it too has a mostly matrix-free computational implementation which leverages recurrence patterns to compress the memory requirements.

\subsection{Nonlinear operations}
\label{sec:nonlinear}
In order to numerically solve nonlinear PDEs, the discrete analogs of nonlinear operations must be well defined. This would be computationally expensive to perform in the wavelet representation because it requires a convolution operation. Instead, inspired by pseudo-spectral solvers, this work performs nonlinear operations point-wise in the physical domain. Specifically, for each nonlinear operation, a BWT maps the wavelet coefficients to their corresponding field values. Then, the nonlinear operation is approximated by acting on field values at each collocation point. Finally, a FWT of the result returns each collocation point to the wavelet domain. While this approach is more computationally efficient, it introduces aliasing errors. An estimate of the magnitude of such errors is provided in \cite{Holmstrom1999}, where it is shown that this process approximates the product of fields, $f(\vec{x})$ and $g(\vec{x})$, with the spatial error
\begin{align}
\label{eqn:product_error}
    ||f(\vec{x}) \times g(\vec{x}) - f_{\varepsilon}(\vec{x}) \times g_{\varepsilon}(\vec{x})||_{\infty} \leq \mathcal{O}(\varepsilon).
\end{align}
\section{Multiresolution wavelet toolkit}
\label{sec:MRWT}
This section outlines the multiresolution wavelet toolkit (MRWT), an algorithm designed to numerically approximate the solution to systems of nonlinear PDEs using a sparse discretization with adaptivity and error control. To illustrate the process, consider a general initial-boundary value problem on a finite domain $\Omega \subset \mathbb{R}^{3}$
\begin{align} 
\label{eqn:general_PDE}
    \frac{\partial \vec{u}}{\partial t} &= \mathcal{F}\left(\vec{x}, t, \vec{u}, \frac{\partial^{\alpha} \vec{u}}{\partial x_{i}^{\alpha}} \right)
    \ 
    \mathrm{in} \ \Omega \times [0,T],
    \nonumber\\
    \vec{u} &= \{u_{1}^{0}, u_{2}^{0}, u_{3}^{0}\} 
    \
    \mathrm{in} \ \Omega \times (t = 0),
    \nonumber\\
    \vec{u} &= \vec{u}_{d} 
    \
    \mathrm{on} \ \partial \Omega_{d} \times [0,T],
    \nonumber\\
    \frac{\partial^{\alpha} \vec{u}}{\partial x_{i}^{\alpha}} &= \vec{u}_{n} 
    \
    \mathrm{on} \ \partial \Omega_{n} \times [0,T].
\end{align}
The order of the PDE sets requirements on the continuity of the wavelet basis, influencing the wavelet order $p$. Once $p$ is set, fields and their spatial derivatives are projected onto the wavelet basis functions using \cref{eqn:matrixNotation} and \cref{eqn:D_product}. During this initialization phase, MRWT precomputes the nonredundant entries in the $\mathbbb{F}$, $\mathbbb{B}$, and ${}^{(\alpha, i)}\mathbbb{D}$ operators. This process transforms the nonlinear PDE into a nonlinear ODE
\begin{align}
\label{eqn:ode}
\frac{\mathrm{d} \vec{u}_{\varepsilon}}{\mathrm{d}t} = \mathcal{F}_{\varepsilon} \left(\vec{x}, t, \vec{u}_{\varepsilon} \right).
\end{align}

Each scalar component of the initial condition (\emph{i.e.}, $u_{1}^{0}$, $u_{2}^{0}$, and $u_{3}^{0}$) is projected individually onto the wavelet basis one resolution level at a time, starting from the coarsest resolution. All wavelet coefficients on $j = 1$ are calculated and each point of type $\lambda = 0$ or with a magnitude greater than or equal to a user-prescribed threshold $\varepsilon$ is deemed essential. The set of essential points for a particular scalar (\emph{e.g.}, $u_{1}$) is defined by 
\begin{align}
\label{eqn:essential}
    \mathbb{E}_{\bullet} = \left\{ \vec{x} \in \Omega : \left[ \left|{}^{0}\mathbb{d}_{\vec{k}}^{1}(\bullet)\right| \ge 0 \; \vee \; \left|{}^{\lambda}\mathbb{d}_{\vec{k}}^{j}(\bullet)\right| \geq \varepsilon \right], \; j \in [1, j_{\mathrm{max}}], \; \lambda \in [1, 7], \; \vec{k} \in \mathbb{Z} \right\},
\end{align}
where $\bullet$ is a marker for a particular scalar field. The set $\mathbb{E}_{\bullet}$ is populated on higher resolution levels through exploring spatial regions near the already discovered essential points. Identifying and calculating all essential points for a particular scalar creates a sparse multiresolution spatial discretization without generating dense grids of uniform spacing. This process is repeated for each scalar component of the initial condition (\emph{i.e.}, $\mathbb{E}_{u_{1}}$, $\mathbb{E}_{u_{2}}$, and $\mathbb{E}_{u_{3}}$). Furthermore, according to \cref{eqn:ferror}, the spatial accuracy of the initial condition is known \emph{a priori} to be $\mathcal{O} \left(\varepsilon\right)$.

Since the solution of the PDE may advect and evolve over time, the points which are presently essential $\mathbb{E}^{n}$ (where $n$ is the current time step) may not be same collection of points which are essential in future times (\emph{e.g.}, $\mathbb{E}^{n + 1}$). Therefore, MRWT uses a predictor-corrector strategy to predict where new collocation points will be required and corrects that prediction during temporal integration until the necessary grid is obtained at each time step. A trial grid is defined for each scalar by including points in the neighborhood of that scalar's essential points such that
\begin{align}
\label{eqn:trial}
    \mathbb{T}_{\bullet} &= \bigcup_{\forall \vec{y} \in \mathbb{E}_{\bullet}} \left\{ \vec{x} \in \Omega : \left| y_{i} - x_{i} \right| < \frac{b_{i} - a_{i}}{p \ 2^{r}}, \; i \in [1, 3], \; j \in [r, \min(r + 1,j_{\mathrm{max}})] \right\},
\end{align}
where the level $r$ is obtained from the coordinate of the essential point $\vec{y} \in \mathbb{E}_{\bullet}$ using the inverse mapping from \cref{eqn:invMap}.

Next, the ${}^{(\alpha, i)}\mathbbb{D}$ operator is used to identify which points influence the derivative calculations at each point in the trial grid $\mathbb{T}$. As shown in \ref{sec:fullDexpression}, the wavelet derivative operator is a multiresolution weighted sum. Calculating a wavelet coefficient corresponding to a spatial derivative on resolution level $r$ will involve contractions with points below, above, and on level $r$. However, because the grid is sparse, many of the collocation points involved in the application of the ${}^{(\alpha, i)}\mathbbb{D}$ operator, given by \cref{eqn:Gammas} may be absent. This affects the validity of the error estimate in \cref{eqn:D_error} as its derivation assumes complete stencils on each resolution level with $f_{j}(\vec{x})$ given by \cref{eqn:fJ}. Since attempting to complete stencils on all levels $j \in [1, j_{\mathrm{max}}]$ for all points would lead to dense grids of uniform spacing (\emph{i.e.}, $f_{j}(\vec{x})$ in \cref{eqn:fJ}), the MRWT algorithm completes derivative stencils for each point in each scalar's trial grid up to its own level $r$ according to
\begin{align}
\label{eqn:grid}
    \mathbb{G}_{\bullet} &= \bigcup_{\forall \vec{y} \in \mathbb{T}_{\bullet}} \left\{ \vec{x} \in \Omega : \left[ \vec{x} = \vec{y} \; \vee \; \left| {}^{\lambda}_{\beta} \mathbb{D}^{j, \vec{k}}_{r, \vec{l}} \right| > 0 \right], \; j \leq r, \; \lambda \in [0, 7], \; \vec{k} \in \mathbb{Z} \right\},
\end{align}
where the type $\beta$, level $r$, and index $\vec{l}$ are obtained from the coordinate of the point $\vec{y} \in \mathbb{T}_{\bullet}$ using the inverse mapping from \cref{eqn:invMap}. 

For each scalar in the system of PDEs, the algorithm constructs $\mathbb{E}_{\bullet} \subset \mathbb{T}_{\bullet} \subset \mathbb{G}_{\bullet}$ and the sparse multiresolution computational grid defined by $\mathbb{G} = \mathbb{G}_{u_{1}} \cup \mathbb{G}_{u_{2}} \cup \mathbb{G}_{u_{3}}$. The grid, $\mathbb{G}_{\bullet}$ contains all of the collocation points that are required to resolve that particular field and its derivatives at discrete times $n$ and $n+1$, with the spatial accuracy bound by \cref{eqn:ferror,eqn:product_error} and approximated by \cref{eqn:D_error}. However, due to overlap between the sets from different scalars, the computational grid may contain more data than is needed for a given scalar. This is corrected by zeroing extraneous coefficient values
\begin{align}
\label{eqn:decouple}
    {}^{\lambda}\mathbb{d}_{\vec{k}}^{j}(\bullet) = 0 \; \forall \ \vec{x} \in \mathbb{G} \setminus \mathbb{G}_{\bullet},
\end{align}
where the type $\lambda$, level $j$, and index $\vec{k}$ are obtained from the coordinate of the point $\vec{x}$ using the inverse mapping from \cref{eqn:invMap}. 
\Cref{eqn:essential,eqn:trial,eqn:grid,eqn:decouple} completely define the spatial discretization used to solve systems of nonlinear PDEs. Leveraging wavelet theory from \cref{sec:wavelets} to construct the sets $\mathbb{E}, \mathbb{T},$ and $\mathbb{G}$ using the wavelet derivative operator $\mathbb{D}$ is one of the novelties of our work. As we show later, this grid construction is necessary to control the spatial accuracy.

An explicit, embedded, Runge-Kutta time integration scheme \cite{RKF45} is used to progress the solution from the time step $n$ to a trial time step $(n+1)^{*}$. This converts the ODE in \cref{eqn:ode} into a system of algebraic equations which update $\vec{u}_{\varepsilon}^{ \ n}$ to the trial time step $\vec{u}_{\varepsilon}^{ \ (n + 1)^{*}}$ while providing an estimate of the temporal error and adjusting the time-step size $\Delta t$ such that the temporal error is less than the user prescribed value. In this work, we set the tolerance to be an order of magnitude less than the spatial error epsilon (\emph{i.e.}, $\mathrm{tol} = \varepsilon / 10$) and adapt the time-step size according to
\begin{align}
\label{eqn:dtNew}
    \Delta t_{\mathrm{new}} = \Delta t \left( \frac{ \mathrm{tol} }{ \| \vec{u}^{ \ n+1}_{K+1} - \vec{u}^{ \ n+1}_{K} \|_{\infty} } \right)^{\frac{1}{K+1}},
\end{align}
where $K$ is the order of the method used to advance $n \rightarrow (n+1)^{*}$. Furthermore, similar to the work done by \cite{dtPercentChange}, we prevent the time-step size from varying too abruptly by bounding the suggested $\Delta t_{\mathrm{new}}$ such that
\begin{align}
\label{eqn:dtBound}
     \Delta t - 10 \%  \Delta t \leq \Delta t_{\mathrm{new}} \leq \Delta t + 0.1 \%  \Delta t.
\end{align}
At each intermediate stage of the Runge-Kutta integration, no boundary conditions are enforced in order to preserve high order temporal convergence \cite{noIntBC}. At the final stage of the Runge-Kutta integration, once the PDE solution has evolved to the trial time $(n+1)^{*}$, boundary conditions are enforced directly (\emph{i.e.}, for Dirichlet data) or by using the penalty formulation \cite{Hesthaven1996} that will modify \cref{eqn:ode} on the boundary (\emph{i.e.}, for Dirichlet and/or Neumann data).

The magnitudes of the wavelet coefficients at the new time will determine if the grid prediction must be corrected. A list of corrections is defined as the set difference between the essential points at the trial time step and the essential points at the current time
\begin{align}
\label{eqn:missing}
    \mathbb{C}^{n} &= \mathbb{E}^{(n + 1)^{*}} \ \setminus \ \mathbb{E}^{n},
\end{align}
If the list of corrections contains any points then these corrections are added to the list of existing essential points
\begin{align}
\label{eqn:addMissing}
    \mathbb{E}^{n} = \mathbb{E}^{n} \ \cup \ \mathbb{C}^{n}.
\end{align}
In this case, the trial time step is discarded and the computational grid at time step $n$ is supplemented with new collocation points to obtain new $\mathbb{T}^{n}_{\bullet}$ and new $\mathbb{G}^{n}_{\bullet}$ according to \cref{eqn:trial,eqn:grid}. The value of the wavelet coefficients stored at these new points is set to zero at time step $n$. After correcting the grid, a new trial time step is performed and \cref{eqn:missing,eqn:addMissing} are re-evaluated. This process of predicting and correcting is repeated iteratively until the list of corrections is an empty set. This grows the sparse computational grid as it evolves with the solution of the PDE such that the spatial error is bounded by \cref{eqn:ferror,eqn:product_error} and approximated by \cref{eqn:D_error} at each time step. 

When the trial time step is accepted (\emph{i.e.}, $\mathbbb{d}^{n + 1} = \mathbbb{d}^{(n + 1)^{*}}$) some wavelet coefficients are no longer needed. Collocation points at the new time are retained only if they are deemed essential (\emph{i.e.}, \cref{eqn:essential}), part of the trial grid (\emph{i.e.}, \cref{eqn:trial}), or needed to complete the derivative stencil (\emph{i.e.}, \cref{eqn:grid}). This procedure prunes the sparse computational grid as it evolves with the solution of the PDE.

MRWT has been implemented as a hybrid MPI+OpenMP framework using modern C++. The multidimensional sparse geometry is stored in sorted coordinate-object format (COO) and is replicated across ranks. Corresponding collocated scalar field data is distributed linearly across ranks and stored in struct-of-array form. This design permits temporal and spatial locality for sparse stencil contractions as well as vectorizable data parallelism for point-wise, right-hand-side kernels. The predictor-corrector algorithm performs bulk insertion and deletion of points relatively infrequently resulting in periodic merge and sort operations. Communication includes both point-to-point messaging and global collectives, and rebalancing is performed after each time step in which the geometry changes. The $\mathbbb{F}$ and $\mathbbb{B}$ stencil operators are matrix-free while ${}^{(\alpha, i)} \mathbbb{D}$ is optimally compressed (\emph{i.e.}, the operator is precomputed, stored sparsely, and its evaluation is carefully scheduled to amortize costs). Preliminary scaling results for the framework appear in \ref{sec:MRWTscaling}. 

\subsection{Code verification}
\label{sec:verification}
Suppose the generic initial-boundary value problem in \cref{eqn:general_PDE}, is specified to be the nonlinear quasi-$2$D Burgers' equation
\begin{align}
\label{eqn:burg2D}
    \frac{\partial u}{\partial t} + (u + c) \left[ \frac{\partial u}{\partial x_{1}} \cos{(\theta)} + \frac{\partial u}{\partial x_{2}} \sin{(\theta)} \right] &= \nu \left[ \frac{\partial^{2} u}{\partial x_{1}^{2}} \cos{^{2}(\theta)} + \frac{\partial^{2} u}{\partial x_{1} \partial x_{2}} \sin{(2 \theta)} + \frac{\partial^{2} u}{\partial x_{2}^{2}} \sin{^{2}(\theta)} \right],
\end{align}
where $\theta$ is a specified constant angle of rotation, projecting the $1$D Burgers' equation into a second spatial dimension. Two analytical solutions exist distinguished by their initial conditions and the value of $c$. In both cases, the boundary conditions are Dirichlet defined by the respective analytical solutions with the following parameters
\begin{align}
    \vec{x} &\in [-1, 1]^{2}, &
    t &\in [0, 1/2], &
    \theta &\in [0, \pi/2], &
    \nu &= 10^{-2}.
\end{align}

First, for $c=0$
\begin{align}
\label{eqn:burgRefining}
    u_{0} &= -\sin{(x_{1} \pi \cos{\theta} + x_{2} \pi \sin{\theta})},
    \nonumber\\
    u(\vec{x},t) &= - \frac{\int \sin(\pi x_{1} \cos(\theta) + \pi x_{2} \sin(\theta) - \pi \eta) \ \exp\left( \frac{-\cos(\pi x_{1} \cos(\theta) + \pi x_{2} \sin(\theta) - \pi \eta)}{2\pi \nu}\right) \ \exp\left( \frac{-\eta^{2}}{4 \nu t} \right) \mathrm{d}\eta}{\int \exp\left( \frac{-\cos(\pi x_{1} \cos(\theta) + \pi x_{2} \sin(\theta) - \pi \eta)}{2\pi \nu}\right) \ \exp\left( \frac{-\eta^{2}}{4 \nu t} \right) \mathrm{d}\eta}.
\end{align}
The MRWT solution of \cref{eqn:burgRefining} can be obtained at early times with a coarse grid, but requires increasingly higher resolution as the gradients become more extreme over time. This formulation verifies MRWT's ability to dynamically refine the sparse multiresolution spatial discretization. 

Alternatively, for $c\neq0$
\begin{align}
\label{eqn:burgWalking}
    u_{0} &= -\tanh{[((x_{1} + 0.5) \cos{\theta} + (x_{2} + 0.5) \sin{\theta})/ 2 \nu]},
    \nonumber\\
    u(\vec{x}, t) &= -\tanh{\left( \frac{(x_{1} + 0.5) \cos{(\theta)} + (x_{2} + 0.5) \sin{(\theta)} - c t}{2 \nu} \right)}.
\end{align}
The MRWT solution of \cref{eqn:burgWalking} requires high resolution only in the vicinity of the steep feature that advects through the spatial domain over time. This formulation verifies MRWT's ability to insert and remove collocation points to resolve evolving features (\emph{e.g.}, shocks).

\begin{figure}[!htb]
\begin{tabular*}{\textwidth}{@{} c @{\extracolsep{\fill}} c @{}}
    \includegraphics[width=0.5\textwidth]{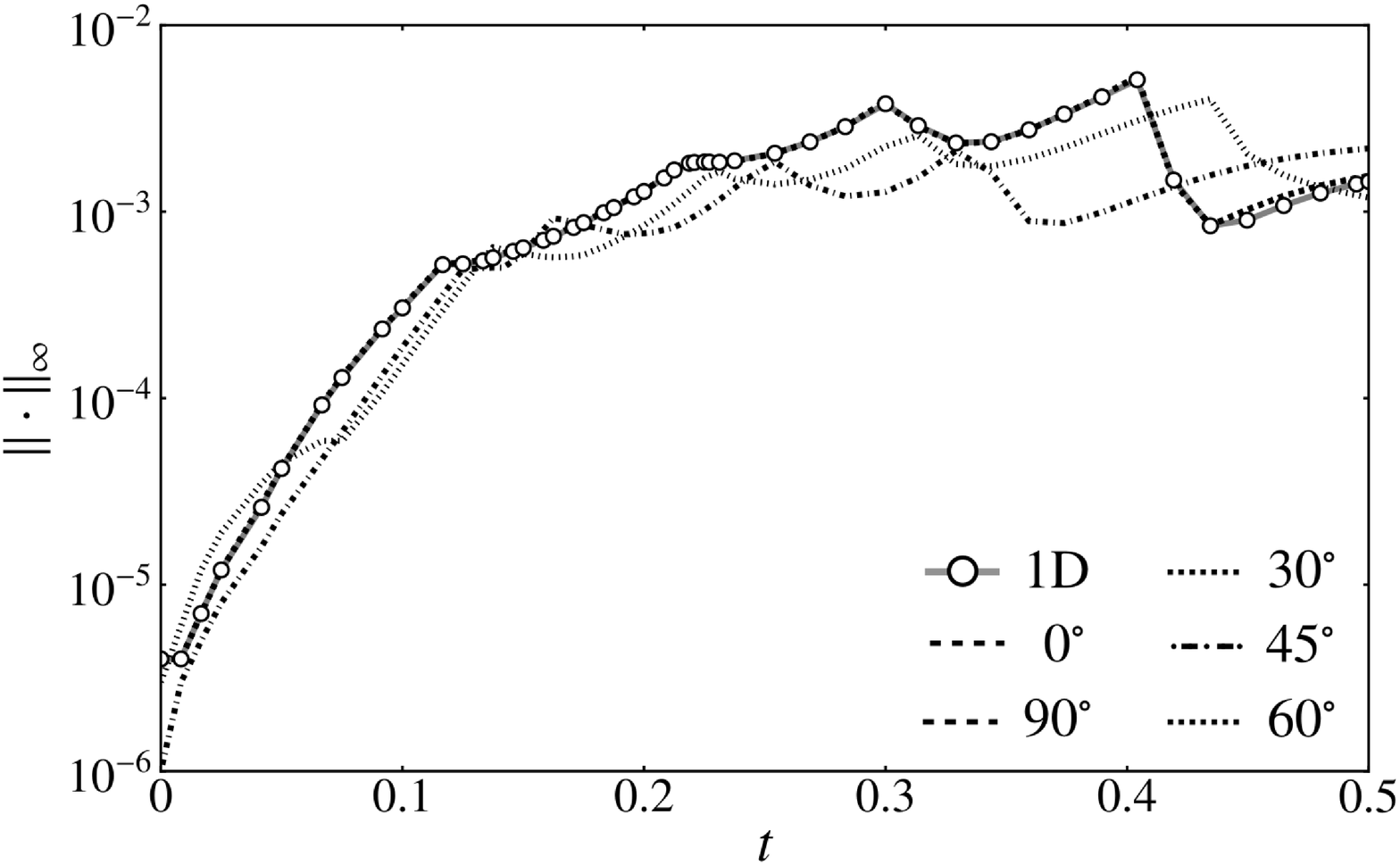}&
    \includegraphics[width=0.5\textwidth]{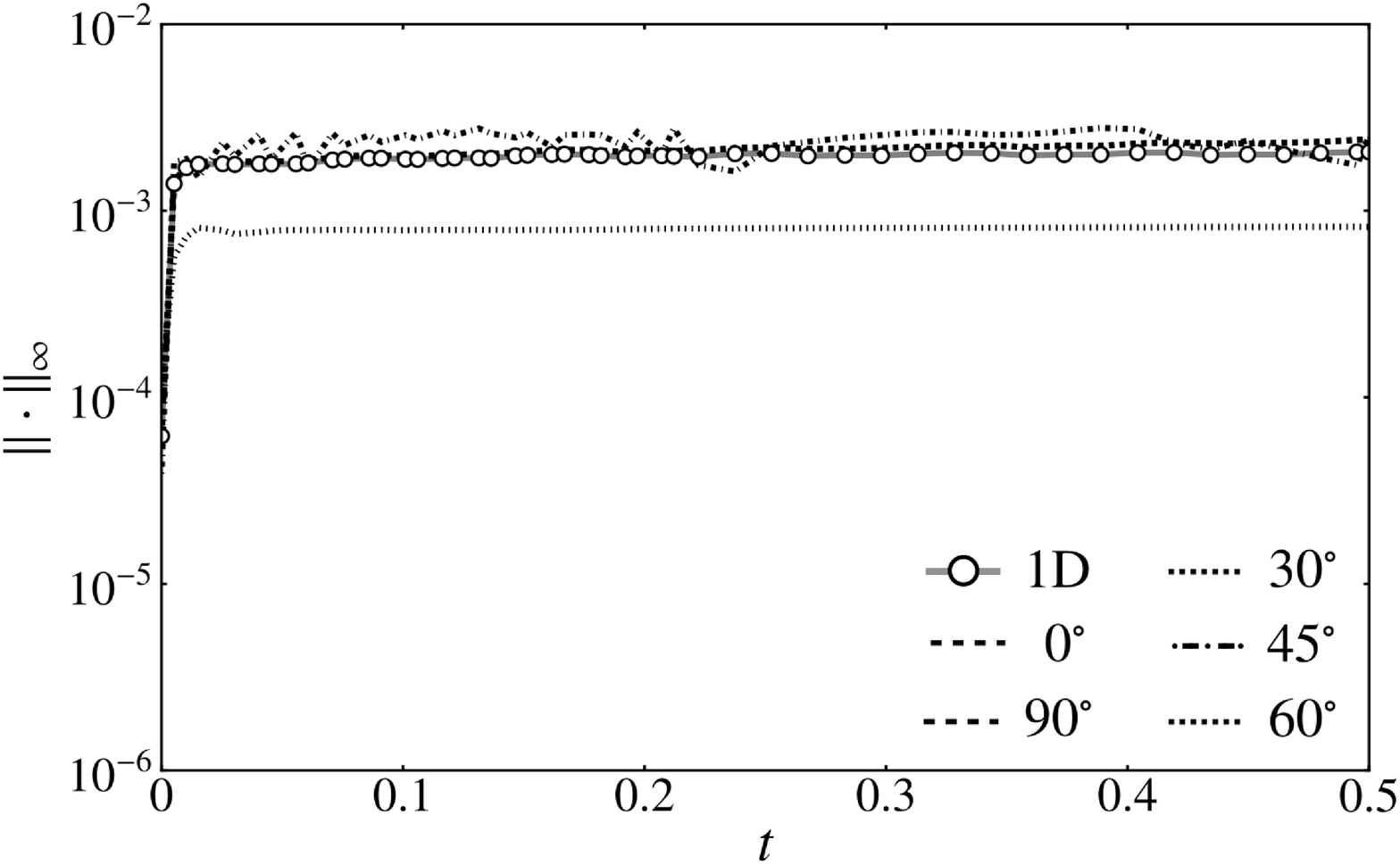}\\
    (a) & (b)
\end{tabular*}
\caption{Spatial accuracy over time for the $2$D Burgers' equation with $p = 6$, $\varepsilon = 10^{-2}$, and a variety of $\theta$ angles using (a) $c=0$ and (b) $c=2$.\label{fig:BurgBoth_error}}
\end{figure}

The MRWT algorithm has been used with the embedded $\mathcal{O}(\Delta t^{4})$ and $\mathcal{O}(\Delta t^{5})$ explicit Runge-Kutta method developed in \cite{RKF45} to solve \cref{eqn:burg2D} for a variety of angles $\theta$. The numerical solution has been compared to the analytical solutions from \cref{eqn:burgRefining,eqn:burgWalking} to evaluate the error at each time step. \Cref{fig:BurgBoth_error} shows that the error in the solution is bounded regardless of the angle $\theta$ and the error is comparable to the one-dimensional analog of \cref{eqn:burg2D}. Note that there is no distinction in the measured errors between the $0^{\circ}$ and $90^{\circ}$ or the $30^{\circ}$ and $60^{\circ}$ rotations, since these solutions are reflections across the $x_{1} = x_{2}$ plane.

\begin{figure}[!htb]
\begin{tabular*}{\textwidth}{@{\extracolsep{\fill}} c @{} c @{}}
    \includegraphics[width=0.4\textwidth]{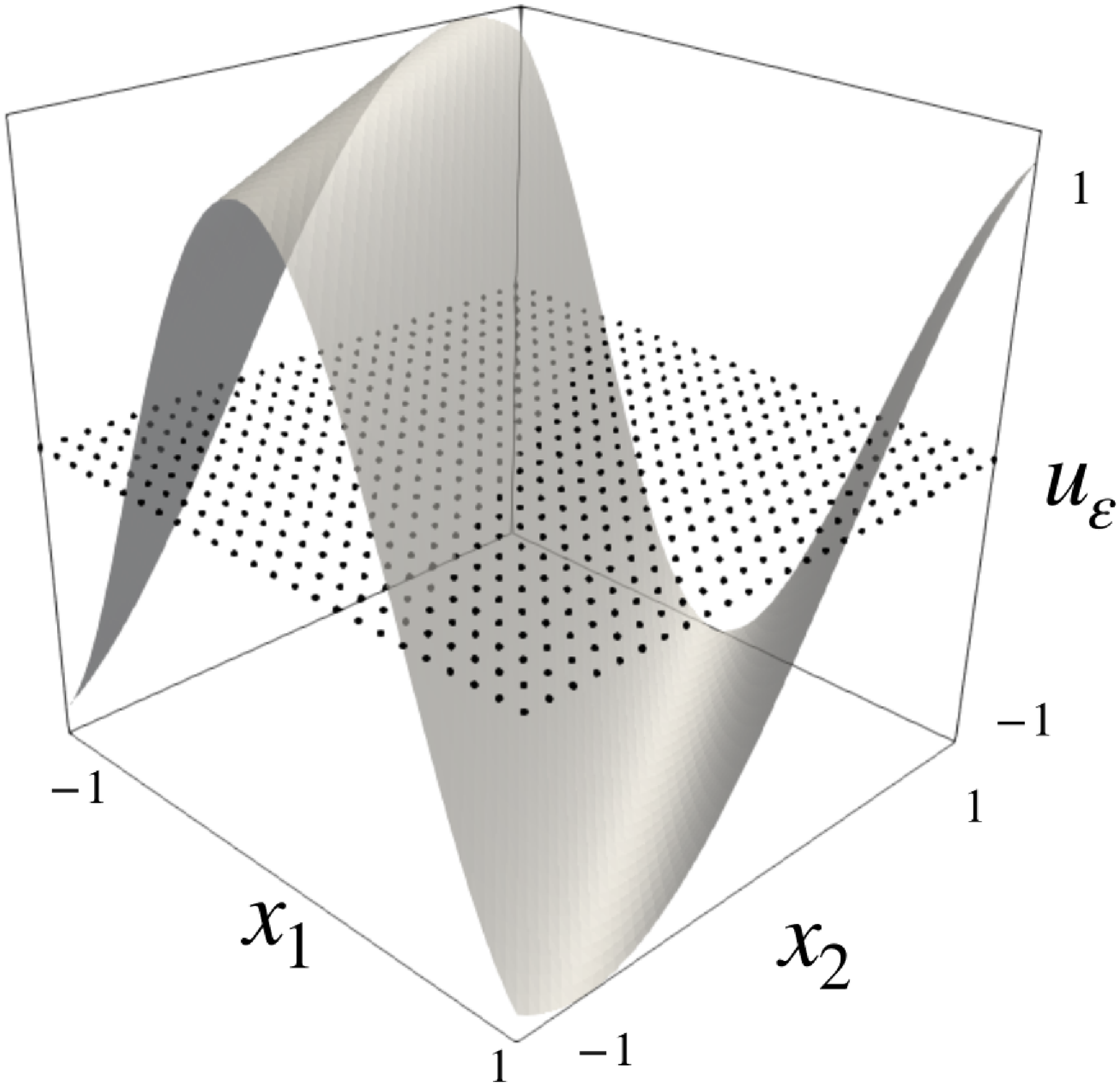}&
    \includegraphics[width=0.4\textwidth]{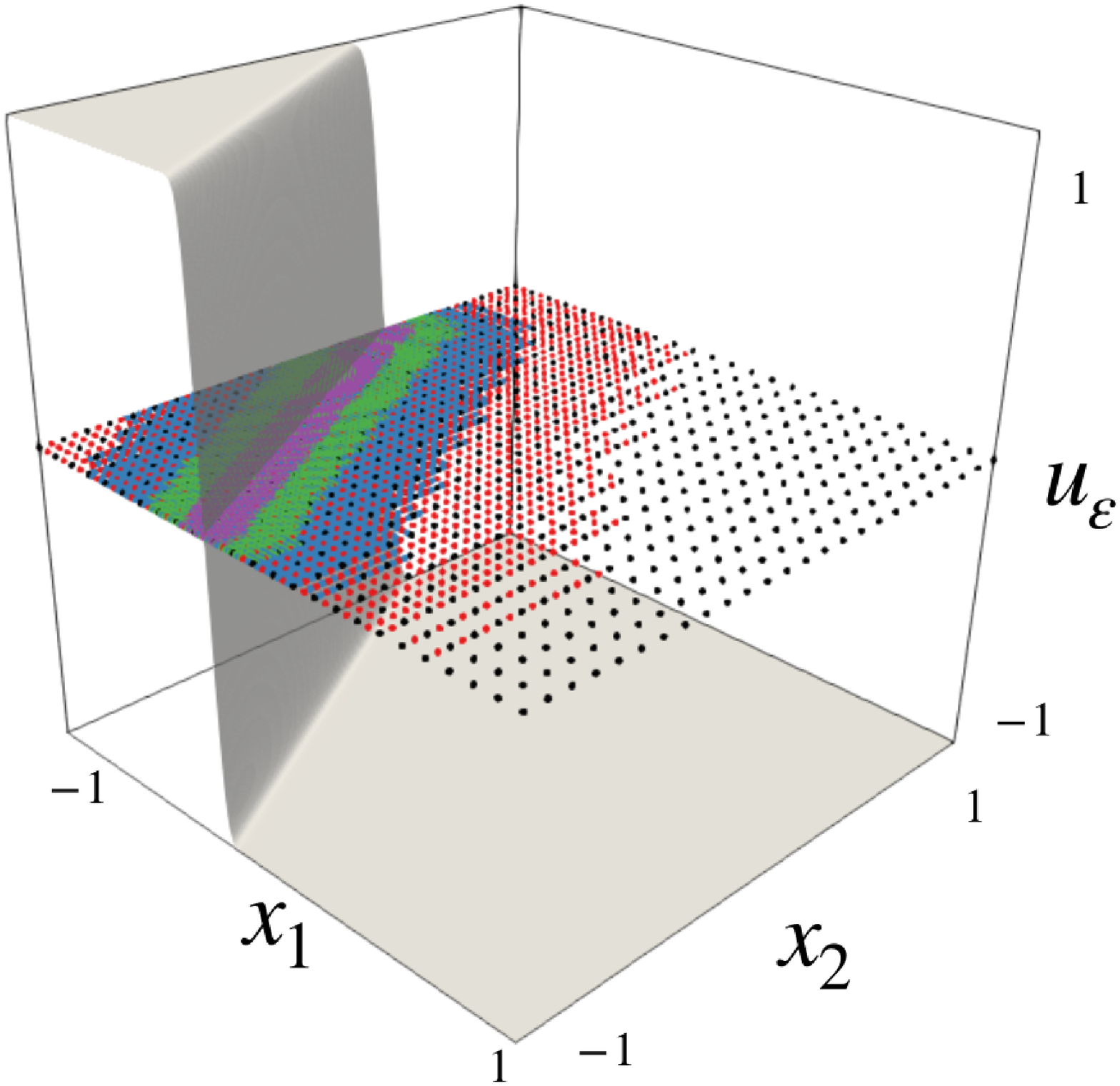}\\
    (a) & (b)\\
    \includegraphics[width=0.4\textwidth]{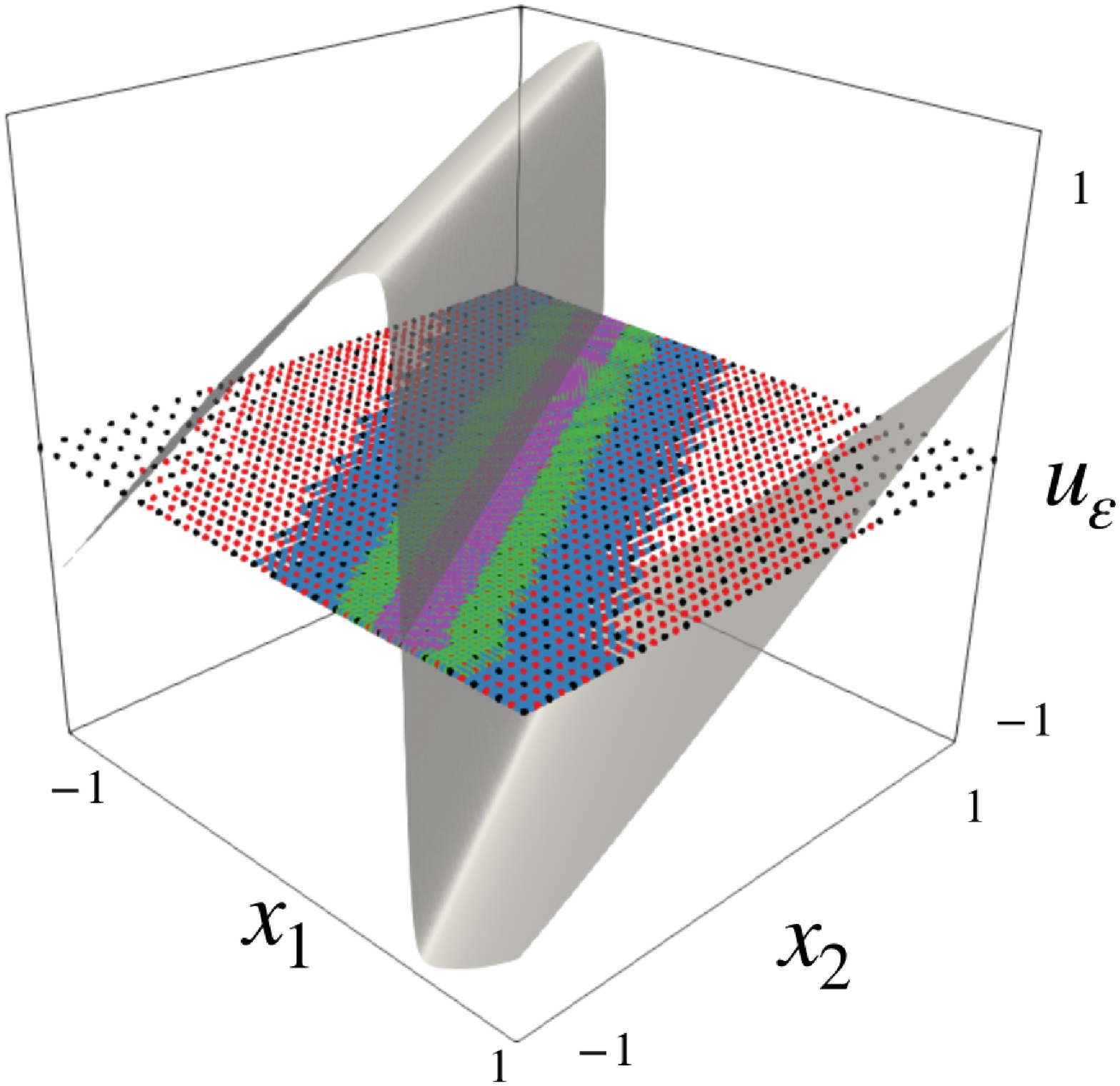}&
    \includegraphics[width=0.4\textwidth]{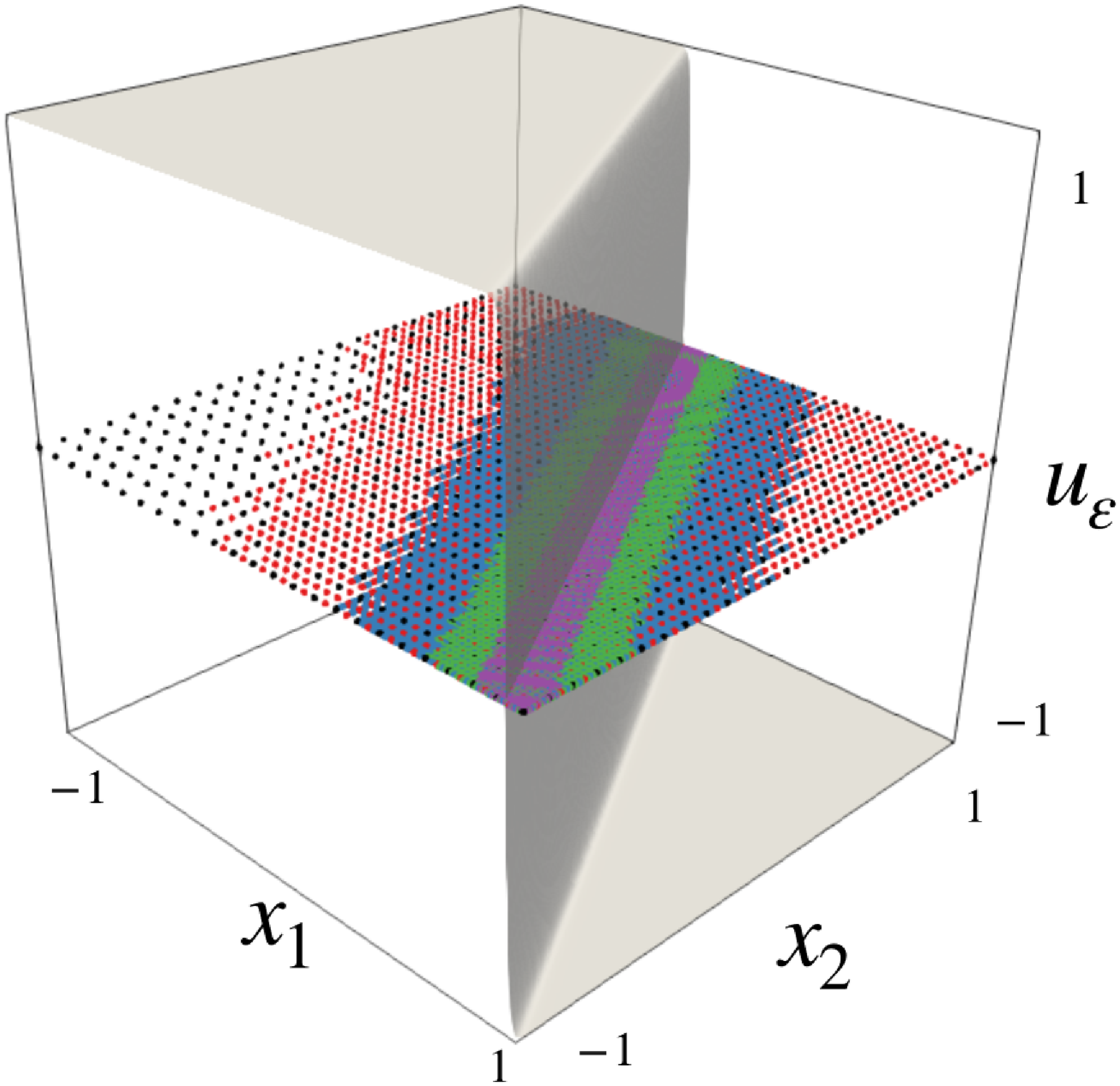}\\
    (c) & (d)
\end{tabular*}
\caption{MRWT solution and sparse multiresolution grid for the $2$D Burgers' equation with $c = 0$ on the left and $c=2$ on the right. Additionally, $\theta = 30^{\circ}$, $p=6$, and $\varepsilon = 10^{-2}$ at: (a) \& (b) $t=0$ s and (c) \& (d) $t=\frac{1}{2}$ s. The grid points are colored according to their resolution level $j$. The reader is referred to the online version of this article for clarity regarding the color in this figure.\label{fig:BurgBoth_MRWT}}
\end{figure}

Additionally, the solution to both versions of \cref{eqn:burg2D} with $\theta = 30^{\circ}$, $p=6$, and $\varepsilon = 10^{-2}$ is shown at various times in \cref{fig:BurgBoth_MRWT} demonstrating that the sparse dynamically adaptive grid provides resolution only where and when it is needed to track features as they develop or move through the domain.

\clearpage

\begin{figure}[!htb]
\begin{tabular*}{\textwidth}{@{} c @{\extracolsep{\fill}} c @{}}
    \includegraphics[width=0.45\textwidth]{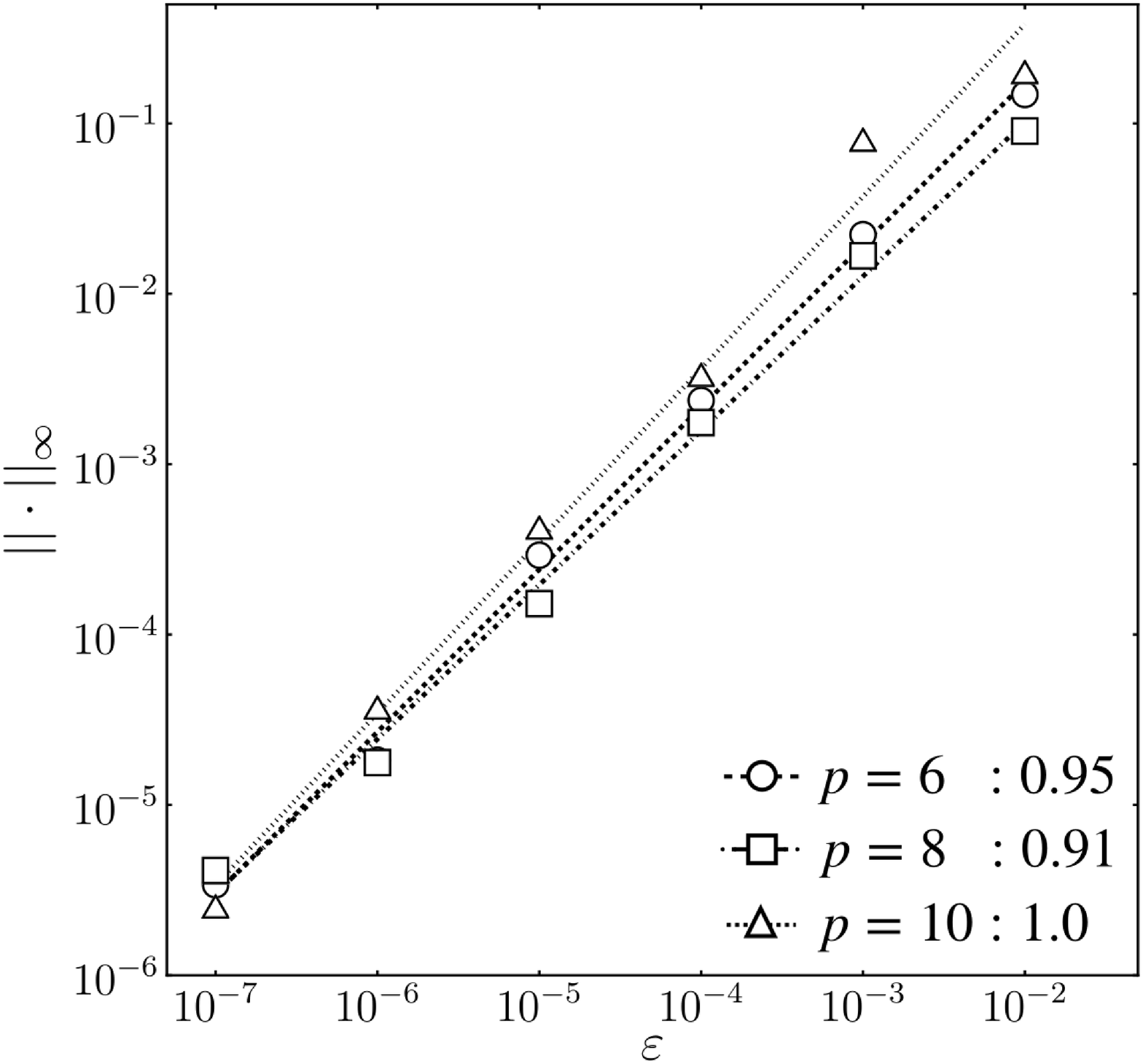}&
    \includegraphics[width=0.465\textwidth]{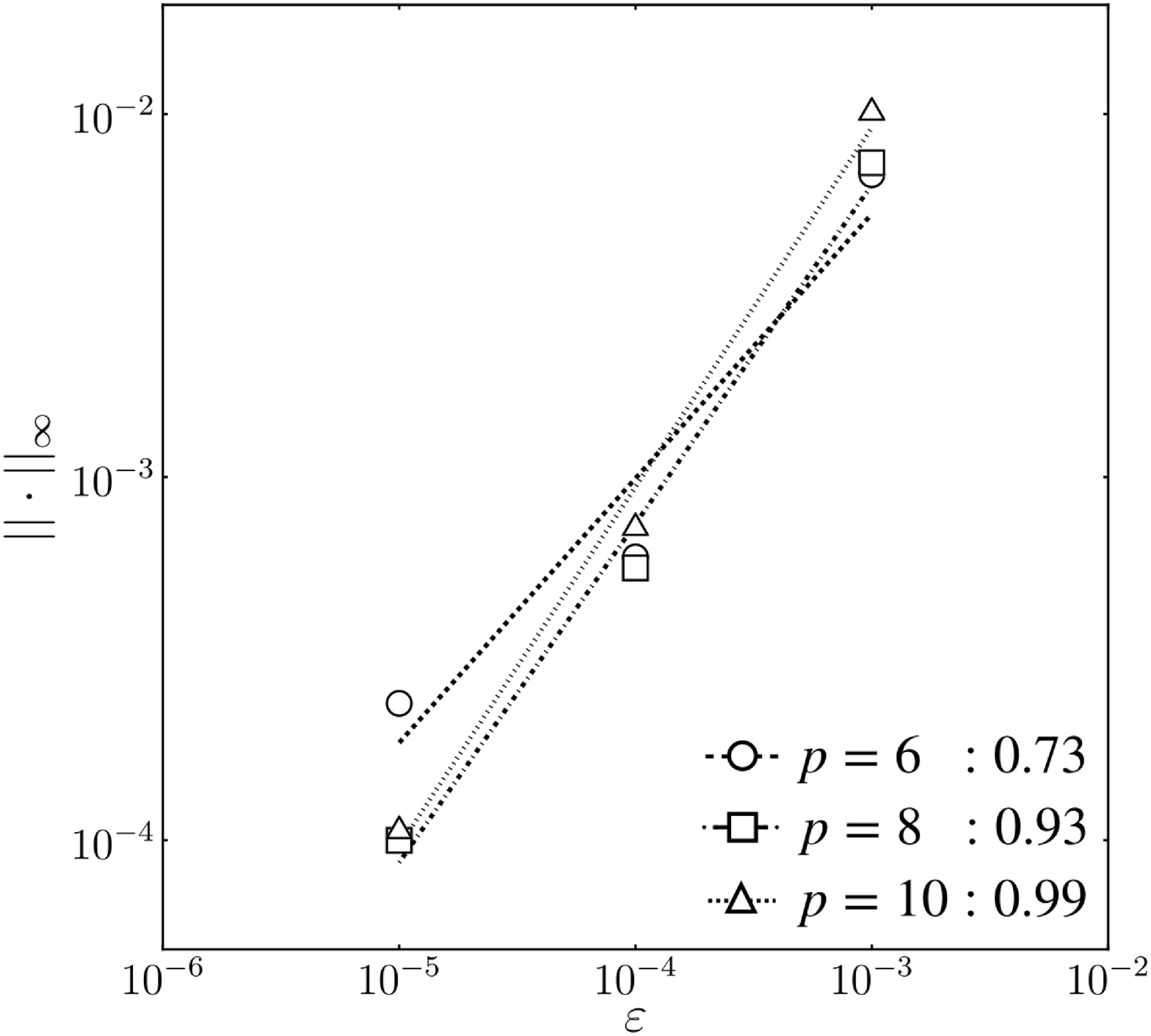}\\
    (a) & (b)
\end{tabular*}
\caption{Spatial convergence rates for the $2$D Burgers' equation with $\theta = 30^{\circ}$ using a variety of $p$ and $\varepsilon$ with (a) $c=0$ and (b) $c=2$. \label{fig:BurgBoth_conv}}
\end{figure}

Furthermore, the analytical solutions of \cref{eqn:burgRefining,eqn:burgWalking} have been used to evaluate the error and corresponding convergence rate of the algorithm by solving \cref{eqn:burg2D} for a range of $\varepsilon$ and $p$ values. \Cref{fig:BurgBoth_conv} shows the convergence rates calculated half-way through the simulation time with $\theta = 30^{\circ}$, which are in agreement with the theoretical estimates from \cref{eqn:ferror,eqn:product_error,eqn:D_error}. In particular, the order of error for the nonlinear PDE should be between $\mathcal{O}(\varepsilon^{1-\frac{\alpha}{p}})$ and $\mathcal{O}(\varepsilon)$. Results presented in \cref{fig:BurgBoth_conv} retain these properties.

In order to test MRWT's ability to solve systems of equations, canonical Navier-Stokes examples such as the Sod shock tube and a Taylor-Sedov blast wave have been presented in \cite{Harnish2018,Harnish2021,HarnishSpecialIssue}. \ref{sec:sedov} contains analysis of $1$D, $2$D, and $3$D simulations of the Taylor-Sedov blast wave. Additional verification of the MRWT algorithm is provided by comparing the location of the simulated shock wave to the theoretical prediction from \cite{sedov}.
\clearpage

\section{Nonlinear dynamic damage modeling}
\label{sec:model}
Having verified the mathematical correctness of the MRWT algorithm, we now demonstrate the multiscale capabilities by modeling high-strain rate damage nucleation and propagation in nonlinear solids. The following subsections will outline the continuum formulation of finite strain kinematics and the constitutive description provided by damage mechanics. Since this work leverages wavelet basis functions for the computational implementation, all governing equations will be defined in the Eulerian frame.
%
%
\subsection{Governing equations}
Let a body in a reference configuration $\mathcal{B}_{0} \subset \mathbb{R}^{3}$ be defined by a continuum of material points (\emph{i.e.}, particles $\vec{X} \in \mathcal{B}_{0}$). After some time $t \in [0, T]$, the body is in a deformed configuration where $\varphi : \mathcal{B}_{0} \times [0, T] \rightarrow \mathbb{R}^{3}$ defines the motion that maps particles to their positions $\vec{x} = \varphi(\vec{X}, t)$ (\emph{i.e.}, places) in the current configuration $\mathcal{B}(t) = \varphi(\mathcal{B}_{0}, t)$. The mechanical effects are governed by the conservation of linear momentum in the Eulerian frame $\Omega \subset \mathcal{B}(t) \subset \mathbb{R}^{3}$
\begin{align}
\label{eqn:momentum}
    \rho \dot{\vec{v}} &= \nabla \cdot \boldsymbol{\sigma} + \rho \vec{b},
    \ 
    \mathrm{in} \ \Omega \times [0,T],
    \nonumber\\
    \vec{v} &= \vec{v}_{0} 
    \
    \mathrm{in} \ \Omega \times (t = 0),
    \nonumber\\
    \vec{v} &= \vec{v}_{d} 
    \
    \mathrm{on} \ \partial \Omega_{d} \times [0,T],
    \nonumber\\
    \frac{\partial^{\alpha} \vec{v}}{\partial x_{i}^{\alpha}} &= \vec{v}_{n} 
    \
    \mathrm{on} \ \partial \Omega_{n} \times [0,T],
\end{align}
where $\rho(\vec{x},t)$ is the density, $\boldsymbol{\sigma}(\vec{x},t)$ is the Cauchy stress tensor, and $\vec{b}(\vec{x},t)$ is a body force per unit mass. 

Inspired by the reference map technique \cite{RMT}, the inverse mapping $\vec{X} = \varphi^{-1}(\vec{x}, t)$ is used to determine which particles occupy specified spatial coordinates (\emph{i.e.}, places $\vec{x}$). However, instead of keeping track of the reference map $\varphi^{-1}(\vec{x}, t)$ directly, we model the displacement $\vec{u}(\vec{x},t) = \vec{x} - \vec{X}$, and like in \cite{RMT} we derive the associated evolution equation from $\dot{\vec{X}} = \vec{0}$ 
\begin{align}
\label{eqn:displacement}
    \dot{\vec{u}} &= \vec{v},
    \ 
    \mathrm{in} \ \Omega \times [0,T],
    \nonumber\\
    \vec{u} &= \vec{u}_{0} 
    \
    \mathrm{in} \ \Omega \times (t = 0),
    \nonumber\\
    \vec{u} &= \vec{u}_{d} 
    \
    \mathrm{on} \ \partial \Omega_{d} \times [0,T],
    \nonumber\\
    \frac{\partial^{\alpha} \vec{u}}{\partial x_{i}^{\alpha}} &= \vec{u}_{n} 
    \
    \mathrm{on} \ \partial \Omega_{n} \times [0,T].
\end{align}
To finalize the Lagrangian-Eulerian coupling of \cref{eqn:momentum,eqn:displacement}, we use the conservation of mass in the Lagrangian frame
\begin{align}
\label{eqn:mass}
    \rho &= \frac{\rho_{0}}{\mathrm{J}}, \mathrm{in} \ \Omega_{0} \times [0,T],
\end{align}
where $\rho_{0}(\vec{X})$ is the density in the reference configuration, $\mathrm{J}(\vec{X}, t) = 1 / \det{\mathbf{f}}$ is the Jacobian of the deformation, $\mathbf{f}(\vec{x}, t) = \boldsymbol{1} - \nabla \vec{u}$ is the inverse deformation gradient, and $\boldsymbol{1}$ is the second order identity tensor.

\clearpage
%
%
\subsection{Constitutive theory}
\label{sec:closureEquations}
The constitutive relations are derived from a thermodynamically consistent description with an energy-based isotropic damage model, applicable under finite strains. In this work, the following potentials describe the strain energy density (per reference volume) and the Helmholtz specific energy
\begin{align}
\label{eqn:potentials}
    W(\omega, \boldsymbol{b}) &= (1 - \omega) \left[ \underbrace{\frac{\kappa}{2} \left( e^{(\mathrm{J} - 1)} - \ln(\mathrm{J}) - 1 \right)}_{\text{volumetric}} + \underbrace{\frac{\mu}{2} \left( \mathrm{J}^{-\frac{2}{3}} \mathrm{tr}(\boldsymbol{b}) - 3 \right)}_{\text{deviatoric}} \right],
    &
    \psi(\omega, \boldsymbol{b}) &= \frac{1}{\rho_{0}} W(\omega, \boldsymbol{b}),
\end{align}
where $\boldsymbol{b}(\vec{x}, t)$ is the left Cauchy-Green strain tensor defined through its inverse $\boldsymbol{b}^{-1}(\vec{x}, t) = \mathbf{f}^{T} \mathbf{f}$, $\kappa$ is the bulk modulus, $\mu$ is the shear modulus, and $\omega \in [0, 1]$ is the scalar damage variable.

Neglecting all but mechanical effects, the second Law of Thermodynamics has the form
\begin{align}
\label{eqn:2nd_law}
    \rho \dot{\psi} - \boldsymbol{\sigma} : \nabla \vec{v} \leq 0.
\end{align}
Combining \cref{eqn:2nd_law} with \cref{eqn:potentials} yields
\begin{align}
\label{eqn:2nd_law_reduced}
    \mathcal{D}_{\omega} = Y \ \dot{\omega} &\geq 0,
    &
    \boldsymbol{\sigma} &= 2 \rho \frac{\partial \psi}{\partial \boldsymbol{b}} \cdot \boldsymbol{b},
\end{align}
where $Y(\vec{x}, t) = - \rho \ \partial \psi/\partial \omega$ is the damage energy release rate (thermodynamic force conjugate to the damage variable $\omega$) and $\mathcal{D}_{\omega}$ is the damage dissipation. From \cref{eqn:potentials,eqn:2nd_law_reduced} the Cauchy stress becomes
\begin{align}
\label{eqn:damage_stress}
    \boldsymbol{\sigma} &= (1 - \omega) \  \left[ \frac{\kappa}{2} \left(e^{(\mathrm{J}-1)} - \frac{1}{\mathrm{J}} \right) \boldsymbol{1} + \mu \ \mathrm{J}^{-5/3} \ \boldsymbol{b}' \ \right],
\end{align}
where $\boldsymbol{b}'(\vec{x}, t)$ is the deviatoric part of $\boldsymbol{b}(\vec{x}, t)$.

The evolution of the damage variable $\omega$ is derived in a similar fashion as in plasticity, \cite{matous2003,matous2008,matous2009,matous2010,matous2015,matous2021} First, an energy-based damage criterion is defined by
\begin{align}
\label{eqn:damage_surface}
    g(Y, \chi) &= G(Y) - \chi \leq 0.
\end{align}
The internal state variable $\chi \in [0, 1]$ sets the threshold for how much energy is required in order to increase the amount of damage in the material. The function $G( Y ) \in [0, 1]$ models the damage process. In this work, it is defined as the Weibull distribution
\begin{align}
\label{eqn:weibull}
    G( Y ) &= 1 - \exp{\left[ - \left( \frac{Y - Y_{\mathrm{in}}}{p_{1} Y_{\mathrm{in}}} \right)^{p_{2}} \right]},
\end{align}
though alternative definitions are admissible (\emph{e.g.}, \cite{historyJu,ju1989,historySimoJu,historySimoJu2}). The material parameter $Y_{\mathrm{in}}$ is the energy barrier required to initiate damage and the parameters $p_{1}$ and $p_{2}$ are scale and shape parameters, which facilitate modeling a wide range of material behaviors.

The irreversible dissipative damage evolves according to two evolution equations
\begin{align}
    \label{eqn:rate_independent}
    \dot{\omega} &= \dot{\varkappa} \frac{\mathrm{d} G(Y)}{\mathrm{d} Y},
    &
    \dot{\chi} &= \dot{\varkappa} \frac{\mathrm{d} G(Y)}{\mathrm{d} Y},
\end{align}
where $\dot{\varkappa} = \dot{Y}$ is a damage consistency parameter determined by enforcing the consistency condition $\dot{g} = 0$. The damage loading and unloading are described by the Kuhn-Tucker complimentary conditions
\begin{align}
    \label{eqn:kuhn_tucker}
    \dot{\varkappa} &\geq 0,
    &
    g(Y, \chi) &\leq 0,
    &
    \dot{\varkappa} \ g(Y, \chi) &= 0.
\end{align}
There are well-posedness and uniqueness problems associated with \cref{eqn:rate_independent} which lead to numerical results that do not converge with respect to the spatial discretization. In order to maintain material ellipticity, nonlocal models \cite{historyBazant} or viscous regularization \cite{justOmega} can be introduced. In this work, a rate dependent viscous damage model is obtained by replacing \cref{eqn:rate_independent} with
\begin{align}
    \label{eqn:rate_dependent}
    \dot{\omega} &= \mu_{\omega} \ \langle \phi(g) \rangle,
    &
    \dot{\chi} &= \mu_{\omega} \ \langle \phi(g) \rangle,
\end{align}
where $\mu_{\omega}$ is the damage viscosity, $\phi(g)$ is the damage flow function, and $\langle \bullet \rangle$ are McAuley brackets. Assuming a linear flow function (\emph{i.e.}, $\phi(g) = g$), the Eulerian description for the damage evolution equations becomes
\begin{align}
\label{eqn:damagePDEs}
    \frac{\partial \omega}{\partial t} &= -\nabla \omega \cdot \vec{v} + \mu_{\omega} \ \langle G(Y) - \chi \rangle,
    &
    \frac{\partial \chi}{\partial t} &= -\nabla \chi \cdot \vec{v} + \mu_{\omega} \ \langle G(Y) - \chi \rangle.
\end{align}
This Eulerian-Lagrangian continuum model for the high-strain rate damage and its numerical implementation using the wavelet solver described in \cref{sec:MRWT} is one of the novelties of this work.
\section{Numerical results}
\label{sec:damage_results}
This section presents results from modeling high-strain rate damage under finite strains using equations from \cref{sec:model} and the numerical method from \cref{sec:MRWT}. Specifically \cref{eqn:momentum,eqn:displacement} model the motion of material in the Eulerian frame, \cref{eqn:mass} provides a kinematically consistent density, and is closed by the stress description in \cref{eqn:damage_stress} and the damage evolution in \cref{eqn:damagePDEs}. \Cref{tab:materialParameters} lists the parameters used to model PMMA \cite{GMelastic} with a peak stress $\sigma_{\mathrm{max}} \approx 100$ MPa and an energy release rate $G_{I} \approx 350$ J/m$^{2}$. These values are consistent with experimental measurements provided in \cite{calibrateDamage}. We note that the viscous damage model introduces the characteristic damage length scale $l_{\omega} \propto c_{l} / \mu_{\omega}$, where $c_{l} = \sqrt{E/\rho_{0}}$ is the longitudinal wave speed and $E$ is Young's modulus \cite{matous2015,matous2021}. In this work, we set the damage viscosity $\mu_{\omega}$ such that the damage length scale is consistent with material data (\emph{i.e.}, $l_{\omega} \approx \mathcal{O}(100 \mu$m)) \cite{PMMAlength}.
\begin{table}[!htb]
    \begin{center}
        \caption{Model parameters for dynamic damage of PMMA.\label{tab:materialParameters}}
        \begin{tabular}{ c c c }
        \toprule
        Variable & Name & Value \\ 
        \midrule
        $\rho_{0}$ & Density & $1183 \ \mathrm{kg/m}^{3}$
        \\
        $\kappa$ & Bulk modulus & $5.86 \ \mathrm{GPa}$
        \\ 
        $\mu$ & Shear modulus & $2.29 \ \mathrm{GPa}$
        \\ 
        $\mu_{\omega}$ & Damage viscosity & $10^{7} \ \mathrm{1/s}$
        \\ 
        $Y_{\mathrm{in}}$ & Initiation threshold & $6 \times 10^{4} \ \mathrm{J/m}^{3}$
        \\ 
        $p_{1}$ & Scale parameter & $30$
        \\ 
        $p_{2}$ & Shape parameter & $1.01$
        \\
        \bottomrule
        \end{tabular}
    \end{center}
\end{table}

\begin{figure}[!htb]
    \centering
    \includegraphics[width=0.5\textwidth]{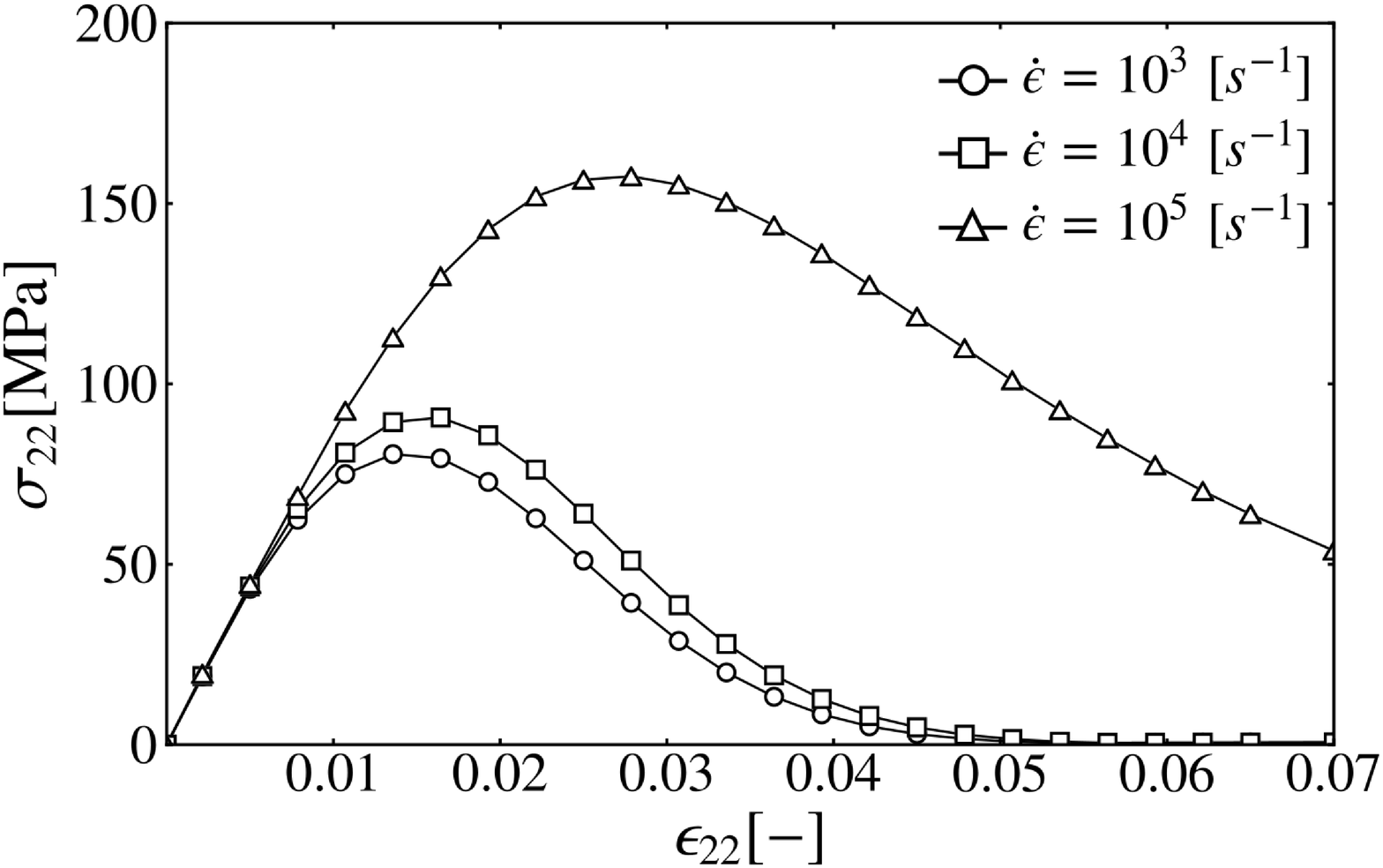}
    \caption{Stress-strain response under uniaxial tension for a variety of strain rates. \label{fig:rateDep}}
\end{figure}

\FloatBarrier
\Cref{fig:rateDep} shows the rate-dependant stress-strain response using material parameters from \cref{tab:materialParameters} in a uniaxial tension setting. These results are consistent with experimental measurements of PMMA exhibiting higher peak tensile strength under dynamic loading \cite{PMMArateDep} with strain rate effects \cite{PMMAsigMax}. 

\begin{figure}[!htb]
    \centering
    \includegraphics[width=0.8\textwidth]{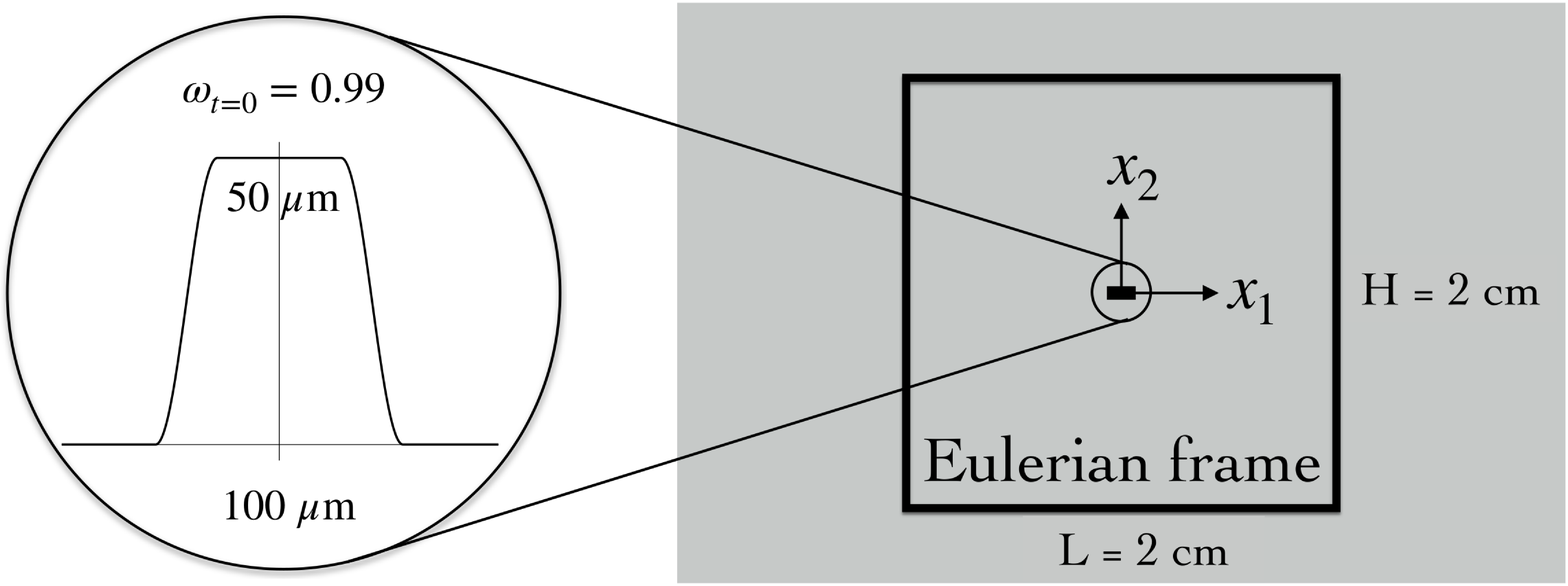}
    \caption{Initial conditions for plane-strain mode-I damage modeling in the Eulerian frame. \label{fig:damageIC}}
\end{figure}

A model of mode-I fracture under plane-strain deformation is defined through the initial conditions
\begin{align}
\vec{v}(\vec{x},0) &= \{0, \dot{\epsilon} \ x_{2} , 0\},
&
\vec{u}(\vec{x},0) &= \vec{0},
&
\rho(\vec{x},0) &= \rho_{0},
\end{align}
with a high strain rate of $\dot{\epsilon} = 5,000$ s$^{-1}$. \Cref{fig:damageIC} depicts the Eulerian domain, defined by a $2$ cm square, containing an initial flaw proportional to the characteristic length of viscous damage (\emph{i.e.}, $\mathcal{O}(100 \ \mu$m)).
In this example, the sparse multiresolution spatial discretization was restricted to no more than $J = 14$ resolution levels which provides continuous resolution between $\mathcal{O}(100 \ \mathrm{nm})$ on the finest level and $\mathcal{O}(1 \ \mathrm{cm})$ on the overall domain (\emph{i.e.}, $\mathcal{O}(10^{5})$ range of length-scales). This range of spatial resolution is impressive in the context of damage simulations as many other solution strategies are either limited to small domains or discretize on a scale that is unable to resolve the characteristic crack thickness (\emph{e.g.}, typical simulations have a range of length-scales around $\mathcal{O}(10^{1})$ \cite{Miller1999,Zhang2007} or $\mathcal{O}(10^{2})$ \cite{Needleman1994,matous2015,matous2021}). 

\FloatBarrier

\begin{figure}[!htb]
    \centering
    \includegraphics[width=0.8\textwidth]{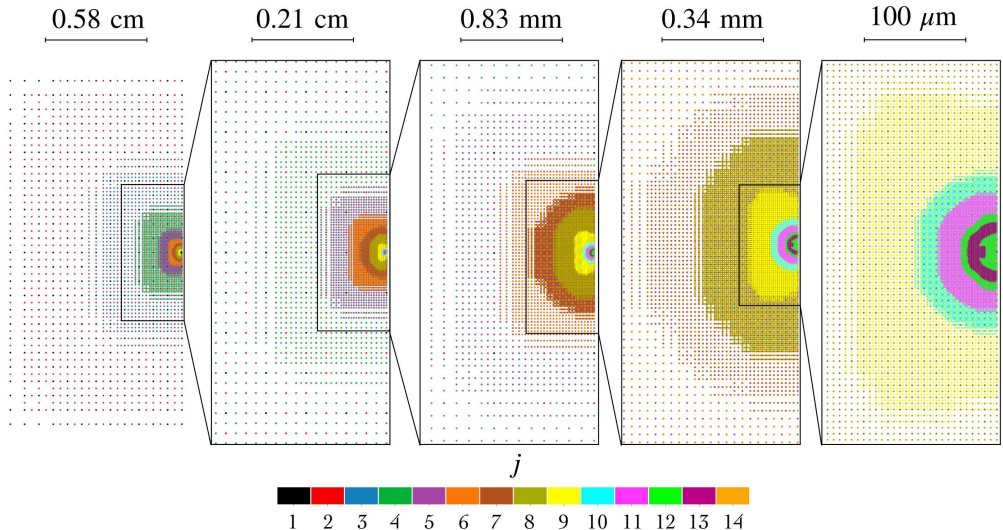}
    \caption{Sparse multiresolution spatial discretization at $t = 0 \ \mathrm{s}$, demonstrating the multiscale capability of MRWT. Note that only half of the domain is shown.\label{fig:gridScale}}
\end{figure}

The remaining figures in this section display simulation results only on the half-domain, obtained from running on $960$ cores. While the embedded Runge-Kutta temporal integration provides a dynamic time-step size according to \cref{eqn:dtNew,eqn:dtBound}, in general $\Delta t \approx \mathcal{O}(10^{-11} \ \mathrm{s})$ requiring approximately $200$ thousand time steps was used to reach $2 \ \mu$s. As shown in \cref{fig:gridScale}, our algorithm adaptively provides resolution only where it is needed as defined by the user-prescribed $\varepsilon$ threshold parameters, which for all results in this section are
\begin{align}
\label{eqn:epsilons}
    \varepsilon_{\omega} &= \varepsilon_{\chi} = 10^{-3},
    &
    \varepsilon_{\vec{u}} &= 10^{-7} \ \mu\mathrm{m},
    &
    \varepsilon_{\vec{v}} &= 10^{-2} \ \mathrm{m/s}.
\end{align}
The MRWT discretization provides dynamic compression of the solution to the PDE. For example, with $p = 6$ the initial condition is discretized on a sparse grid of $35,689$ collocation points on the full domain, whereas a uniform grid with $J = 14$ resolution levels would have $38,655,098,881$ points. Furthermore, this extreme compression continues throughout the simulation where after $2 \ \mu$s of deformation, our sparse multiresolution discretization contains only $810,633$ collocation points on the full domain.

\begin{figure}[!htb]
\begin{tabular*}{\textwidth}{@{} c @{\extracolsep{\fill}} c @{}}
    \includegraphics[width=0.48\textwidth]{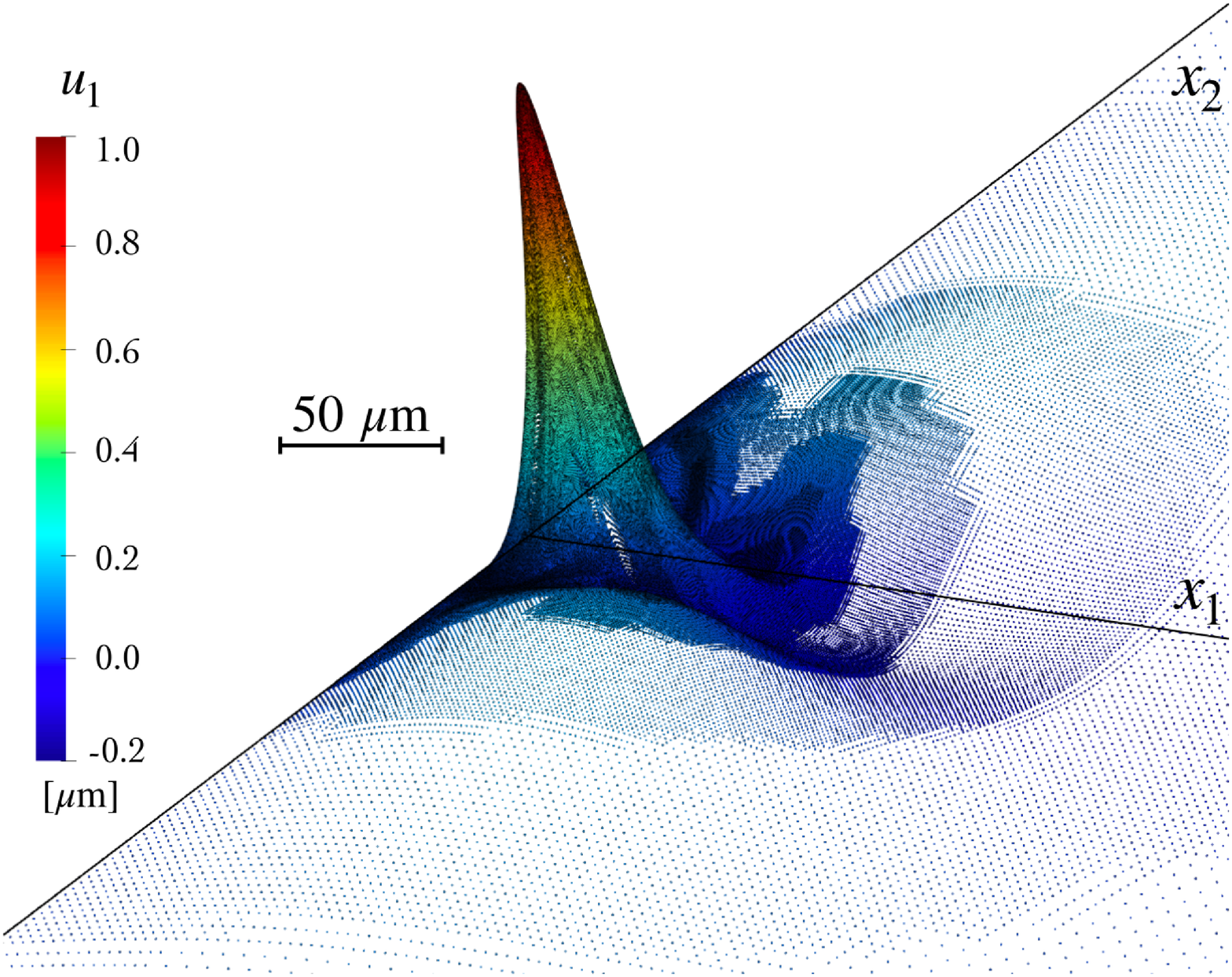}&
    \includegraphics[width=0.48\textwidth]{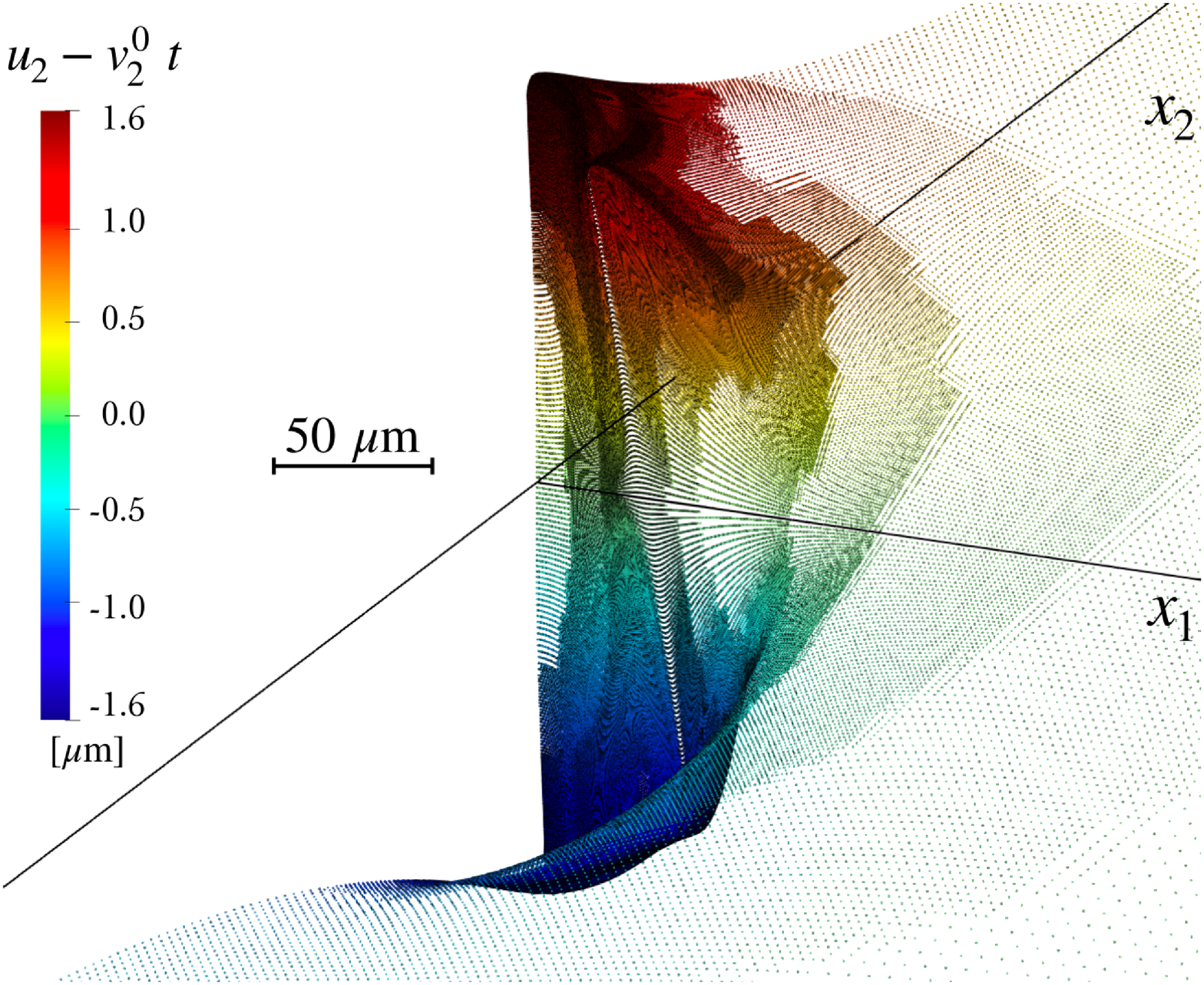}\\
    (a) & (b)
\end{tabular*}
\caption{Components (a) $u_{x}$ and (b) $u_{y}$ of the displacement vector at $t=2.0 \ \mu$s .\label{fig:u_MRWT}}
\end{figure}

The deformation around the initial flaw induces strain localization as exhibited by the concentrated displacement in \cref{fig:u_MRWT}. Note that the background displacement of mode-I opening (\emph{i.e.}, $u_{2}^{\mathrm{bg}} \approx v_{2}^{0} \times t$) is subtracted from the $u_{2}$ component of displacement to reveal solution features (\emph{e.g.}, pronounced displacement of material points in the vicinity of the flaw). 

\begin{figure}[!htb]
\begin{tabular*}{\textwidth}{@{} c @{\extracolsep{\fill}} c @{}}
    \includegraphics[width=0.48\textwidth]{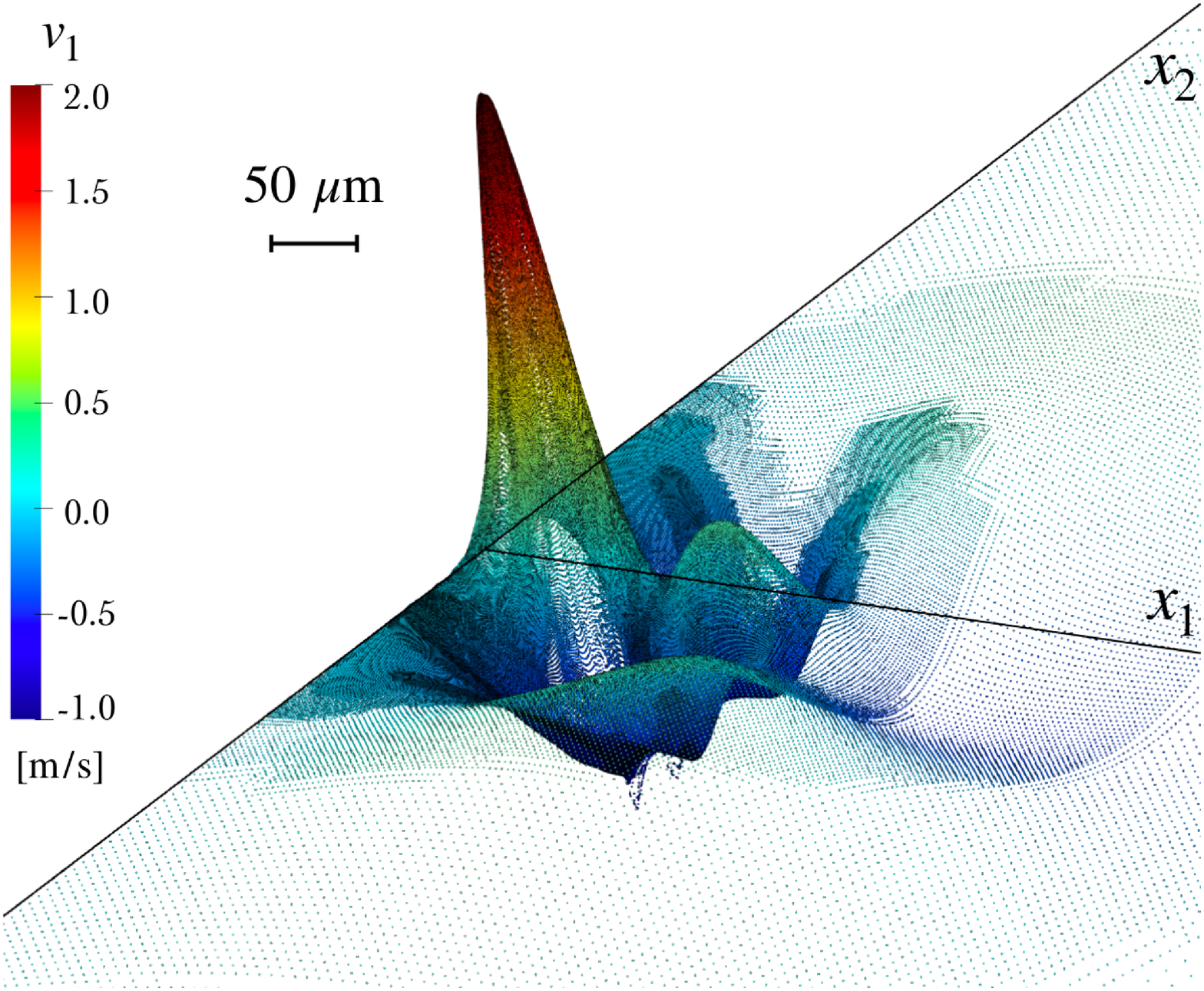}&
    \includegraphics[width=0.48\textwidth]{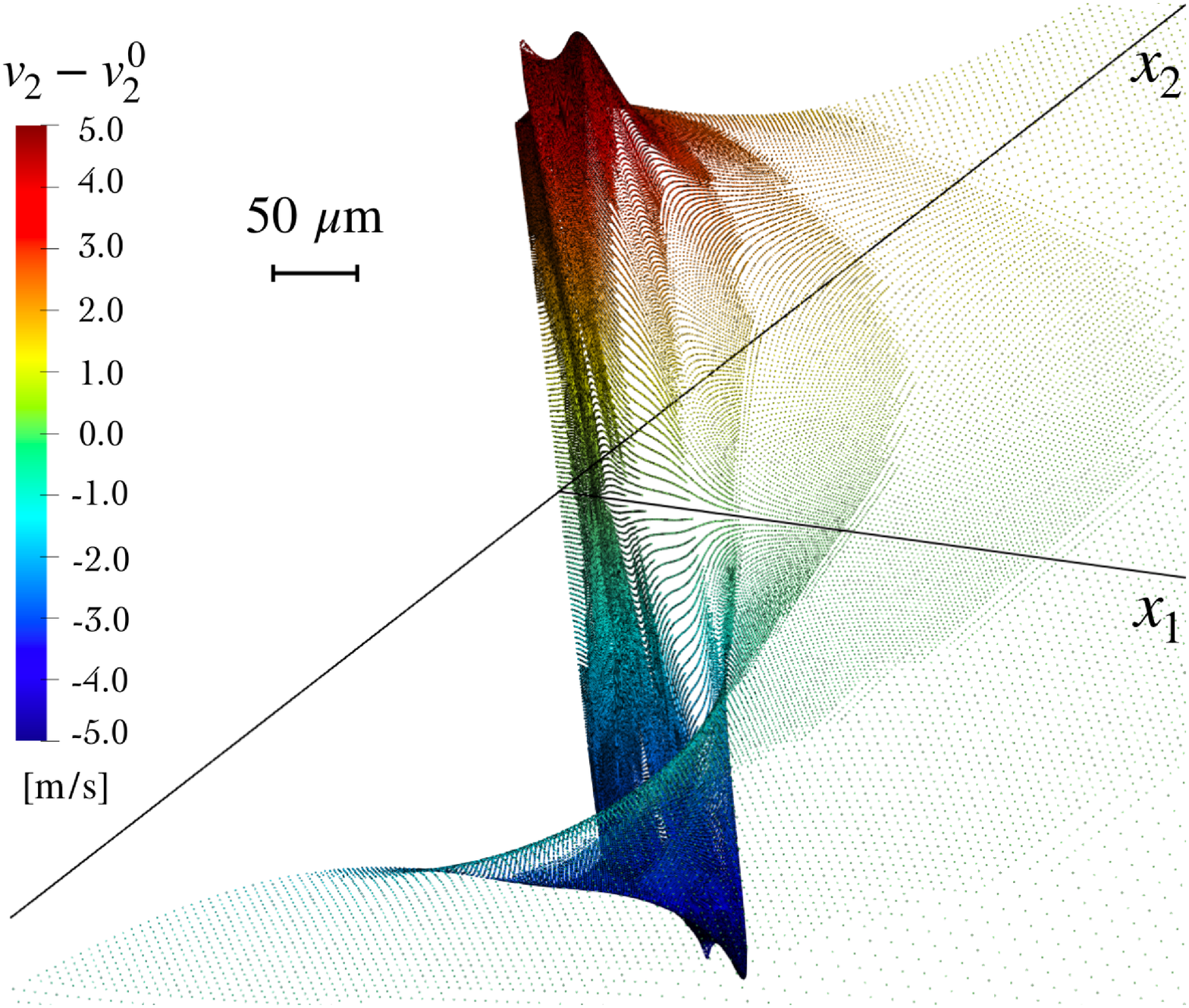}\\
    (a) & (b)
\end{tabular*}
\caption{Components (a) $v_{x}$ and (b) $v_{y}$ of the velocity vector at $t=2.0 \ \mu$s .\label{fig:v_MRWT}}
\end{figure}

\Cref{fig:v_MRWT} illustrates the dynamics of this process with significant spikes in the velocity field (\emph{e.g.}, $\mathcal{O}(1) \ \mathrm{m/s}$) over a very small region (\emph{e.g.}, $\mathcal{O}(10^{-5}) \ \mathrm{m}$). Again, the background velocity of mode-I opening (\emph{i.e.} $v_{2}^{\mathrm{bg}} \approx v_{2}^{0}$) is subtracted from the $v_{2}$ component of velocity to reveal solution features (\emph{e.g.}, symmetric local maxima/minima above/below the damage profile).

\begin{figure}[!htb]
\centering
    \includegraphics[width=0.55\textwidth]{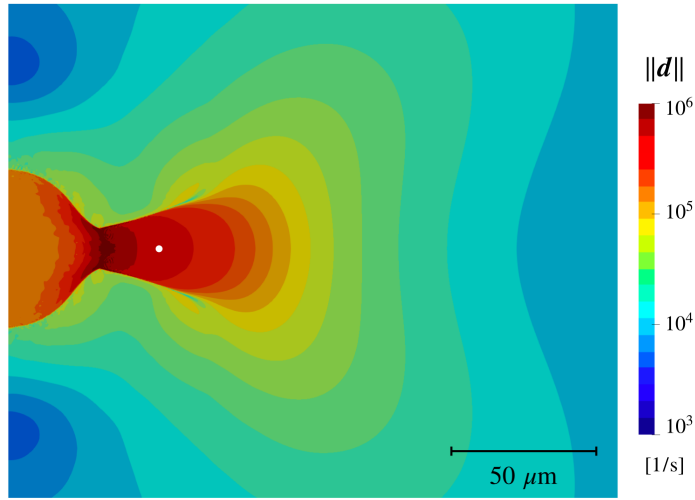}
\caption{Frobenius norm of the rate of the deformation tensor $\| \boldsymbol{d} \|$ at $t=2.0 \ \mu$s. The white circle indicates the spatial coordinate $\vec{x} = [52\sfrac{1}{12},0,0] \ \mu$m where the constitutive response is probed. \label{fig:d_MRWT}}
\end{figure}

\begin{figure}[!htb]
\begin{tabular*}{\textwidth}{@{} c @{\extracolsep{\fill}} c @{}}
    \includegraphics[width=0.48\textwidth]{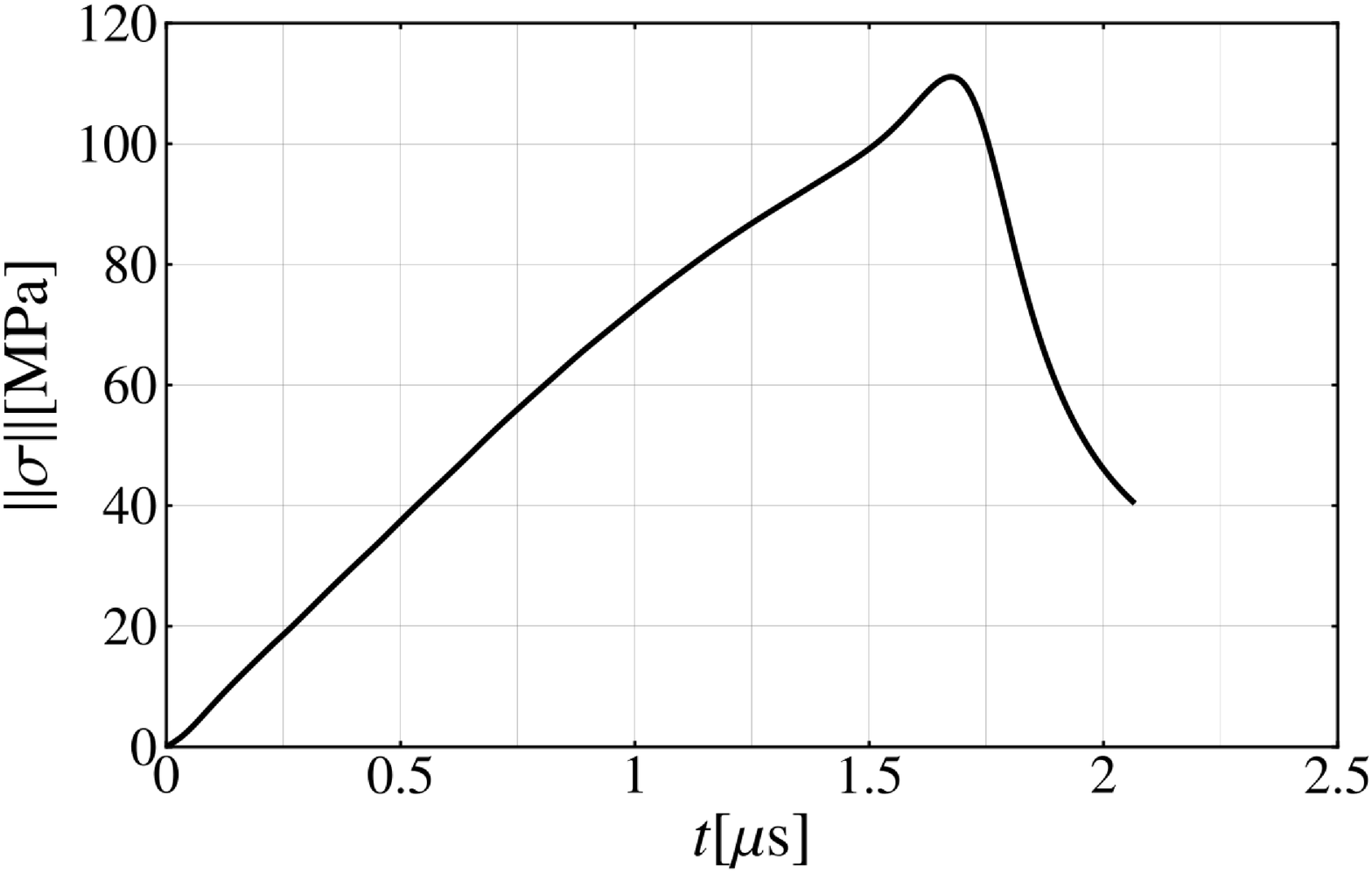}&
    \includegraphics[width=0.49
\textwidth]{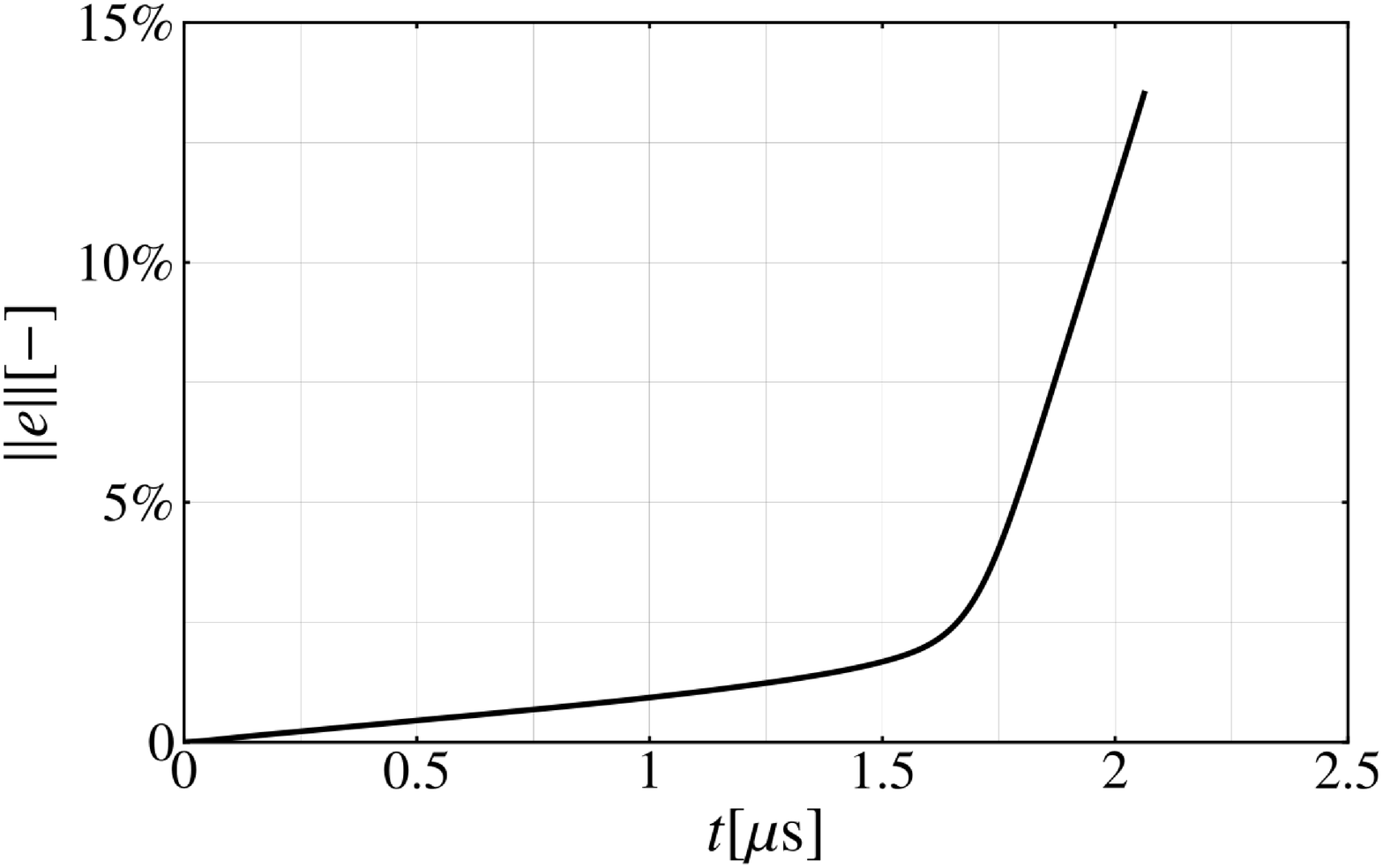}\\
    (a) & (b)
\end{tabular*}
\caption{MRWT computed values of (a) Frobenius norm of the Cauchy stress and (b) Frobenius norm of the Almansi-Hamel strain over time, obtained at a specific spatial location $\vec{x} = [52\sfrac{1}{12},0,0] \ \mu$m chosen because it is outside predamaged region of the initial condition. \label{fig:x_stress_strain_MRWT}}
\end{figure}

\Cref{fig:d_MRWT} reveals the extreme nature of the rate of the deformation tensor (\emph{i.e.}, $\boldsymbol{d} = ( \ \nabla \vec{v} \ )^{\mathrm{sym}}$). Since this quantity is proportional to the strain rate, the MRWT simulation suggests that in the damaged regions the material experiences strain rates multiple orders of magnitude higher than the imposed loading rate of $\dot{\epsilon} = 5,000 \ \mathrm{s}^{-1}$. This observation is consistent with the model depicted in \cref{fig:rateDep}.
This consideration is made more clear by observing the temporal evolution of stress and strain at a particular coordinate outside the region of the initial flaw (\emph{i.e.}, $\vec{x} = [52\sfrac{1}{12},0,0] \ \mu$m). \Cref{fig:x_stress_strain_MRWT} plots the stress and strain at this location over time. In this figure, the precipitous decline in the stress reflects the brittle nature of PMMA. Furthermore, there are $2$ distinct strain rates separated by the limit point in the stress-strain curve around $t \approx 1.6 \ \mu$s (\emph{i.e.}, $\dot{\epsilon} \approx 10^{4} \ \mathrm{s}^{-1}$ during hardening and $\dot{\epsilon} \approx 3 \times 10^{5} \ \mathrm{s}^{-1}$ during softening).

\FloatBarrier

\begin{figure}[!htb]
\begin{tabular*}{\textwidth}{@{} c @{\extracolsep{\fill}} c @{}}
    \includegraphics[width=0.48\textwidth]{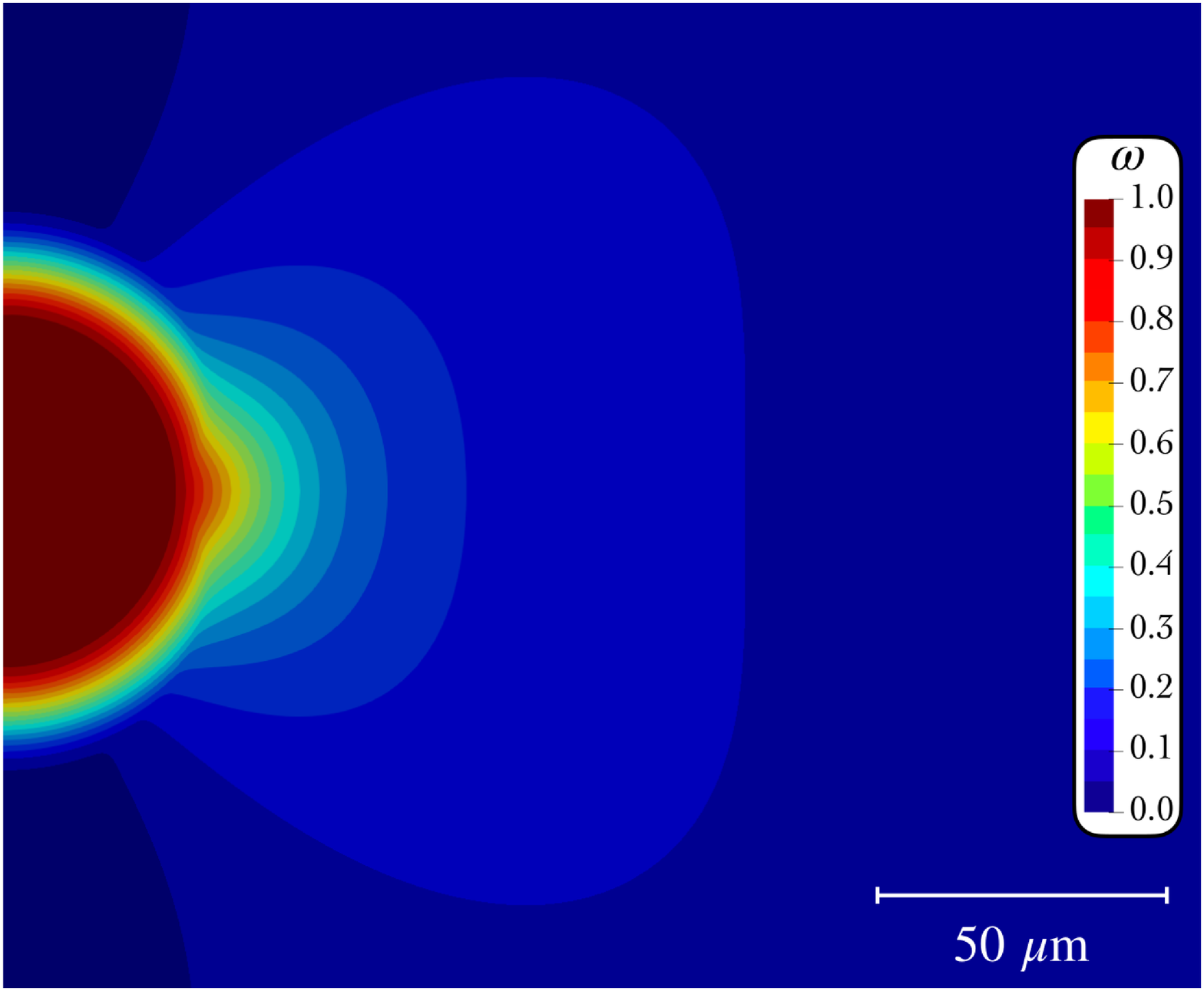}&
    \includegraphics[width=0.48\textwidth]{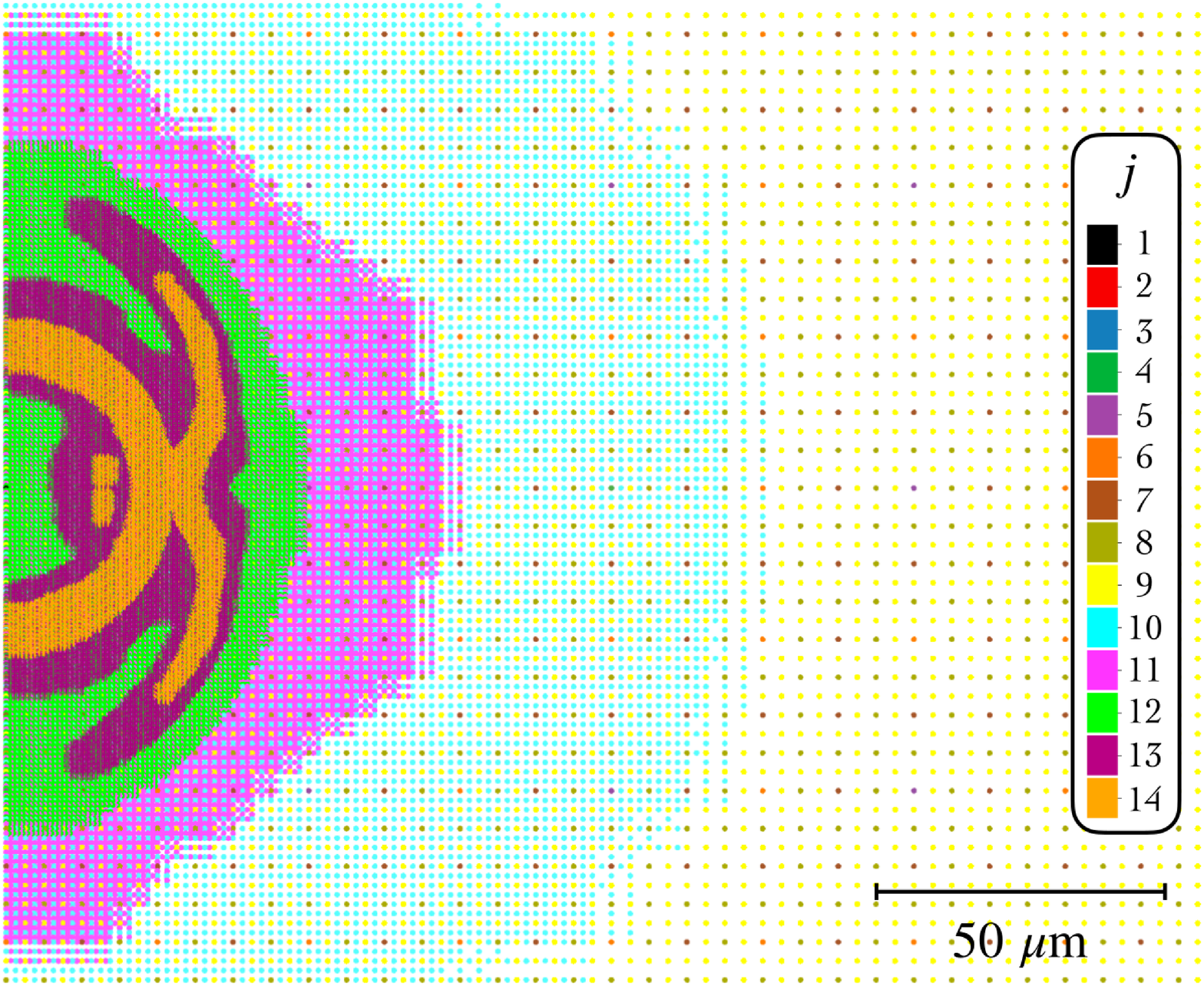}\\
    (a) & (b)\\
    \includegraphics[width=0.48\textwidth]{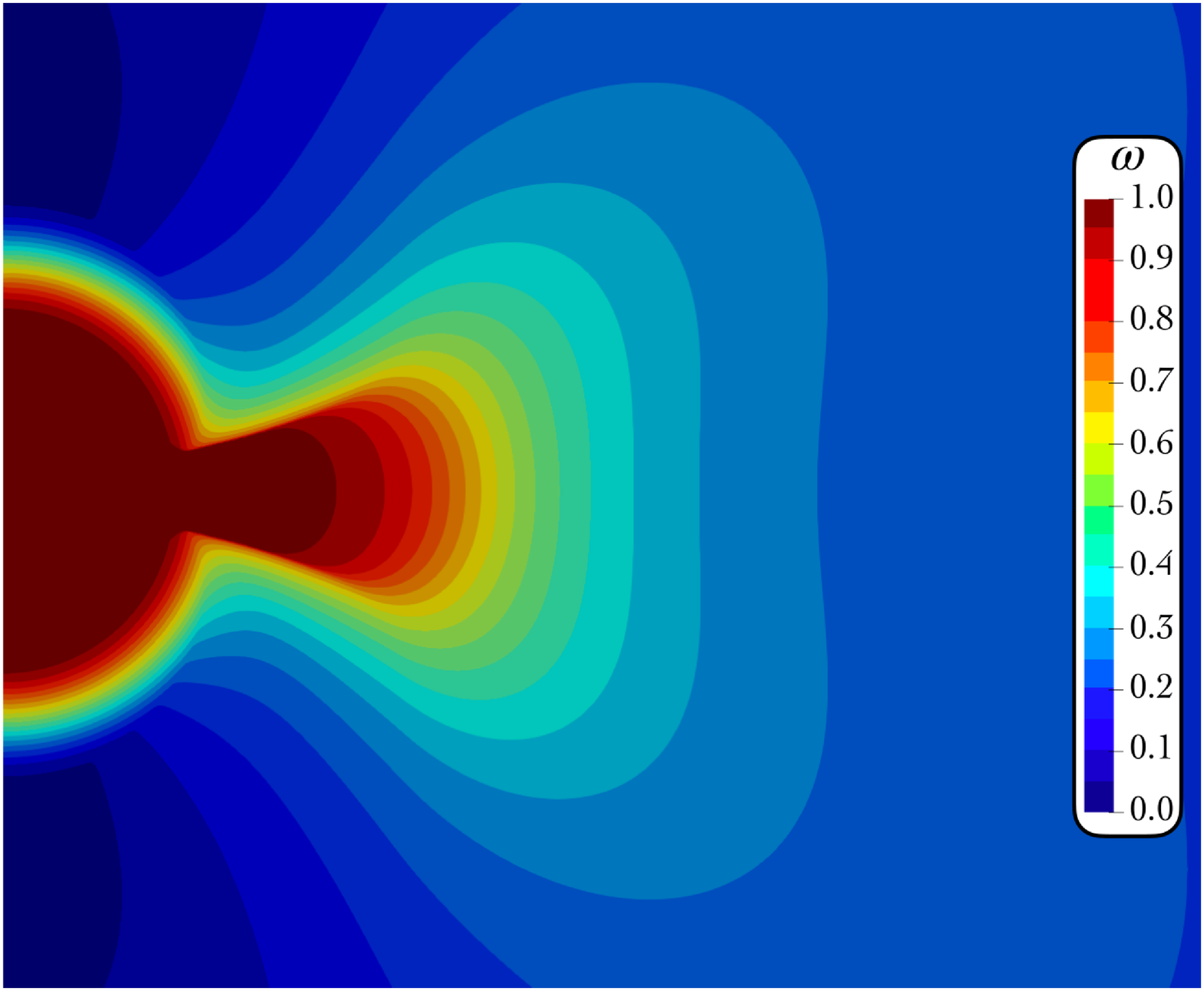}&
    \includegraphics[width=0.48\textwidth]{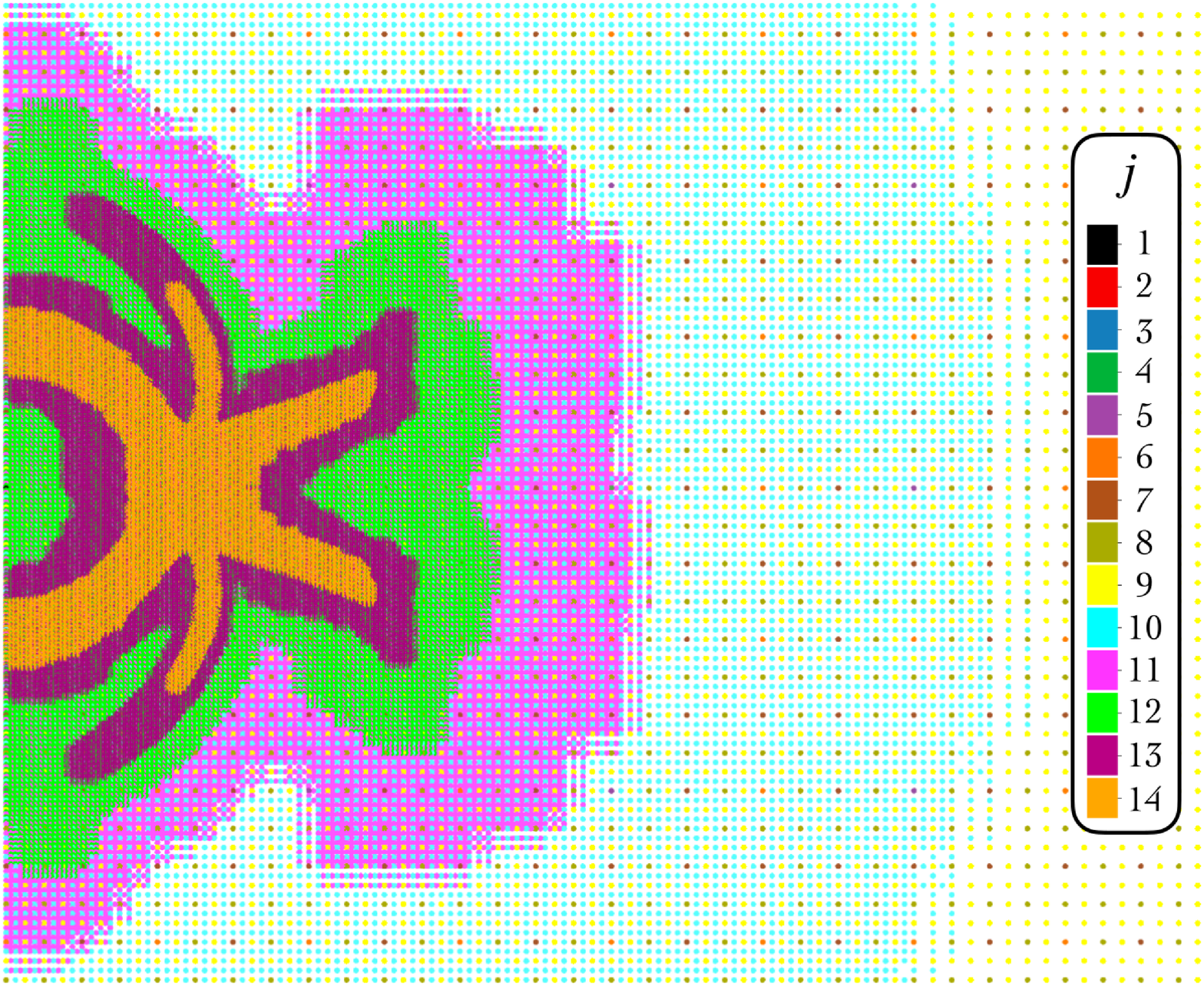}\\
    (c) & (d)
\end{tabular*}
\caption{MRWT generated damage profile and sparse multiresolution grid for the dynamic damage of PMMA at various times: (a) \& (b) $t=1.5 \ \mu$s and (c) \& (d) $t=2.0 \ \mu$s. The grid points are colored according to their resolution level $j$. The reader is referred to the online version of this article for clarity regarding the color in this figure.\label{fig:damage_grid_MRWT}}
\end{figure}

\begin{figure}[!htb]
\begin{tabular*}{\textwidth}{@{} c @{\extracolsep{\fill}} c @{}}
    \includegraphics[width=0.48\textwidth]{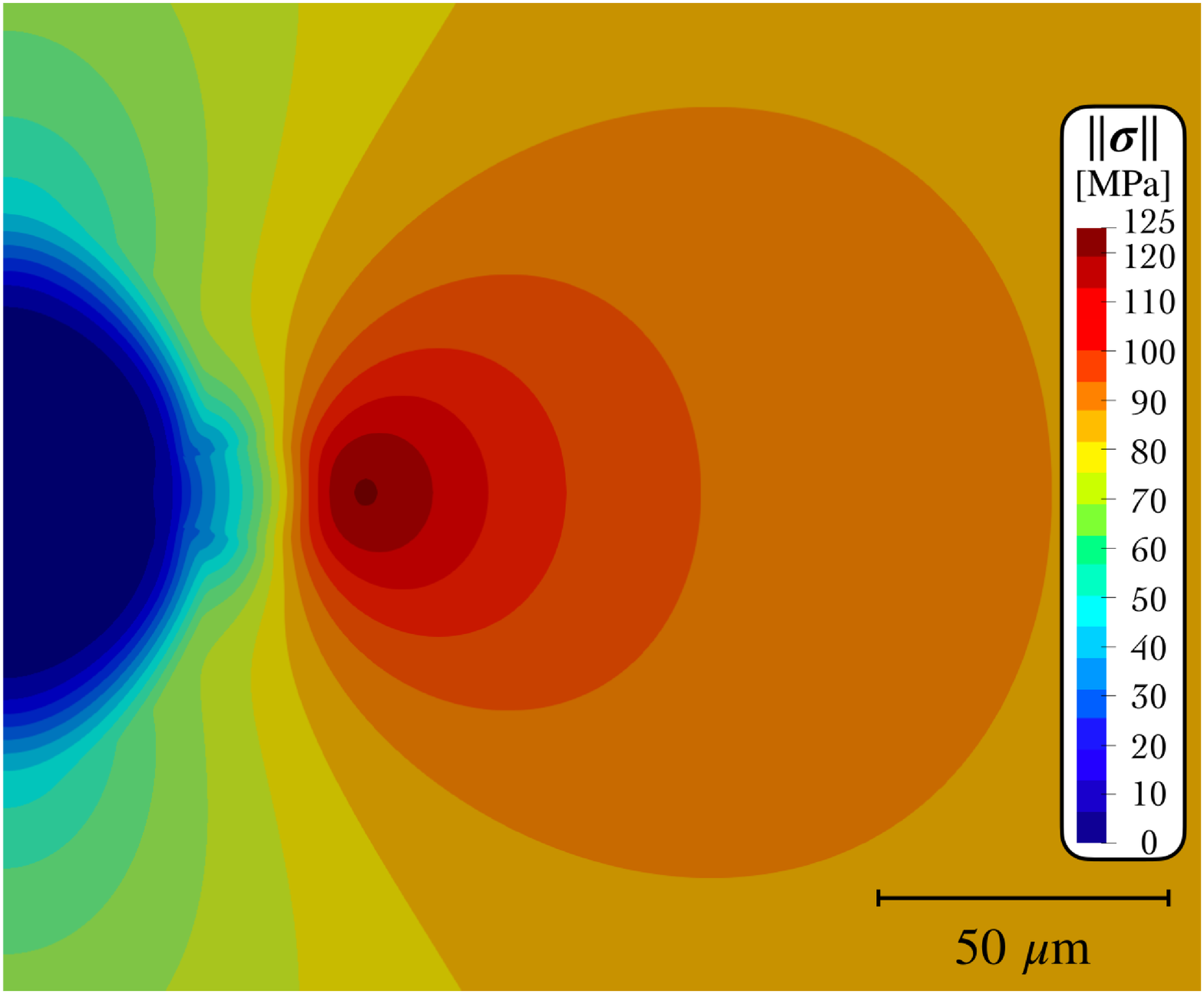}&
    \includegraphics[width=0.48\textwidth]{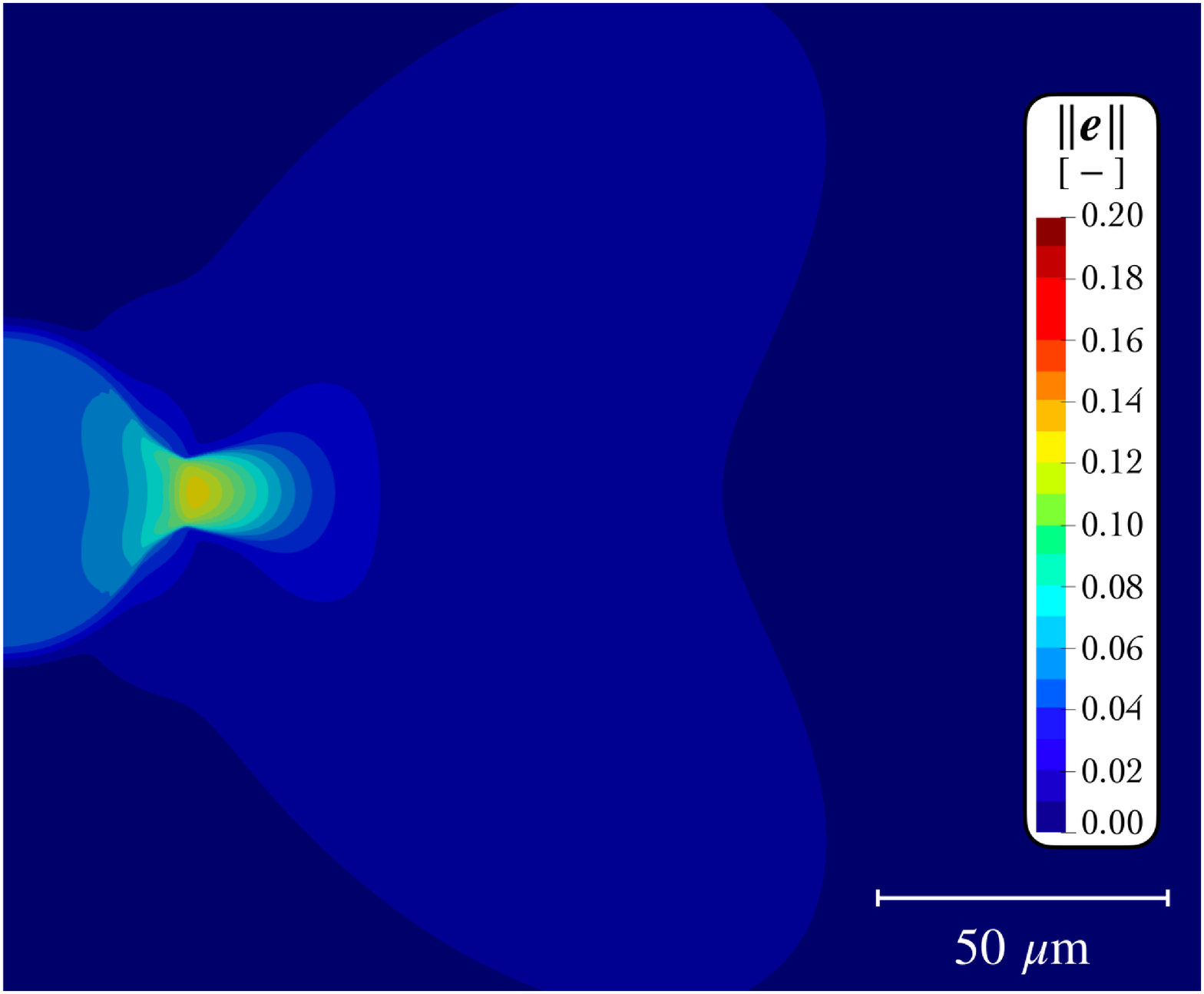}\\
    (a) & (b)\\
    \includegraphics[width=0.48\textwidth]{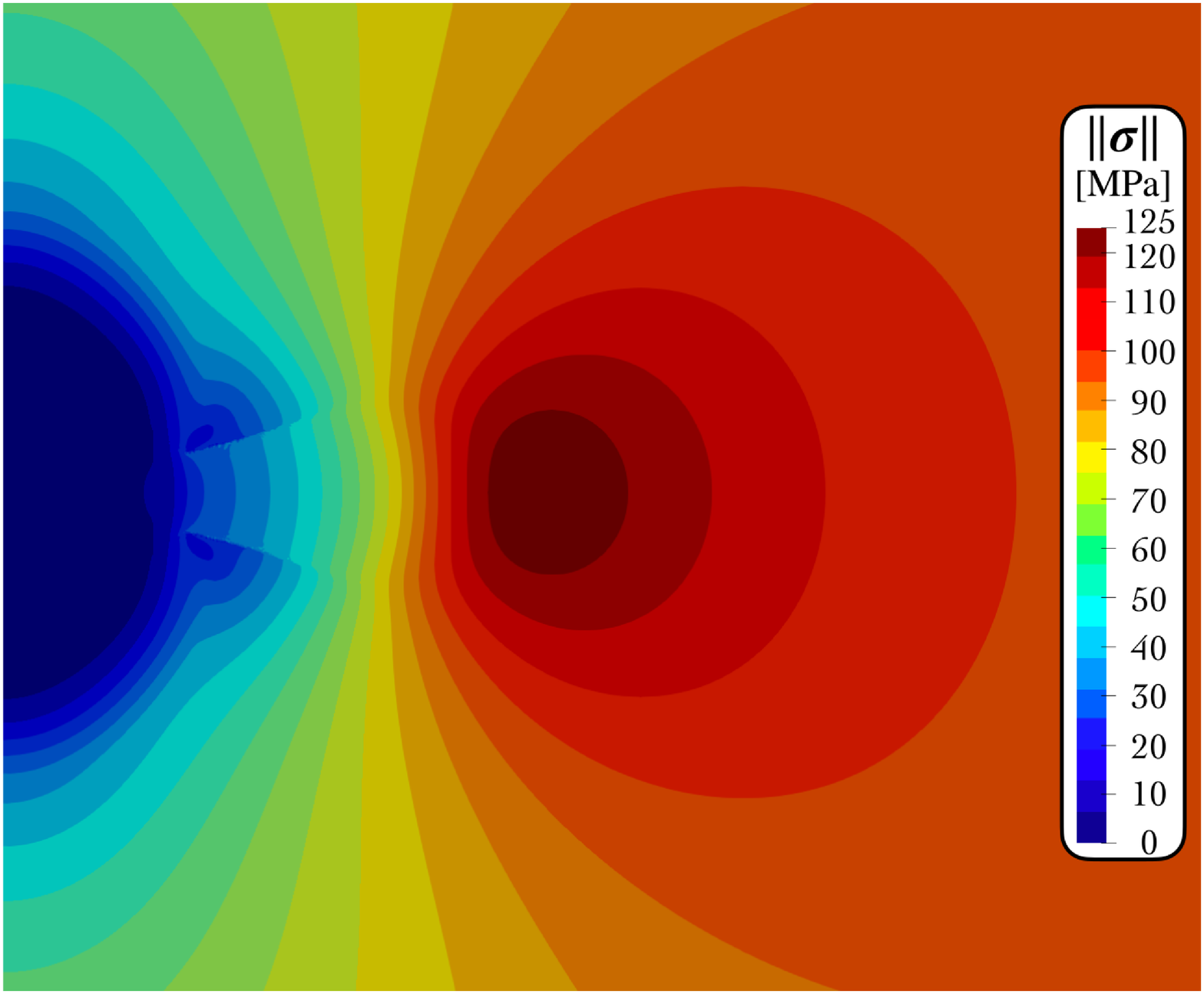}&
    \includegraphics[width=0.48\textwidth]{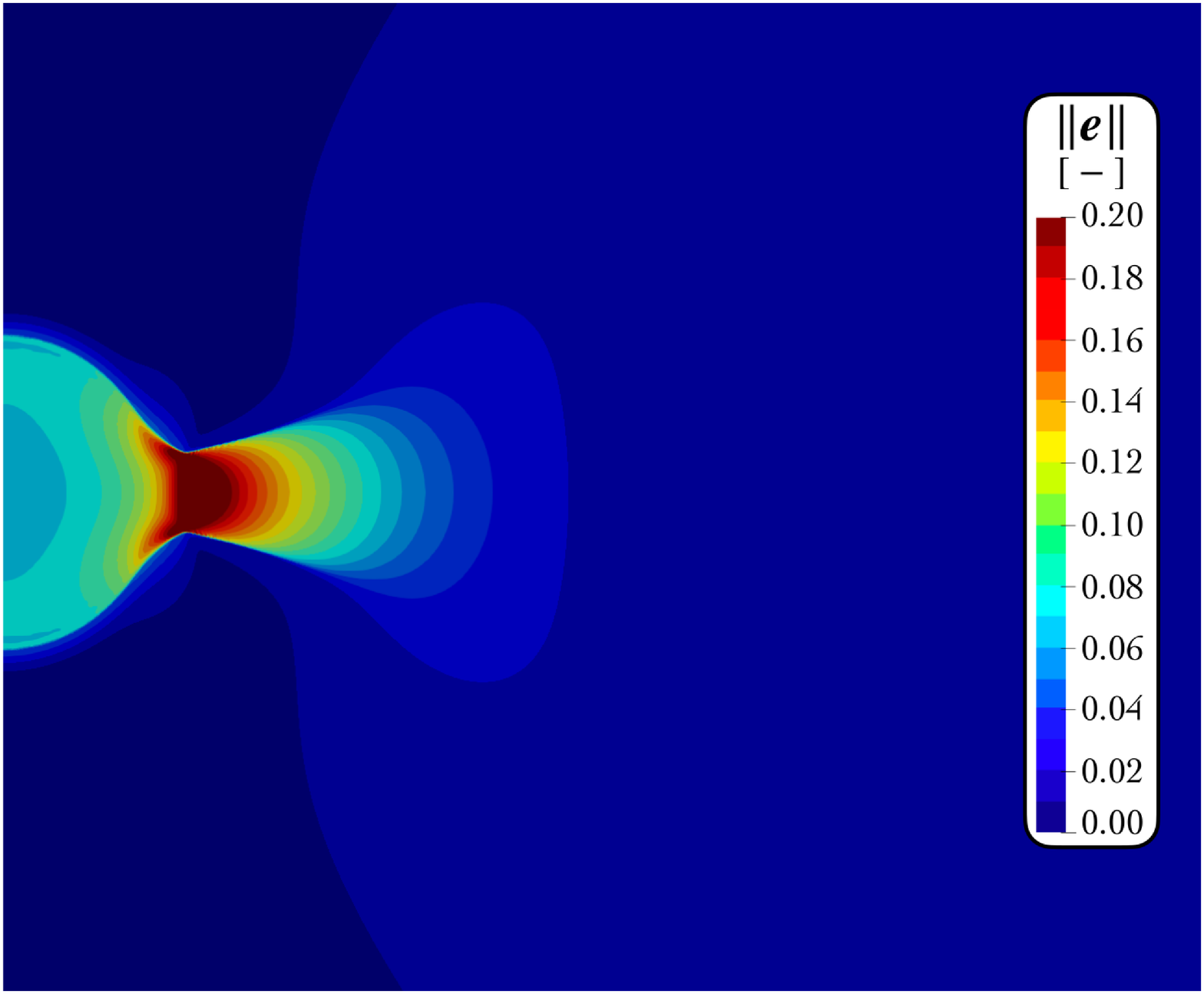}\\
    (c) & (d)
\end{tabular*}
\caption{MRWT generated stress and strain profiles for the dynamic damage of PMMA at various times: (a) \& (b) $t=1.75 \ \mu$s and (c) \& (d) $t=2.0 \ \mu$s.\label{fig:stress_strain_MRWT}}
\end{figure}

This strain localization around the initial flaw provides the energy to initiate and grow damage in a process zone. \Cref{fig:damage_grid_MRWT} shows the damage profile as it evolves and illustrates the dynamic adaptivity of the MRWT spatial discretization. A stress concentration develops in front of the process zone which propagates in a direction orthogonal to the applied deformation. \Cref{fig:stress_strain_MRWT,fig:x_stress_strain_MRWT} demonstrate that the damaged material in the wake of the stress concentration has a reduced capacity to support a load and experiences significantly higher strains, reaching approximately $20 \%$.

According to linear fracture mechanics \cite{Freund}, material deformed in this configuration will generate a crack of length $a(t)$ that grows with a speed $\dot{a}(t)$ which accelerates until reaching the Rayleigh speed $c_{R}$. For the material parameters in \cref{tab:materialParameters}, $c_{R} \approx 1,295$ m/s. An estimate for the length of the crack from the MRWT simulation is provided by measuring the distance from the origin to the peak stress (\emph{e.g.}, \cref{fig:stress_strain_MRWT}). From the initial condition, the crack starts with a nonzero length $a_{0} \approx 46.4 \ \mu$m. After $2.2 \ \mu$s the crack has grown to a length of $a \approx 135.1 \ \mu$m. The corresponding speed of the crack from the MRWT simulation is $\dot{a} \approx 280$ m/s, indicating the early stage of the damage process. The crack speed was computed from the slopes of quadratic functions locally fit to the crack length over time. 

\begin{figure}[!htb]
\centering
    \includegraphics[width=0.75\textwidth]{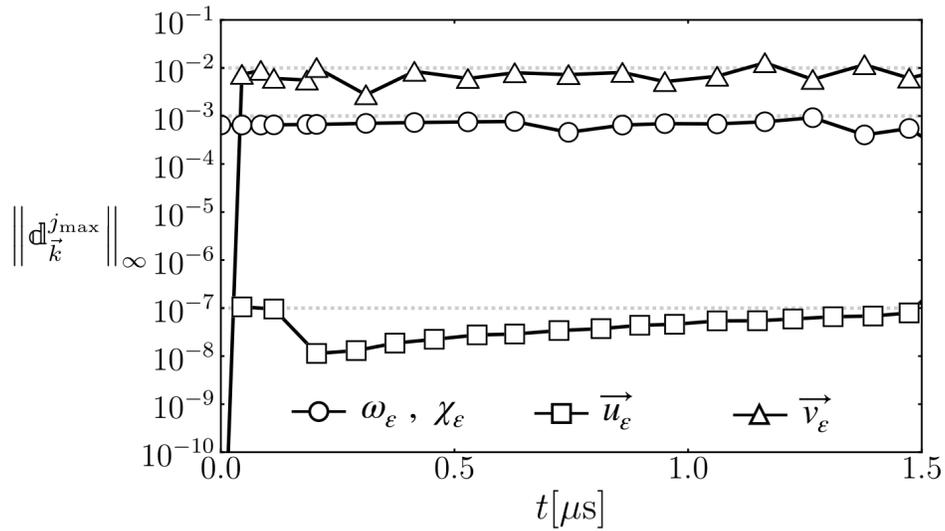}
\caption{Solution verification of dynamic damage using \cref{eqn:errorEstimate}, dashed lines indicate each field's epsilon threshold according to \cref{eqn:epsilons}.\label{fig:dMags}}
\end{figure}

Lastly, our novel numerical method can provide assessment of the
quality of the solution (\emph{i.e.}, the solution
verification). Whereas traditional numerical techniques require a
sequence of solutions on multiple grids, this work provides
self-verification from a single simulation. Specifically, with $J =
14$ resolution levels available to our algorithm, the error of the
solution of the PDEs is bound by \cref{eqn:ferror,eqn:product_error}
and approximated by \cref{eqn:D_error} as long as $j_{\mathrm{max}} <
J$. Beyond this limit (\emph{i.e.}, $j_{\mathrm{max}} \geq J$), our
algorithm provides automatic estimates of the error through the
magnitude of the wavelet coefficients on the highest resolution level,
\begin{align}
\label{eqn:errorEstimate}
    \mathrm{error}(\bullet) \propto \left\| {}^{\lambda}\mathbb{d}_{\vec{k}}^{J}(\bullet) \right\|_{\infty}.
\end{align}
\Cref{fig:dMags} shows the error estimates from the damage simulation
for the first $1.5 \ \mu$s with $j_{\mathrm{max}} < J$. In this time
interval, each field has an error proportional to their respective epsilons
according to \cref{eqn:epsilons}. After $2 \ \mu$s of deformation when
$j_{\mathrm{max}}$ was no longer guaranteed, the relative percent
error for each field was: $2.7 \%$ for the damage variables, $0.027
\%$ for the displacement, and $9.7 \%$ for the velocity.

\section{Conclusions}
Our proposed wavelet based algorithm is well suited for computational science and engineering applications requiring dynamic resolution across multiple spatial and temporal scales. Whereas traditional numerical methods require \emph{a posteriori} remeshing/refining procedures, the wavelet basis functions provide intrinsic adaptivity with resolution only where and when it is needed according to the user-prescribed accuracy. Although other wavelet based solvers have been developed in the past, our novel algorithm is the first to the best of our knowledge which uses the wavelet operators (\emph{e.g.}, the direct differentiation of the wavelet basis) to inform the grid adaption strategy. In this work, the spatial discretization is guided by the wavelet error estimates for fields, derivatives, and nonlinear operations. The error of the solution of the PDE is controlled in space and time by combining the predictor/corrector grid adaption procedure with embedded Runge-Kutta time integration. Furthermore, the matrix formulation of the wavelet operators allows the MRWT algorithm to be trivially applied to $N$ spatial dimensions and data compression is enhanced by using the minimum number of points on the coarse grid, thereby maximizing the number of collocation points which have the potential to be omitted from the discretization. MRWT has been verified by showing spatial convergence of the solution of nonlinear PDEs at rates that are consistent with the aforementioned error estimates. Additionally, the capabilities of MRWT have been demonstrated by simulating multiscale applications, such as the Navier-Stokes equations.

Finally, we have demonstrated our wavelet numerical method on a difficult high-strain rate dynamic damage problem. The sparse multiresolution discretization provides continuous resolution between $O(100 \ \mathrm{nm})$ on the finest level and $O(1 \ \mathrm{cm})$ on the overall domain. As shown, the simulation of damage nucleation and propagation is well resolved on this discretization spanning $5$ orders of magnitude without imposing significant computational overhead. By providing highly resolved kinematic quantities and associated constitutive response, our algorithm can inform material scientists of the extreme dynamics involved in material failure, thereby facilitating the design of robust materials. Furthermore, since the wavelet algorithm provides accurate simulations of nonlinear dynamic damage it is a convenient tool for the development of more accurate macroscopic damage models (\emph{e.g.}, anisotropic). This work illustrates that our wavelet solver is capable of simulating complex systems with the resolution required to validate material models against experimental results.

\section*{Acknowledgments}
This work was supported by Los Alamos National Laboratory (LANL) under award numbers $369229$ and $370985$ and by the Department of Energy, National Nuclear Security Administration, under award number DENA$0002377$ as part of the Predictive Science Academic Alliance Program II.

This research was also supported in part by Lilly Endowment, Inc., through its support for the Indiana University Pervasive Technology Institute. The authors acknowledge the Indiana University Pervasive Technology Institute for providing supercomputing resources that have contributed to the research results reported within this paper.

\clearpage

\appendix
\section{Wavelet derivatives}
\label{sec:fullDexpression}
This section provides additional details on differentiation in \cref{eqn:D_basis} and \cref{eqn:D_product}. For sufficiently differentiable basis, the derivatives of the basis are also continuous functions and are projected onto the original basis using \cref{eqn:FWT_BWT} with global index vectors $\vec{m}$ and $\vec{n}$
\begin{align}
\label{eqn:D_details_A}
    \frac{\partial^{\alpha}}{\partial x_{i}^{\alpha}} f(\vec{x}) &\approx \mathbb{d}_{m_{1} \ldots \ m_{N}} \frac{\partial^{\alpha} \ \Psi(\vec{x})_{m_{1} \ldots \ m_{N}}}{\partial x_{i}^{\alpha}},
    \nonumber \\
    \frac{\partial^{\alpha}}{\partial x_{i}^{\alpha}} f(\vec{x}) &\approx \mathbb{d}_{m_{1} \ldots \ m_{N}} \left( \mathbb{D}_{n_{1} \ldots \ n_{N}} \ \Psi(\vec{x})_{n_{1} \ldots \ n_{N}} \right)_{m_{1} \ldots \ m_{N}},
    \nonumber \\
    \frac{\partial^{\alpha}}{\partial x_{i}^{\alpha}} f(\vec{x}) &\approx \mathbb{d}_{n_{1} \ldots \ n_{N}} \left( \mathbb{D}_{m_{1} \ldots \ m_{N}} \ \Psi(\vec{x})_{m_{1} \ldots \ m_{N}} \right)_{n_{1} \ldots \ n_{N}}.
\end{align}
Note that in \cref{eqn:D_details_A} the containers $\mathbbb{D}$ and $\boldsymbol{\Psi}(\vec{x})$ exist for each collocation point indexed by $\vec{n}$. While the entries in $\mathbbb{D}$ may vary for each $\vec{n}$, the basis $\boldsymbol{\Psi}(\vec{x})$ are identical, therefore combining like terms yields
\begin{align}
\label{eqn:D_details_B}
    \frac{\partial^{\alpha}}{\partial x_{i}^{\alpha}} f(\vec{x}) &\approx \mathbb{d}_{n_{1} \ldots \ n_{N}} \mathbb{D}_{m_{1} \ldots \ m_{N}}^{n_{1} \ldots \ n_{N}} \ \Psi(\vec{x})_{m_{1} \ldots \ m_{N}},
    \nonumber \\
    \frac{\partial^{\alpha}}{\partial x_{i}^{\alpha}} f(\vec{x}) &\approx \mathbb{D}_{m_{1} \ldots \ m_{N}}^{n_{1} \ldots \ n_{N}} \ \mathbb{d}_{n_{1} \ldots \ n_{N}} \ \Psi(\vec{x})_{m_{1} \ldots \ m_{N}},
    \nonumber \\
    \frac{\partial^{\alpha}}{\partial x_{i}^{\alpha}} f(\vec{x}) &\approx \widetilde{\mathbb{d}}_{m_{1} \ldots \ m_{N}} \ \Psi(\vec{x})_{m_{1} \ldots \ m_{N}}.
\end{align}
Where the $\widetilde{\mathbbb{d}}$ values are wavelet coefficients for the derivative of the field, obtained by contracting wavelet coefficients of the field with connection coefficients (\emph{i.e.}, wavelet coefficients of the derivatives of the basis)
\begin{align}
\label{eqn:D_details_C}
    \widetilde{\mathbb{d}}_{m_{1} \ldots \ m_{N}} &= \mathbb{D}_{m_{1} \ldots \ m_{N}}^{n_{1} \ldots \ n_{N}} \ \mathbb{d}_{n_{1} \ldots \ n_{N}}.
\end{align}
It is helpful to express \cref{eqn:D_details_C} for just one entry in $\widetilde{\mathbbb{d}}$ (\emph{i.e.}, specifying the index $\vec{m}$ identifies a particular ${}^{\lambda}\widetilde{\mathbb{d}}_{\vec{k}}^{ \ j}$). This corresponds to the value of the $\lambda$-type wavelet coefficient for the derivative of the field on resolution level $j$ at spatial index $\vec{k}$
\begin{align}
\label{eqn:D_details_D}
    {}^{\lambda}\widetilde{\mathbb{d}}_{\vec{k}}^{ \ j} &= \sum_{r=1}^{j_{\mathrm{max}}} \sum_{\beta=0}^{2^{N}-1} {}_{\lambda}^{\beta} \mathbb{D}_{j, \vec{k}}^{r, \vec{l}} \; \; {}^{\beta}\mathbb{d}_{\vec{l}}^{r}
    \nonumber\\
    {}^{\lambda}\widetilde{\mathbb{d}}_{\vec{k}}^{ \ j} &= \sum_{r=1}^{j_{\mathrm{max}}} \sum_{\beta=0}^{2^{N}-1}\int \frac{\partial^{\alpha}}{\partial x_{i}^{\alpha}} \bigg( {}^{\beta}\Psi_{\vec{l}}^{r}(\vec{x}) \bigg) \ {}^{\lambda}\widetilde{\Psi}_{\vec{k}}^{j}(\vec{x}) \ \mathrm{d}\vec{x} \; \; {}^{\beta}\mathbb{d}_{\vec{l}}^{r}
    \nonumber\\
    {}^{\lambda}\widetilde{\mathbb{d}}_{\vec{k}}^{ \ j} &= \sum_{r=1}^{j - 1} \sum_{\beta=0}^{2^{N}-1}\int \frac{\partial^{\alpha}}{\partial x_{i}^{\alpha}} \bigg( {}^{\beta}\Psi_{\vec{l}}^{r}(\vec{x}) \bigg) \ {}^{\lambda}\widetilde{\Psi}_{\vec{k}}^{j}(\vec{x}) \ \mathrm{d}\vec{x} \; \; {}^{\beta}\mathbb{d}_{\vec{l}}^{r} \ \ldots 
    \nonumber\\
    &+\sum_{\beta=0}^{2^{N}-1} \int \frac{\partial^{\alpha}}{\partial x_{i}^{\alpha}} \bigg( {}^{\beta}\Psi_{\vec{l}}^{j}(\vec{x}) \bigg) \ {}^{\lambda}\widetilde{\Psi}_{\vec{k}}^{j}(\vec{x}) \ \mathrm{d}\vec{x} \; \; {}^{\beta}\mathbb{d}_{\vec{l}}^{j} \ \ldots 
    \nonumber\\
    &+ \sum_{r=j + 1}^{j_{\mathrm{max}}} \sum_{\beta=1}^{2^{N} - 1} \int \frac{\partial^{\alpha}}{\partial x_{i}^{\alpha}} \bigg( {}^{\beta}\Psi_{\vec{l}}^{r}(\vec{x}) \bigg) \ {}^{\lambda}\widetilde{\Psi}_{\vec{k}}^{j}(\vec{x}) \ \mathrm{d}\vec{x} \; \; {}^{\beta}\mathbb{d}_{\vec{l}}^{r} \ .
\end{align}
The integrals of \cref{eqn:D_details_D} are solved exactly by leveraging the properties of the wavelet basis functions, and will vary depending on the types of wavelet coefficients $\lambda$ and $\beta$ as well as the comparison between levels $j$ and $r$. For example, consider calculating the derivative coefficients of a $\lambda = 1$ type of $\mathbb{d}$ coefficient at a specified level $j$ for the case of the mixed derivative (\emph{i.e.}, $\partial^{\alpha_{1} + \alpha_{2}}/\partial x_{1}^{\alpha_{1}} \partial x_{2}^{\alpha_{2}}$) in $3$ spatial dimensions
\begin{align}
\label{eqn:longD}
    {}^{1}\widetilde{\mathbb{d}}_{k_{1} k_{2} k_{3}}^{ \ j} &= \big( \widetilde{\boldsymbol{g}} \cdot {}^{\alpha_{1}} \boldsymbol{\chi} \cdot \boldsymbol{h}^{j} \big)_{k_{1} l_{1}} \ \big( {}^{\alpha_{2}} \boldsymbol{\chi} \cdot \boldsymbol{h}^{j - 1} \big)_{k_{2} l_{2}} \ \big( \boldsymbol{\delta} \cdot \boldsymbol{h}^{j - 1}  \big)_{k_{3} l_{3}} \ {}^{0}\mathbb{d}_{l_{1} l_{2} l_{3}}^{ \ 1} \ldots \nonumber\\
    &+ \sum_{r = 1}^{j - 1} \bigg[ \big( \widetilde{\boldsymbol{g}} \cdot {}^{\alpha_{1}} \boldsymbol{\chi} \cdot \boldsymbol{h}^{j-r} \cdot \boldsymbol{g} \big)_{k_{1} l_{1}} \ \big( {}^{\alpha_{2}} \boldsymbol{\chi} \cdot \boldsymbol{h}^{j-r} \big)_{k_{2} l_{2}} \ \big( \boldsymbol{\delta} \cdot \boldsymbol{h}^{j-r} \big)_{k_{3} l_{3}} \ {}^{1}\mathbb{d}_{l_{1} l_{2} l_{3}}^{ \ r} \ldots \nonumber\\
    &+ \big( \widetilde{\boldsymbol{g}} \cdot {}^{\alpha_{1}} \boldsymbol{\chi} \cdot \boldsymbol{h}^{j-r + 1} \big)_{k_{1} l_{1}} \ \big( {}^{\alpha_{2}} \boldsymbol{\chi} \cdot \boldsymbol{h}^{j-r - 1} \cdot \boldsymbol{g} \big)_{k_{2} l_{2}} \ \big( \boldsymbol{\delta} \cdot \boldsymbol{h}^{j-r} \big)_{k_{3} l_{3}} \ {}^{2}\mathbb{d}_{l_{1} l_{2} l_{3}}^{ \ r} \ldots \nonumber\\
    &+ \big( \widetilde{\boldsymbol{g}} \cdot {}^{\alpha_{1}} \boldsymbol{\chi} \cdot \boldsymbol{h}^{j-r} \cdot \boldsymbol{g} \big)_{k_{1} l_{1}} \ \big( {}^{\alpha_{2}} \boldsymbol{\chi} \cdot \boldsymbol{h}^{j-r - 1} \cdot \boldsymbol{g} \big)_{k_{2} l_{2}} \ \big( \boldsymbol{\delta} \cdot \boldsymbol{h}^{j-r} \big)_{k_{3} l_{3}} \ {}^{3}\mathbb{d}_{l_{1} l_{2} l_{3}}^{ \ r} \ldots \nonumber\\
    &+ \big( \widetilde{\boldsymbol{g}} \cdot {}^{\alpha_{1}} \boldsymbol{\chi} \cdot \boldsymbol{h}^{j-r + 1} \big)_{k_{1} l_{1}} \ \big( {}^{\alpha_{2}} \boldsymbol{\chi} \cdot \boldsymbol{h}^{j-r} \big)_{k_{2} l_{2}} \ \big( \boldsymbol{\delta} \cdot \boldsymbol{h}^{j-r - 1} \cdot \boldsymbol{g} \big)_{k_{3} l_{3}} \ {}^{4}\mathbb{d}_{l_{1} l_{2} l_{3}}^{ \ r} \ldots \nonumber\\
    &+ \big( \widetilde{\boldsymbol{g}} \cdot {}^{\alpha_{1}} \boldsymbol{\chi} \cdot \boldsymbol{h}^{j-r} \cdot \boldsymbol{g} \big)_{k_{1} l_{1}} \ \big( {}^{\alpha_{2}} \boldsymbol{\chi} \cdot \boldsymbol{h}^{j-r} \big)_{k_{2} l_{2}} \ \big( \boldsymbol{\delta} \cdot \boldsymbol{h}^{j-r - 1} \cdot \boldsymbol{g} \big)_{k_{3} l_{3}} \ {}^{5}\mathbb{d}_{l_{1} l_{2} l_{3}}^{ \ r} \ldots \nonumber\\
    &+ \big( \widetilde{\boldsymbol{g}} \cdot {}^{\alpha_{1}} \boldsymbol{\chi} \cdot \boldsymbol{h}^{j-r + 1} \big)_{k_{1} l_{1}} \ \big( {}^{\alpha_{2}} \boldsymbol{\chi} \cdot \boldsymbol{h}^{j-r - 1} \cdot \boldsymbol{g} \big)_{k_{2} l_{2}} \ \big( \boldsymbol{\delta} \cdot \boldsymbol{h}^{j-r - 1} \cdot \boldsymbol{g} \big)_{k_{3} l_{3}} \ {}^{6}\mathbb{d}_{l_{1} l_{2} l_{3}}^{ \ r} \ldots \nonumber\\
    &+ \big( \widetilde{\boldsymbol{g}} \cdot {}^{\alpha_{1}} \boldsymbol{\chi} \cdot \boldsymbol{h}^{j-r} \cdot \boldsymbol{g} \big)_{k_{1} l_{1}} \ \big( {}^{\alpha_{2}} \boldsymbol{\chi} \cdot \boldsymbol{h}^{j-r - 1} \cdot \boldsymbol{g} \big)_{k_{2} l_{2}} \ \big( \boldsymbol{\delta} \cdot \boldsymbol{h}^{j-r - 1} \cdot \boldsymbol{g} \big)_{k_{3} l_{3}} \ {}^{7}\mathbb{d}_{l_{1} l_{2} l_{3}}^{ \ r} \bigg] \ldots \nonumber \\
    &+ \big( \widetilde{\boldsymbol{g}} \cdot {}^{\alpha_{1}} \boldsymbol{\chi} \cdot \boldsymbol{g} \big)_{k_{1} l_{1}} \ \big( {}^{\alpha_{2}} \boldsymbol{\chi} \big)_{k_{2} l_{2}} \ \big( \boldsymbol{\delta} \big)_{k_{3} l_{3}} \ {}^{1}\mathbb{d}_{l_{1} l_{2} l_{3}}^{ \ j} \ldots \nonumber\\
    &+ \big( \widetilde{\boldsymbol{g}} \cdot {}^{\alpha_{1}} \boldsymbol{\chi} \cdot \boldsymbol{h} \big)_{k_{1} l_{1}} \ \big( \widetilde{\boldsymbol{h}} \cdot {}^{\alpha_{2}} \boldsymbol{\chi} \cdot \boldsymbol{g} \big)_{k_{2} l_{2}} \ \big( \boldsymbol{\delta} \big)_{k_{3} l_{3}} \ {}^{2}\mathbb{d}_{l_{1} l_{2} l_{3}}^{ \ j} \ldots \nonumber\\
    &+ \big( \widetilde{\boldsymbol{g}} \cdot {}^{\alpha_{1}} \boldsymbol{\chi} \cdot \boldsymbol{g} \big)_{k_{1} l_{1}} \ \big( \widetilde{\boldsymbol{h}} \cdot {}^{\alpha_{2}} \boldsymbol{\chi} \cdot \boldsymbol{g} \big)_{k_{2} l_{2}} \ \big( \boldsymbol{\delta} \big)_{k_{3} l_{3}} \ {}^{3}\mathbb{d}_{l_{1} l_{2} l_{3}}^{ \ j} \ldots \nonumber\\
    &+ \sum_{r = j + 1}^{j_{\mathrm{max}}} \bigg[ \big( \widetilde{\boldsymbol{g}} \cdot  \widetilde{\boldsymbol{h}}^{r - j} \cdot {}^{\alpha_{1}} \boldsymbol{\chi} \cdot \boldsymbol{g} \big)_{k_{1} l_{1}} \ \big( \widetilde{\boldsymbol{h}}^{r - j} \cdot {}^{\alpha_{2}} \boldsymbol{\chi} \big)_{k_{2} l_{2}} \ \big( \widetilde{\boldsymbol{h}}^{r - j} \cdot \boldsymbol{\delta} \big)_{k_{3} l_{3}} \ {}^{1}\mathbb{d}_{l_{1} l_{2} l_{3}}^{ \ r} \ldots \nonumber\\
    &+ \big( \widetilde{\boldsymbol{g}} \cdot  \widetilde{\boldsymbol{h}}^{r - j - 1} \cdot {}^{\alpha_{1}} \boldsymbol{\chi} \big)_{k_{1} l_{1}} \ \big( \widetilde{\boldsymbol{h}}^{r - j + 1} \cdot {}^{\alpha_{2}} \boldsymbol{\chi} \cdot \boldsymbol{g} \big)_{k_{2} l_{2}} \ \big( \widetilde{\boldsymbol{h}}^{r - j} \cdot \boldsymbol{\delta} \big)_{k_{3} l_{3}} \ {}^{2}\mathbb{d}_{l_{1} l_{2} l_{3}}^{ \ r} \ldots \nonumber\\
    &+ \big( \widetilde{\boldsymbol{g}} \cdot  \widetilde{\boldsymbol{h}}^{r - j} \cdot {}^{\alpha_{1}} \boldsymbol{\chi} \cdot \boldsymbol{g} \big)_{k_{1} l_{1}} \ \big( \widetilde{\boldsymbol{h}}^{r - j + 1} \cdot {}^{\alpha_{2}} \boldsymbol{\chi} \cdot \boldsymbol{g} \big)_{k_{2} l_{2}} \ \big( \widetilde{\boldsymbol{h}}^{r - j} \cdot \boldsymbol{\delta} \big)_{k_{3} l_{3}} \ {}^{3}\mathbb{d}_{l_{1} l_{2} l_{3}}^{ \ r} \bigg] .
\end{align}
In \cref{eqn:longD}, $\boldsymbol{\delta}$ is the identity matrix, the $\boldsymbol{h}, \boldsymbol{g}, \widetilde{\boldsymbol{h}}, \widetilde{\boldsymbol{g}}$ matrices are defined by the wavelet filter coefficients (also shown in \cite{Harnish2018} to be subsets of the the $\mathbbb{F}$ and $\mathbbb{B}$ operators). The exponent notation implies repeated matrix multiplication, and the $\boldsymbol{\chi}$ matrix is defined by solving an eigenvector problem as described in \cite{Harnish2018,Harnish2021}.

\clearpage

\section{MRWT scaling performance}
\label{sec:MRWTscaling}
\begin{figure}[!htb]
\centering
    \includegraphics[width=0.5\textwidth]{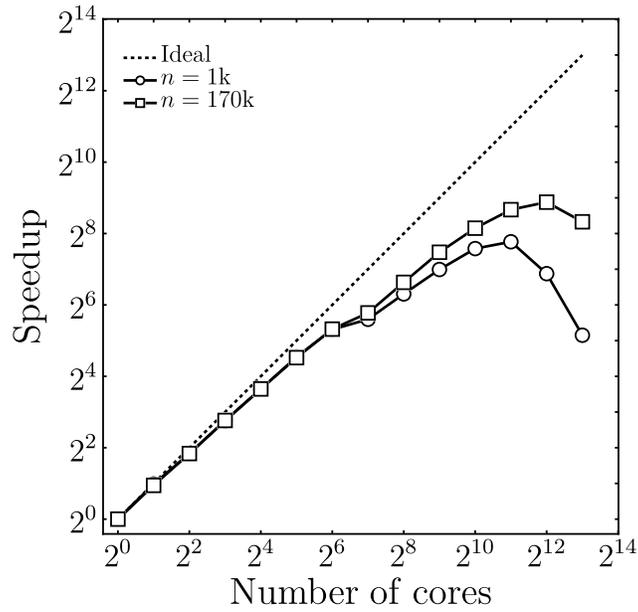}
\caption{Strong scaling performance of hybrid OpenMP/MPI parallel implementation of MRWT.\label{fig:scaling_MRWT}}
\end{figure}

Strong scaling of the MRWT has been evaluated at both early and late
times in the simulation of damage nucleation and propagation from
\cref{sec:damage_results}. The time required to complete an embedded
$\mathcal{O}(\Delta t^{4})$ and $\mathcal{O}(\Delta t^{5})$ explicit
Runge-Kutta integration step (\emph{i.e.}, $n \rightarrow n + 1$) was
measured using Indiana University’s BigRed200, an HPE Cray EX
supercomputer where each node is equipped with $256$ GB of memory and
two $64$-core AMD EPYC $7742$ processors \cite{IUpti}. Results bind
one MPI process to each socket and continue to use OpenMP within each
process. Simultaneous multithreading is not evaluated. The associated
speedup is defined with respect to the run-time using only a single
core. As shown in \cref{fig:scaling_MRWT}, the early time-step
(\emph{i.e.}, $n = 1,000$ time-steps, $t \approx 0.03 \ \mu$s, and
$\mathcal{O}(100$k) points) strong scales from $225$ s to $1.03$ s per
step at $2,048$ cores, at which point communication dominates. The
later time-step (\emph{i.e.}, $n = 170,000$ time-steps, $t \approx
1.75 \ \mu$s, and $\mathcal{O}(840$k) points) scales from $2,542$ s to
$5.4$ s per step at $4,096$ cores.

\clearpage

\section{Taylor-Sedov blast wave}
\label{sec:sedov}
\begin{figure}[!htb]
\begin{tabular*}{\textwidth}{@{} c @{\extracolsep{\fill}} c @{}}
    \includegraphics[width=0.48\textwidth]{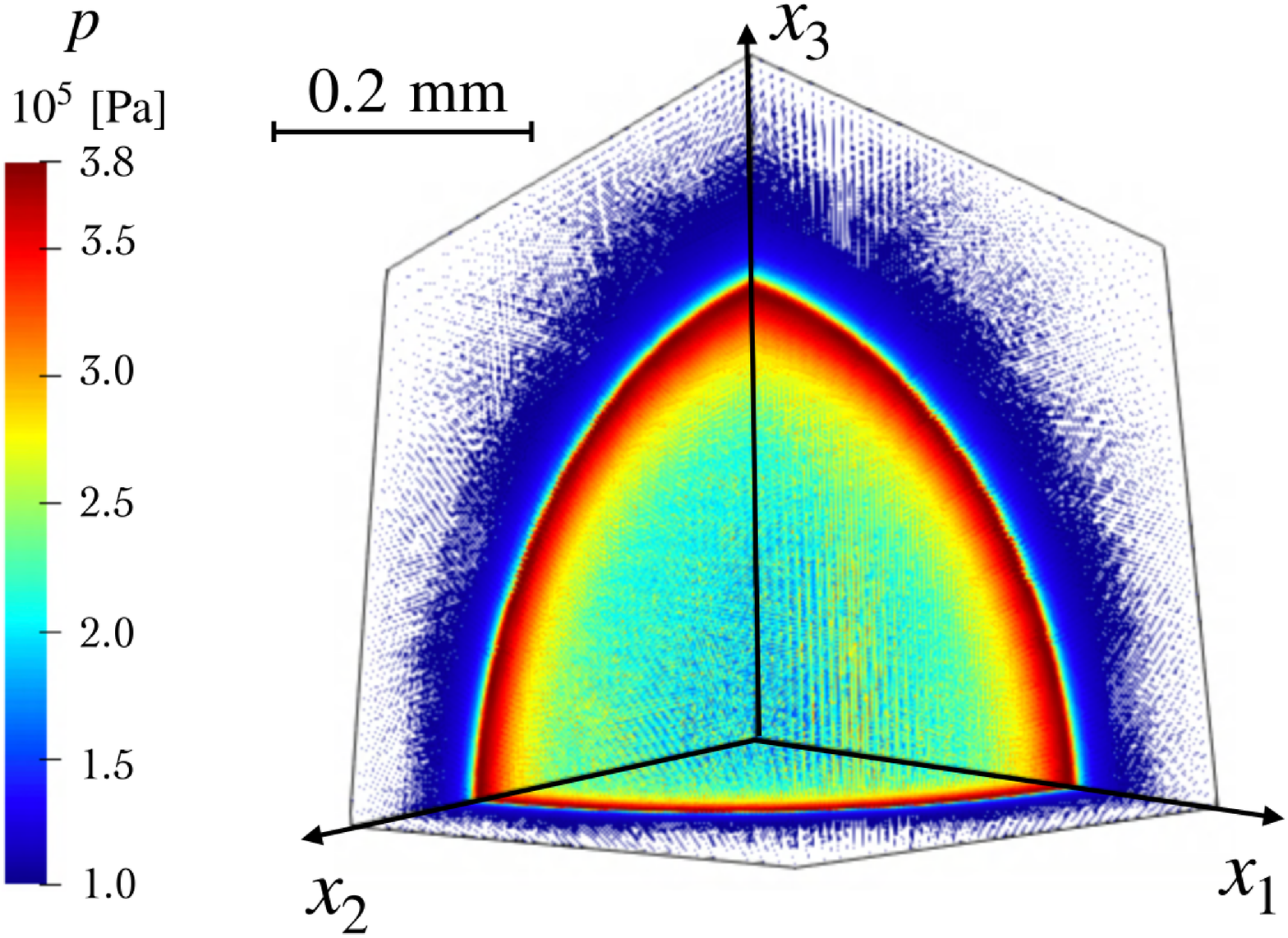}&
    \includegraphics[width=0.48\textwidth]{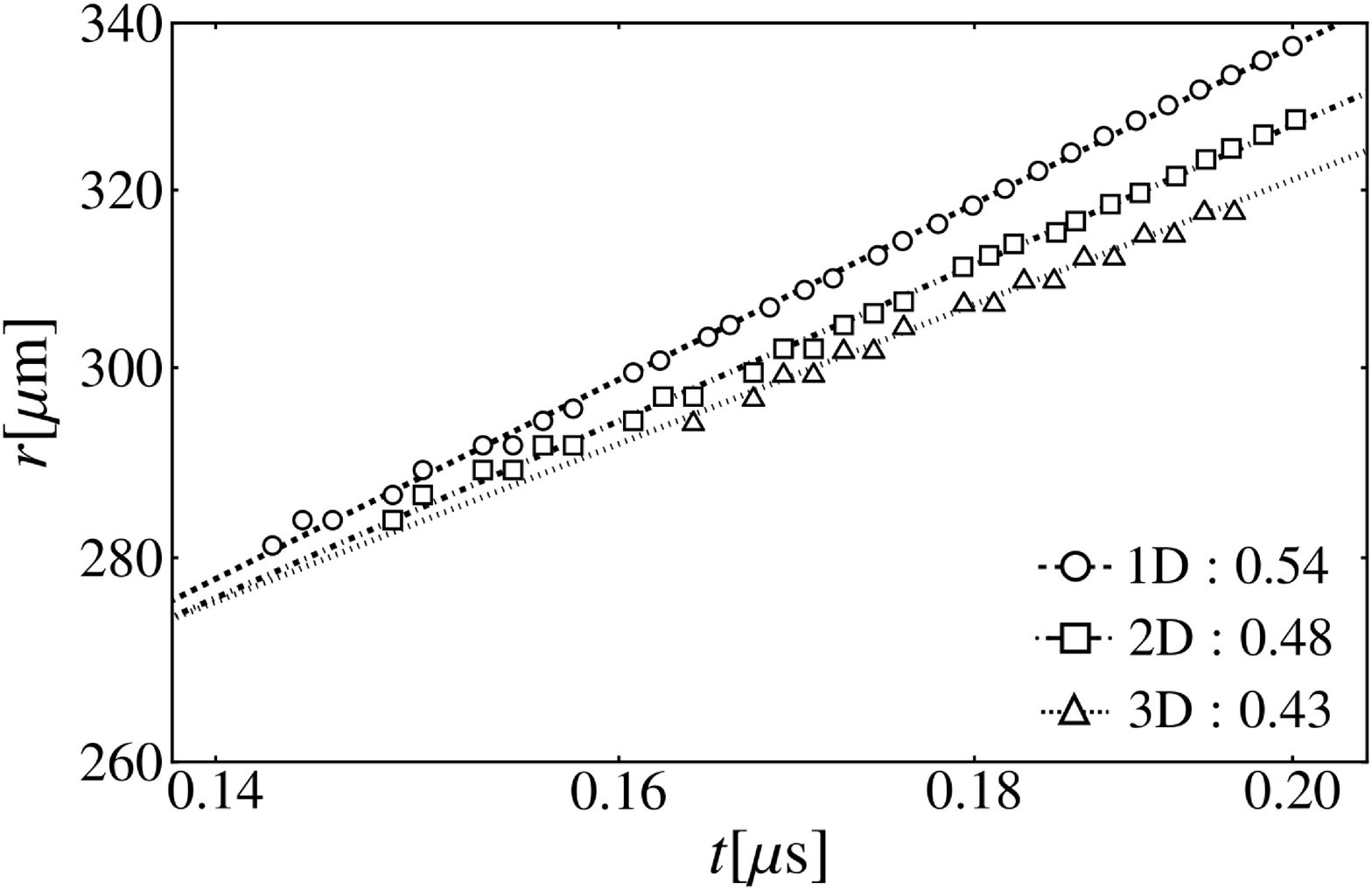}\\
    (a) & (b)
\end{tabular*}
\caption{(a) MRWT simulation of the pressure in a $3$D Taylor-Sedov
  blast wave at $t = 0.2 \ \mu$s. (b) MRWT computed radius of the
  shock wave as a function of time, lines are fitted to the discrete
  data and the corresponding slope should match the exponent in
  \cref{eqn:sedov_radius}. Specifically: $1$D $\rightarrow
  0.\overline{66}$, $2$D $\rightarrow 0.5$, and $3$D $\rightarrow
  0.4$. \label{fig:sedov_analysis}} 
\end{figure}

Details of the Taylor-Sedov problem can be found in previous
publications \cite{Harnish2021,HarnishSpecialIssue}. In summary, the
Navier-Stokes equations are solved with Fourier's law of heat
conduction and model of a Newtonian fluid as a calorically perfect
ideal gas with constant material properties matching those of dry
air. The initial condition is a pressure spike in an otherwise
undisturbed medium and over time this potential energy is converted
into kinetic energy forming a shock wave at a radius $r(t)$
theoretically propagating outward such that \cite{sedov} 
\begin{align}
\label{eqn:sedov_radius}
    r(t) \propto t^{\frac{2}{2 + N}},
\end{align}
where $N$ is the number of spatial dimensions. In our simulations, the
radius of the shock wave is defined as the location of the most
extreme pressure gradient. \Cref{fig:sedov_analysis} plots the
pressure field from the $3$D simulation using MRWT along with the
radius calculations from $1$D, $2$D, and $3$D (see
\cref{eqn:sedov_radius}). The slopes match well the theoretical
predictions indicating that our wavelet method provides accurate
solutions to systems of nonlinear PDEs in multiple dimensions. Note
that \Cref{fig:sedov_analysis}(a) was obtained with $p = 6$ and
threshold parameters
\begin{align}
    \varepsilon_{\rho} &= 10^{-2} \ \| \ \rho \ \|_{\infty},
    &
    \varepsilon_{v_{i}} &= 10^{-2} \ \| \  v_{i} \ \|_{\infty},
    &
    \varepsilon_{e} &= 10^{-2} \ \| \ e \ \|_{\infty},
\end{align}
where $\rho$ is the density, $v_{i}$ is the $i^{\mathrm{th}}$ component of the velocity, and $e$ is the specific internal energy.

\clearpage

\bibliography{references}

\begin{thebibliography}{127}
\expandafter\ifx\csname natexlab\endcsname\relax\def\natexlab#1{#1}\fi
\providecommand{\url}[1]{\texttt{#1}}
\providecommand{\href}[2]{#2}
\providecommand{\path}[1]{#1}
\providecommand{\DOIprefix}{doi:}
\providecommand{\ArXivprefix}{arXiv:}
\providecommand{\URLprefix}{URL: }
\providecommand{\Pubmedprefix}{pmid:}
\providecommand{\doi}[1]{\href{http://dx.doi.org/#1}{\path{#1}}}
\providecommand{\Pubmed}[1]{\href{pmid:#1}{\path{#1}}}
\providecommand{\bibinfo}[2]{#2}
\ifx\xfnm\relax \def\xfnm[#1]{\unskip,\space#1}\fi
\bibitem[{Irwin(1957)}]{historyIrwin}
\bibinfo{author}{G.~R. Irwin},
\newblock \bibinfo{title}{{Analysis of stresses and strains near the end of a
  crack traversing a plate}},
\newblock \bibinfo{journal}{Journal of Applied Mechanics} \bibinfo{volume}{24}
  (\bibinfo{year}{1957}) \bibinfo{pages}{361--364}.
\bibitem[{Dugdale(1960)}]{historyDugdale}
\bibinfo{author}{D.~S. Dugdale},
\newblock \bibinfo{title}{{Yielding of steel sheets containing slits}},
\newblock \bibinfo{journal}{Journal of the Mechanics and Physics of Solids}
  \bibinfo{volume}{8} (\bibinfo{year}{1960}) \bibinfo{pages}{100--104}.
\bibitem[{Willis(1967)}]{historyWillis}
\bibinfo{author}{J.~R. Willis},
\newblock \bibinfo{title}{{A comparison of the fracture criteria of Griffith
  and Barenblatt}},
\newblock \bibinfo{journal}{Journal of the Mechanics and Physics of Solids}
  \bibinfo{volume}{15} (\bibinfo{year}{1967}) \bibinfo{pages}{151--162}.
\bibitem[{Spanoudakis and Young(1984)}]{historySpanoudakis}
\bibinfo{author}{J.~Spanoudakis}, \bibinfo{author}{R.~J. Young},
\newblock \bibinfo{title}{{Crack propagation in a glass particle-filled epoxy
  resin part 1. effect of particle volume fraction and size}},
\newblock \bibinfo{journal}{Journal of Materials Science} \bibinfo{volume}{19}
  (\bibinfo{year}{1984}) \bibinfo{pages}{473--486}.
\bibitem[{Ba{\v{z}}ant(1986)}]{historyBazant2}
\bibinfo{author}{Z.~P. Ba{\v{z}}ant},
\newblock \bibinfo{title}{{Mechanics of distributed cracking}},
\newblock \bibinfo{journal}{Applied Mechanics Reviews} \bibinfo{volume}{39}
  (\bibinfo{year}{1986}) \bibinfo{pages}{675--705}.
\bibitem[{Ju(1989)}]{historyJu}
\bibinfo{author}{J.~W. Ju},
\newblock \bibinfo{title}{{Energy-based coupled elastoplastic damage models at
  finite strains}},
\newblock \bibinfo{journal}{Journal of Engineering Mechanics}
  \bibinfo{volume}{115} (\bibinfo{year}{1989}) \bibinfo{pages}{2507--2525}.
\bibitem[{Needleman(1990)}]{historyNeedleman}
\bibinfo{author}{A.~Needleman},
\newblock \bibinfo{title}{{An analysis of tensile decohesion along an
  interface}},
\newblock \bibinfo{journal}{Journal of the Mechanics and Physics of Solids}
  \bibinfo{volume}{38} (\bibinfo{year}{1990}) \bibinfo{pages}{289--324}.
\bibitem[{Fish and Yu(2001)}]{historyFish}
\bibinfo{author}{J.~Fish}, \bibinfo{author}{Q.~Yu},
\newblock \bibinfo{title}{{Multiscale damage modelling for composite materials:
  Theory and computational framework}},
\newblock \bibinfo{journal}{International Journal for Numerical Methods in
  Engineering} \bibinfo{volume}{52} (\bibinfo{year}{2001})
  \bibinfo{pages}{161--191}.
\bibitem[{Ba{\v{z}}ant and Jir{\'{a}}sek(2002)}]{historyBazant}
\bibinfo{author}{Z.~P. Ba{\v{z}}ant}, \bibinfo{author}{M.~Jir{\'{a}}sek},
\newblock \bibinfo{title}{{Nonlocal integral formulations of plasticity and
  damage: survey of progress}},
\newblock \bibinfo{journal}{Journal of Engineering Mechanics}
  \bibinfo{volume}{128} (\bibinfo{year}{2002}) \bibinfo{pages}{1119--1149}.
\bibitem[{Buehler et~al.(2003)Buehler, Abraham, and Gao}]{historyBuehler}
\bibinfo{author}{M.~J. Buehler}, \bibinfo{author}{F.~F. Abraham},
  \bibinfo{author}{H.~Gao},
\newblock \bibinfo{title}{{Hyperelasticity governs dynamic fracture at a
  critical length scale}},
\newblock \bibinfo{journal}{Nature} \bibinfo{volume}{426}
  (\bibinfo{year}{2003}) \bibinfo{pages}{141--146}.
\bibitem[{Voyiadjis et~al.(2004)Voyiadjis, Abu Al-Rub, and
  Palazotto}]{historyVoyiadjis}
\bibinfo{author}{G.~Z. Voyiadjis}, \bibinfo{author}{R.~K. Abu Al-Rub},
  \bibinfo{author}{A.~N. Palazotto},
\newblock \bibinfo{title}{{Thermodynamic framework for coupling of non-local
  viscoplasticity and non-local anisotropic viscodamage for dynamic
  localization problems using gradient theory}},
\newblock \bibinfo{journal}{International Journal of Plasticity}
  \bibinfo{volume}{20} (\bibinfo{year}{2004}) \bibinfo{pages}{981--1038}.
\bibitem[{Kitey and Tippur(2005)}]{historyKitey}
\bibinfo{author}{R.~Kitey}, \bibinfo{author}{H.~V. Tippur},
\newblock \bibinfo{title}{{Role of particle size and filler-matrix adhesion on
  dynamic fracture of glass-filled epoxy part I. Macromeasurements}},
\newblock \bibinfo{journal}{Acta Materialia} \bibinfo{volume}{53}
  (\bibinfo{year}{2005}) \bibinfo{pages}{1153--1165}.
\bibitem[{Park et~al.(2009)Park, Paulino, and Roesler}]{historyPark}
\bibinfo{author}{K.~Park}, \bibinfo{author}{G.~H. Paulino},
  \bibinfo{author}{J.~R. Roesler},
\newblock \bibinfo{title}{{A unified potential-based cohesive model of
  mixed-mode fracture}},
\newblock \bibinfo{journal}{Journal of the Mechanics and Physics of Solids}
  \bibinfo{volume}{57} (\bibinfo{year}{2009}) \bibinfo{pages}{891--908}.
\bibitem[{Huynh and Belytschko(2009)}]{historyHuynh}
\bibinfo{author}{D.~B.~P. Huynh}, \bibinfo{author}{T.~Belytschko},
\newblock \bibinfo{title}{{The extended finite element method for fracture in
  composite materials}},
\newblock \bibinfo{journal}{International Journal for Numerical Methods in
  Engineering} \bibinfo{volume}{77} (\bibinfo{year}{2009})
  \bibinfo{pages}{214--239}.
\bibitem[{Kim et~al.(2010)Kim, Pereira, and Duarte}]{historyKim}
\bibinfo{author}{D.~J. Kim}, \bibinfo{author}{J.~P. Pereira},
  \bibinfo{author}{C.~A. Duarte},
\newblock \bibinfo{title}{{Analysis of three-dimensional fracture mechanics
  problems: A two-scale approach using coarse-generalized FEM meshes}},
\newblock \bibinfo{journal}{International Journal for Numerical Methods in
  Engineering} \bibinfo{volume}{81} (\bibinfo{year}{2010})
  \bibinfo{pages}{335--365}.
\bibitem[{Hankinson(1921)}]{Hankinsonphenom}
\bibinfo{author}{R.~Hankinson}, \bibinfo{title}{{Investigation of Crushing
  Strength of Spruce at Varying Angles of Grain}}, \bibinfo{type}{Technical
  Report}, Air Service Information Circular, \bibinfo{year}{1921}.
\bibitem[{Hill(1948)}]{Hillphenom}
\bibinfo{author}{R.~Hill},
\newblock \bibinfo{title}{{A theory of the yielding and plastic flow of
  anisotropic metals}},
\newblock \bibinfo{journal}{Proceedings of the Royal Society of London}
  \bibinfo{volume}{193} (\bibinfo{year}{1948}) \bibinfo{pages}{281--297}.
\bibitem[{Drucker and Prager(1952)}]{DPphenom}
\bibinfo{author}{D.~C. Drucker}, \bibinfo{author}{W.~Prager},
\newblock \bibinfo{title}{{Soil mechanics and plastic analysis or limit
  design}},
\newblock \bibinfo{journal}{Quarterly of Applied Mathematics}
  \bibinfo{volume}{10} (\bibinfo{year}{1952}) \bibinfo{pages}{157--165}.
\bibitem[{Tsai and Wu(1971)}]{TWphenom}
\bibinfo{author}{S.~W. Tsai}, \bibinfo{author}{E.~M. Wu},
\newblock \bibinfo{title}{{A general theory of strength for anisotropic
  materials}},
\newblock \bibinfo{journal}{Journal of Composite Materials} \bibinfo{volume}{5}
  (\bibinfo{year}{1971}) \bibinfo{pages}{58--80}.
\bibitem[{Hoek and Brown(1980)}]{HBphenom}
\bibinfo{author}{E.~Hoek}, \bibinfo{author}{E.~Brown},
\newblock \bibinfo{title}{{Empirical strength criterion for rock masses}},
\newblock \bibinfo{journal}{Journal of the Geotechnical Engineering Division}
  \bibinfo{volume}{106} (\bibinfo{year}{1980}) \bibinfo{pages}{1013--1035}.
\bibitem[{Johnson and Holmquist(1994)}]{JHphenom}
\bibinfo{author}{G.~R. Johnson}, \bibinfo{author}{T.~J. Holmquist},
\newblock \bibinfo{title}{{An improved computational constitutive model for
  brittle materials}},
\newblock in: \bibinfo{booktitle}{AIP Conference Proceedings},
  \bibinfo{year}{1994}, pp. \bibinfo{pages}{981--984}.
  \DOIprefix\doi{10.1063/1.46199}.
\bibitem[{Staat(2021)}]{MCphenom}
\bibinfo{author}{M.~Staat},
\newblock \bibinfo{title}{{An extension strain type Mohr–Coulomb criterion}},
\newblock \bibinfo{journal}{Rock Mechanics and Rock Engineering}
  \bibinfo{volume}{54} (\bibinfo{year}{2021}) \bibinfo{pages}{6207--6233}.
\bibitem[{Griffith(1921)}]{Griffith}
\bibinfo{author}{A.~A. Griffith},
\newblock \bibinfo{title}{{The phenomena of rupture and flow in solids}},
\newblock \bibinfo{journal}{Philosophical Transactions of the Royal Society of
  London} \bibinfo{volume}{221} (\bibinfo{year}{1921})
  \bibinfo{pages}{163--198}.
\bibitem[{Barenblatt(1961)}]{historyBarenblatt}
\bibinfo{author}{G.~I. Barenblatt},
\newblock \bibinfo{title}{{The mathematical theory of equilibrium cracks formed
  in brittle fracture}},
\newblock \bibinfo{journal}{Zhurnal Prikladnoy Mekhaniki i Tekhnicheskoy
  Fiziki} \bibinfo{volume}{4} (\bibinfo{year}{1961}) \bibinfo{pages}{3--56}.
\bibitem[{Freund(1990)}]{Freund}
\bibinfo{author}{L.~B. Freund}, \bibinfo{title}{{Dynamic Fracture Mechanics}},
  \bibinfo{publisher}{Cambridge University Press},
  \bibinfo{address}{Cambridge}, \bibinfo{year}{1990}.
  \DOIprefix\doi{10.1017/CBO9780511546761}.
\bibitem[{Krajcinovic(1995)}]{historyKrajcinovic}
\bibinfo{author}{D.~Krajcinovic},
\newblock \bibinfo{title}{{Continuum damage mechanics: when and how?}},
\newblock \bibinfo{journal}{International Journal of Damage Mechanics}
  \bibinfo{volume}{4} (\bibinfo{year}{1995}) \bibinfo{pages}{217--229}.
\bibitem[{Lemaitre(1971)}]{fatigue1971}
\bibinfo{author}{J.~Lemaitre},
\newblock \bibinfo{title}{{Evaluation of dissipation and damage in metals
  submitted to dynamic loading}},
\newblock in: \bibinfo{booktitle}{International Conference on the Mechanical
  Behavior of Materials}, \bibinfo{year}{1971}, pp. \bibinfo{pages}{1--8}.
\bibitem[{Leckie and Hayhurst(1974)}]{creep1974}
\bibinfo{author}{F.~A. Leckie}, \bibinfo{author}{D.~R. Hayhurst},
\newblock \bibinfo{title}{{Creep rupture of structures}},
\newblock \bibinfo{journal}{Proceedings of the Royal Society of London. Series
  A, Mathematical and Physical Sciences} \bibinfo{volume}{340}
  (\bibinfo{year}{1974}) \bibinfo{pages}{323--347}.
\bibitem[{Lemaitre(1985)}]{ductile1985}
\bibinfo{author}{J.~Lemaitre},
\newblock \bibinfo{title}{{A continuous damage mechanics model for ductile
  fracture}},
\newblock \bibinfo{journal}{Journal of Engineering Materials and Technology}
  \bibinfo{volume}{107} (\bibinfo{year}{1985}) \bibinfo{pages}{83--89}.
\bibitem[{Dragon and Chihab(1985)}]{ductile1985b}
\bibinfo{author}{A.~Dragon}, \bibinfo{author}{A.~Chihab},
\newblock \bibinfo{title}{{On finite damage: ductile fracture-damage
  evolution}},
\newblock \bibinfo{journal}{Mechanics of Materials} \bibinfo{volume}{4}
  (\bibinfo{year}{1985}) \bibinfo{pages}{95--106}.
\bibitem[{Dragon(1985)}]{ductile1985c}
\bibinfo{author}{A.~Dragon},
\newblock \bibinfo{title}{{Plasticity and ductile fracture damage: study of
  void growth in metals}},
\newblock \bibinfo{journal}{Engineering Fracture Mechanics}
  \bibinfo{volume}{21} (\bibinfo{year}{1985}) \bibinfo{pages}{875--885}.
\bibitem[{Krajcinovic and Fonseka(1981)}]{brittle1981}
\bibinfo{author}{D.~Krajcinovic}, \bibinfo{author}{G.~U. Fonseka},
\newblock \bibinfo{title}{{The continuous damage theory of brittle materials}},
\newblock \bibinfo{journal}{Journal of Applied Mechanics} \bibinfo{volume}{48}
  (\bibinfo{year}{1981}) \bibinfo{pages}{809--815}.
\bibitem[{Krajcinovic(1983)}]{brittle1983}
\bibinfo{author}{D.~Krajcinovic},
\newblock \bibinfo{title}{{Constitutive equations for damaging materials}},
\newblock \bibinfo{journal}{Journal of Applied Mechanics} \bibinfo{volume}{50}
  (\bibinfo{year}{1983}) \bibinfo{pages}{355--360}.
\bibitem[{Marigo(1985)}]{brittle1985}
\bibinfo{author}{J.~J. Marigo},
\newblock \bibinfo{title}{{Modelling of brittle and fatigue damage for elastic
  material by growth of microvoids}},
\newblock \bibinfo{journal}{Engineering Fracture Mechanics}
  \bibinfo{volume}{21} (\bibinfo{year}{1985}) \bibinfo{pages}{861--874}.
\bibitem[{Lemaitre(1984)}]{isotropic1984}
\bibinfo{author}{J.~Lemaitre},
\newblock \bibinfo{title}{{How to use damage mechanics}},
\newblock \bibinfo{journal}{Nuclear Engineering and Design}
  \bibinfo{volume}{80} (\bibinfo{year}{1984}) \bibinfo{pages}{233--245}.
\bibitem[{Kachanov(1980)}]{anisotropic1980}
\bibinfo{author}{M.~Kachanov},
\newblock \bibinfo{title}{{Continuum model of medium with cracks}},
\newblock \bibinfo{journal}{Journal of the Engineering Mechanics Division}
  \bibinfo{volume}{106} (\bibinfo{year}{1980}) \bibinfo{pages}{1039--1051}.
\bibitem[{Chaboche(1981)}]{anisotropic1981}
\bibinfo{author}{J.-L. Chaboche},
\newblock \bibinfo{title}{{Continuous damage mechanics -- a tool to describe
  phenomena before crack initiation}},
\newblock \bibinfo{journal}{Nuclear Engineering and Design}
  \bibinfo{volume}{64} (\bibinfo{year}{1981}) \bibinfo{pages}{233--247}.
\bibitem[{Kamrin et~al.(2012)Kamrin, Rycroft, and Nave}]{RMT}
\bibinfo{author}{K.~Kamrin}, \bibinfo{author}{C.~H. Rycroft},
  \bibinfo{author}{J.~C. Nave},
\newblock \bibinfo{title}{{Reference map technique for finite-strain elasticity
  and fluid-solid interaction}},
\newblock \bibinfo{journal}{Journal of the Mechanics and Physics of Solids}
  \bibinfo{volume}{60} (\bibinfo{year}{2012}) \bibinfo{pages}{1952--1969}.
\bibitem[{Matou{\v{s}} et~al.(2008)Matou{\v{s}}, Kulkarni, and
  Geubelle}]{matous2008}
\bibinfo{author}{K.~Matou{\v{s}}}, \bibinfo{author}{M.~G. Kulkarni},
  \bibinfo{author}{P.~H. Geubelle},
\newblock \bibinfo{title}{{Multiscale cohesive failure modeling of
  heterogeneous adhesives}},
\newblock \bibinfo{journal}{Journal of the Mechanics and Physics of Solids}
  \bibinfo{volume}{56} (\bibinfo{year}{2008}) \bibinfo{pages}{1511--1533}.
\bibitem[{Kulkarni et~al.(2010)Kulkarni, Matou{\v{s}}, and
  Geubelle}]{matous2010}
\bibinfo{author}{M.~G. Kulkarni}, \bibinfo{author}{K.~Matou{\v{s}}},
  \bibinfo{author}{P.~H. Geubelle},
\newblock \bibinfo{title}{{Coupled multi-scale cohesive modeling of failure in
  heterogeneous adhesives}},
\newblock \bibinfo{journal}{International Journal for Numerical Methods in
  Engineering} \bibinfo{volume}{84} (\bibinfo{year}{2010})
  \bibinfo{pages}{916--946}.
\bibitem[{Mosby and Matou{\v{s}}(2015)}]{matous2015}
\bibinfo{author}{M.~Mosby}, \bibinfo{author}{K.~Matou{\v{s}}},
\newblock \bibinfo{title}{{On mechanics and material length scales of failure
  in heterogeneous interfaces using a finite strain high performance solver}},
\newblock \bibinfo{journal}{Modelling and Simulation in Materials Science and
  Engineering} \bibinfo{volume}{23} (\bibinfo{year}{2015}).
\bibitem[{Simo and Ju(1987{\natexlab{a}})}]{historySimoJu}
\bibinfo{author}{J.~C. Simo}, \bibinfo{author}{J.~W. Ju},
\newblock \bibinfo{title}{{Strain- and stress-based continuum damage models --
  I. Formulation}},
\newblock \bibinfo{journal}{Int. J. Solids Structures} \bibinfo{volume}{23}
  (\bibinfo{year}{1987}{\natexlab{a}}) \bibinfo{pages}{821--840}.
\bibitem[{Simo and Ju(1987{\natexlab{b}})}]{historySimoJu2}
\bibinfo{author}{J.~C. Simo}, \bibinfo{author}{J.~W. Ju},
\newblock \bibinfo{title}{{Strain- and stress-based continuum damage models --
  II. Computational aspects}},
\newblock \bibinfo{journal}{Int. J. Solids Structures} \bibinfo{volume}{23}
  (\bibinfo{year}{1987}{\natexlab{b}}) \bibinfo{pages}{841--869}.
\bibitem[{Xue et~al.(2010)Xue, Pontin, Zok, and Hutchinson}]{zoneXue}
\bibinfo{author}{Z.~Xue}, \bibinfo{author}{M.~G. Pontin},
  \bibinfo{author}{F.~W. Zok}, \bibinfo{author}{J.~W. Hutchinson},
\newblock \bibinfo{title}{{Calibration procedures for a computational model of
  ductile fracture}},
\newblock \bibinfo{journal}{Engineering Fracture Mechanics}
  \bibinfo{volume}{77} (\bibinfo{year}{2010}) \bibinfo{pages}{492--509}.
\bibitem[{Volokh(2012)}]{zoneSteel}
\bibinfo{author}{K.~Y. Volokh},
\newblock \bibinfo{title}{{Characteristic length of damage localization in
  steel}},
\newblock \bibinfo{journal}{Engineering Fracture Mechanics}
  \bibinfo{volume}{94} (\bibinfo{year}{2012}) \bibinfo{pages}{85--86}.
\bibitem[{Volokh(2013)}]{zoneConcrete}
\bibinfo{author}{K.~Y. Volokh},
\newblock \bibinfo{title}{{Characteristic length of damage localization in
  concrete}},
\newblock \bibinfo{journal}{Mechanics Research Communications}
  \bibinfo{volume}{51} (\bibinfo{year}{2013}) \bibinfo{pages}{29--31}.
\bibitem[{Volokh(2011)}]{zoneRubber}
\bibinfo{author}{K.~Y. Volokh},
\newblock \bibinfo{title}{{Characteristic length of damage localization in
  rubber}},
\newblock \bibinfo{journal}{International Journal of Fracture}
  \bibinfo{volume}{168} (\bibinfo{year}{2011}) \bibinfo{pages}{113--116}.
\bibitem[{Krajcinovic and Fanella(1986)}]{damageScales}
\bibinfo{author}{D.~Krajcinovic}, \bibinfo{author}{D.~Fanella},
\newblock \bibinfo{title}{{A micromechanical damage model for concrete}},
\newblock \bibinfo{journal}{Engineering Fracture Mechanics}
  \bibinfo{volume}{25} (\bibinfo{year}{1986}) \bibinfo{pages}{585--596}.
\bibitem[{Pearson and Yee(1991)}]{expPearson}
\bibinfo{author}{R.~A. Pearson}, \bibinfo{author}{A.~F. Yee},
\newblock \bibinfo{title}{{Influence of particle size and particle size
  distribution on toughening mechanisms in rubber-modified epoxies}},
\newblock \bibinfo{journal}{Journal of Materials Science} \bibinfo{volume}{26}
  (\bibinfo{year}{1991}) \bibinfo{pages}{3828--3844}.
\bibitem[{Dekkers and Heikens(1983)}]{expDekkers}
\bibinfo{author}{M.~E.~J. Dekkers}, \bibinfo{author}{D.~Heikens},
\newblock \bibinfo{title}{{The effect of interfacial adhesion on the tensile
  behavior of polystyrene-glass-bead composites}},
\newblock \bibinfo{journal}{Journal of Applied Polymer Science}
  \bibinfo{volume}{28} (\bibinfo{year}{1983}) \bibinfo{pages}{3809--3815}.
\bibitem[{Huang et~al.(2004)Huang, Yuan, Jiang, An, Jiang, and Li}]{expHuang}
\bibinfo{author}{L.~Huang}, \bibinfo{author}{Q.~Yuan},
  \bibinfo{author}{W.~Jiang}, \bibinfo{author}{L.~An},
  \bibinfo{author}{S.~Jiang}, \bibinfo{author}{R.~K. Li},
\newblock \bibinfo{title}{{Mechanical and thermal properties of glass
  bead-filled nylon-6}},
\newblock \bibinfo{journal}{Journal of Applied Polymer Science}
  \bibinfo{volume}{94} (\bibinfo{year}{2004}) \bibinfo{pages}{1885--1890}.
\bibitem[{Kinlock and Hunston(1987)}]{expKinlock}
\bibinfo{author}{A.~J. Kinlock}, \bibinfo{author}{D.~L. Hunston},
\newblock \bibinfo{title}{{Effect of volume fraction of dispersed rubbery phase
  on the toughness of rubber-toughened epoxy polymers}},
\newblock \bibinfo{journal}{Journal of Materials Science Letters}
  \bibinfo{volume}{6} (\bibinfo{year}{1987}) \bibinfo{pages}{131--139}.
\bibitem[{Radford(1971)}]{expRadford}
\bibinfo{author}{K.~C. Radford},
\newblock \bibinfo{title}{{The mechanical properties of an epoxy resin with a
  second phase dispersion}},
\newblock \bibinfo{journal}{Journal of Materials Science} \bibinfo{volume}{6}
  (\bibinfo{year}{1971}) \bibinfo{pages}{1286--1291}.
\bibitem[{Nakamura et~al.(1992)Nakamura, Yamaguchi, {Masayoshi Okubo}, and
  Matsumoto}]{expNakamura}
\bibinfo{author}{Y.~Nakamura}, \bibinfo{author}{M.~Yamaguchi},
  \bibinfo{author}{{Masayoshi Okubo}}, \bibinfo{author}{T.~Matsumoto},
\newblock \bibinfo{title}{{Effects of particle size on mechanical and impact
  properties of epoxy resin filled with spherical silica}},
\newblock \bibinfo{journal}{Journal of Applied Polymer Science}
  \bibinfo{volume}{45} (\bibinfo{year}{1992}) \bibinfo{pages}{1281--1289}.
\bibitem[{Singh et~al.(2002)Singh, Zhang, and Chan}]{expSingh}
\bibinfo{author}{R.~P. Singh}, \bibinfo{author}{M.~Zhang},
  \bibinfo{author}{D.~Chan},
\newblock \bibinfo{title}{{Toughening of a brittle thermosetting polymer:
  Effects of reinforcement particle size and volume fraction}},
\newblock \bibinfo{journal}{Journal of Materials Science} \bibinfo{volume}{37}
  (\bibinfo{year}{2002}) \bibinfo{pages}{781--788}.
\bibitem[{Arag{\'{o}}n et~al.(2013)Arag{\'{o}}n, Soghrati, and
  Geubelle}]{compChallenges2}
\bibinfo{author}{A.~M. Arag{\'{o}}n}, \bibinfo{author}{S.~Soghrati},
  \bibinfo{author}{P.~H. Geubelle},
\newblock \bibinfo{title}{{Effect of in-plane deformation on the cohesive
  failure of heterogeneous adhesives}},
\newblock \bibinfo{journal}{Journal of the Mechanics and Physics of Solids}
  \bibinfo{volume}{61} (\bibinfo{year}{2013}) \bibinfo{pages}{1600--1611}.
\bibitem[{Ha and Bobaru(2011)}]{peridynamicDamage}
\bibinfo{author}{Y.~D. Ha}, \bibinfo{author}{F.~Bobaru},
\newblock \bibinfo{title}{{Characteristics of dynamic brittle fracture captured
  with peridynamics}},
\newblock \bibinfo{journal}{Engineering Fracture Mechanics}
  \bibinfo{volume}{78} (\bibinfo{year}{2011}) \bibinfo{pages}{1156--1168}.
\bibitem[{Ma et~al.(2010)Ma, Wang, and Li}]{SPHdamage}
\bibinfo{author}{G.~W. Ma}, \bibinfo{author}{X.~J. Wang},
  \bibinfo{author}{Q.~M. Li},
\newblock \bibinfo{title}{{Modeling strain rate effect of heterogeneous
  materials using SPH method}},
\newblock \bibinfo{journal}{Rock Mechanics and Rock Engineering}
  \bibinfo{volume}{43} (\bibinfo{year}{2010}) \bibinfo{pages}{763--776}.
\bibitem[{Pedersen et~al.(2008)Pedersen, Simone, and Sluys}]{FEMdamageA}
\bibinfo{author}{R.~R. Pedersen}, \bibinfo{author}{A.~Simone},
  \bibinfo{author}{L.~J. Sluys},
\newblock \bibinfo{title}{{An analysis of dynamic fracture in concrete with a
  continuum visco-elastic visco-plastic damage model}},
\newblock \bibinfo{journal}{Engineering Fracture Mechanics}
  \bibinfo{volume}{75} (\bibinfo{year}{2008}) \bibinfo{pages}{3782--3805}.
\bibitem[{Zhou and Molinari(2004)}]{FEMdamageB}
\bibinfo{author}{F.~Zhou}, \bibinfo{author}{J.~F. Molinari},
\newblock \bibinfo{title}{{Dynamic crack propagation with cohesive elements: A
  methodology to address mesh dependency}},
\newblock \bibinfo{journal}{International Journal for Numerical Methods in
  Engineering} \bibinfo{volume}{59} (\bibinfo{year}{2004})
  \bibinfo{pages}{1--24}.
\bibitem[{Belytschko et~al.(2003)Belytschko, Chen, Xu, and Zi}]{XFEMdamage}
\bibinfo{author}{T.~Belytschko}, \bibinfo{author}{H.~Chen},
  \bibinfo{author}{J.~Xu}, \bibinfo{author}{G.~Zi},
\newblock \bibinfo{title}{{Dynamic crack propagation based on loss of
  hyperbolicity and a new discontinuous enrichment}},
\newblock \bibinfo{journal}{International Journal for Numerical Methods in
  Engineering} \bibinfo{volume}{58} (\bibinfo{year}{2003})
  \bibinfo{pages}{1873--1905}.
\bibitem[{Matou{\v{s}}(2003)}]{matous2003}
\bibinfo{author}{K.~Matou{\v{s}}},
\newblock \bibinfo{title}{{Damage evolution in particulate composite
  materials}},
\newblock \bibinfo{journal}{International Journal of Solids and Structures}
  \bibinfo{volume}{40} (\bibinfo{year}{2003}) \bibinfo{pages}{1489--1503}.
\bibitem[{Lee et~al.(2021)Lee, Ramos, and Matou{\v{s}}}]{matous2021}
\bibinfo{author}{S.~Lee}, \bibinfo{author}{K.~Ramos},
  \bibinfo{author}{K.~Matou{\v{s}}},
\newblock \bibinfo{title}{{Numerical study of damage in particulate composites
  during high-strain rate loading using novel damage model}},
\newblock \bibinfo{journal}{Mechanics of Materials} \bibinfo{volume}{160}
  (\bibinfo{year}{2021}).
\bibitem[{Kulkarni et~al.(2009)Kulkarni, Geubelle, and
  Matou{\v{s}}}]{matous2009}
\bibinfo{author}{M.~G. Kulkarni}, \bibinfo{author}{P.~H. Geubelle},
  \bibinfo{author}{K.~Matou{\v{s}}},
\newblock \bibinfo{title}{{Multi-scale modeling of heterogeneous adhesives:
  Effect of particle decohesion}},
\newblock \bibinfo{journal}{Mechanics of Materials} \bibinfo{volume}{41}
  (\bibinfo{year}{2009}) \bibinfo{pages}{573--583}.
\bibitem[{Abedi et~al.(2010)Abedi, Hawker, Haber, and
  Matou{\v{s}}}]{vspike2010}
\bibinfo{author}{R.~Abedi}, \bibinfo{author}{M.~A. Hawker},
  \bibinfo{author}{R.~B. Haber}, \bibinfo{author}{K.~Matou{\v{s}}},
\newblock \bibinfo{title}{{An adaptive spacetime discontinuous Galerkin method
  for cohesive models of elastodynamic fracture}},
\newblock \bibinfo{journal}{International Journal for Numerical Methods in
  Engineering} \bibinfo{volume}{81} (\bibinfo{year}{2010})
  \bibinfo{pages}{1207--1241}.
\bibitem[{Berger and Oliger(1984)}]{AMR_def}
\bibinfo{author}{M.~J. Berger}, \bibinfo{author}{J.~Oliger},
\newblock \bibinfo{title}{{Adaptive mesh refinement for hyperbolic partial
  differential equations}},
\newblock \bibinfo{journal}{Journal of Computational Physics}
  \bibinfo{volume}{53} (\bibinfo{year}{1984}) \bibinfo{pages}{484--512}.
\bibitem[{Fatkullin and Hesthaven(2001)}]{AMR_hesthaven}
\bibinfo{author}{I.~Fatkullin}, \bibinfo{author}{J.~S. Hesthaven},
\newblock \bibinfo{title}{{Adaptive high-order finite-difference method for
  nonlinear wave problems}},
\newblock \bibinfo{journal}{Journal of Scientific Computing}
  \bibinfo{volume}{16} (\bibinfo{year}{2001}) \bibinfo{pages}{47--67}.
\bibitem[{Brandt(1977)}]{multigrid_def}
\bibinfo{author}{A.~Brandt},
\newblock \bibinfo{title}{{Multi-level adaptive solutions to boundary-value
  problems}},
\newblock \bibinfo{journal}{Mathematics of Computation} \bibinfo{volume}{31}
  (\bibinfo{year}{1977}) \bibinfo{pages}{333--390}.
\bibitem[{Hackbusch(1978)}]{multigrid_hackbusch}
\bibinfo{author}{W.~Hackbusch},
\newblock \bibinfo{title}{{On the multi-grid method applied to difference
  equations}},
\newblock \bibinfo{journal}{Computing} \bibinfo{volume}{20}
  (\bibinfo{year}{1978}) \bibinfo{pages}{291--306}.
\bibitem[{Yushu and Matou{\v{s}}(2020)}]{matousDewen}
\bibinfo{author}{D.~Yushu}, \bibinfo{author}{K.~Matou{\v{s}}},
\newblock \bibinfo{title}{{The image-based multiscale multigrid solver,
  preconditioner, and reduced order model}},
\newblock \bibinfo{journal}{Journal of Computational Physics}
  \bibinfo{volume}{406} (\bibinfo{year}{2020}).
\bibitem[{Benek et~al.(1986)Benek, Steger, Dougherty, and Buning}]{chimera}
\bibinfo{author}{J.~A. Benek}, \bibinfo{author}{J.~L. Steger},
  \bibinfo{author}{F.~C. Dougherty}, \bibinfo{author}{P.~G. Buning},
  \bibinfo{title}{{Chimera: A Grid-Embedding Technique}},
  \bibinfo{type}{Technical Report}, AEDC-TR-85-64, \bibinfo{year}{1986}.
\bibitem[{Dong and Karniadakis(2003)}]{karniadakis}
\bibinfo{author}{S.~S. Dong}, \bibinfo{author}{G.~E. Karniadakis},
\newblock \bibinfo{title}{{P-refinement and P-threads}},
\newblock \bibinfo{journal}{Computer Methods in Applied Mechanics and
  Engineering} \bibinfo{volume}{192} (\bibinfo{year}{2003})
  \bibinfo{pages}{2191--2201}.
\bibitem[{Gui and Babu{\v{s}}ka(1986{\natexlab{a}})}]{FEM_hp1}
\bibinfo{author}{W.~Gui}, \bibinfo{author}{I.~Babu{\v{s}}ka},
\newblock \bibinfo{title}{{The h, p and h-p versions of the finite element
  method in 1 dimension part I. the error analysis of the p-version}},
\newblock \bibinfo{journal}{Numerische Mathematik} \bibinfo{volume}{49}
  (\bibinfo{year}{1986}{\natexlab{a}}) \bibinfo{pages}{577--612}.
\bibitem[{Gui and Babu{\v{s}}ka(1986{\natexlab{b}})}]{FEM_hp2}
\bibinfo{author}{W.~Gui}, \bibinfo{author}{I.~Babu{\v{s}}ka},
\newblock \bibinfo{title}{{The h, p and h-p versions of the finite element
  method in 1 dimension part II. the error analysis of the h-and h-p
  versions}},
\newblock \bibinfo{journal}{Numerische Mathematik} \bibinfo{volume}{49}
  (\bibinfo{year}{1986}{\natexlab{b}}) \bibinfo{pages}{613--657}.
\bibitem[{Rajagopal and Sivakumar(2007)}]{FEM_rh}
\bibinfo{author}{A.~Rajagopal}, \bibinfo{author}{S.~M. Sivakumar},
\newblock \bibinfo{title}{{A combined r-h adaptive strategy based on material
  forces and error assessment for plane problems and bimaterial interfaces}},
\newblock \bibinfo{journal}{Computational Mechanics} \bibinfo{volume}{41}
  (\bibinfo{year}{2007}) \bibinfo{pages}{49--72}.
\bibitem[{Jawerth and Sweldens(1994)}]{Jawerth1994}
\bibinfo{author}{B.~Jawerth}, \bibinfo{author}{W.~Sweldens},
\newblock \bibinfo{title}{{An overview of wavlet based multiresolution
  analyses}},
\newblock \bibinfo{journal}{SIAM Review} \bibinfo{volume}{36}
  (\bibinfo{year}{1994}) \bibinfo{pages}{377--412}.
\bibitem[{Schneider and Vasilyev(2010)}]{Schneider2010}
\bibinfo{author}{K.~Schneider}, \bibinfo{author}{O.~V. Vasilyev},
\newblock \bibinfo{title}{{Wavelet methods in computational fluid dynamics}},
\newblock \bibinfo{journal}{Annual Review of Fluid Mechanics}
  \bibinfo{volume}{42} (\bibinfo{year}{2010}) \bibinfo{pages}{473--503}.
\bibitem[{Liandrat and Tchamitchian(1990)}]{Liandrat1990}
\bibinfo{author}{J.~Liandrat}, \bibinfo{author}{P.~Tchamitchian},
  \bibinfo{title}{{Resolution of the 1D Regularized Burgers Equation Using a
  Spatial Wavelet Approximation}}, \bibinfo{type}{Technical Report}, ICASE
  90-83, \bibinfo{year}{1990}.
\bibitem[{Beylkin and Keiser(1997)}]{Beylkin1997}
\bibinfo{author}{G.~Beylkin}, \bibinfo{author}{J.~M. Keiser},
\newblock \bibinfo{title}{{On the adaptive numerical solution of nonlinear
  partial differential equations in wavelet bases}},
\newblock \bibinfo{journal}{Journal of Computational Physics}
  \bibinfo{volume}{132} (\bibinfo{year}{1997}) \bibinfo{pages}{233--259}.
\bibitem[{Bertoluzza(1996)}]{Bertoluzza1996a}
\bibinfo{author}{S.~Bertoluzza},
\newblock \bibinfo{title}{{Adaptive wavelet collocation method for the solution
  of Burgers equation}},
\newblock \bibinfo{journal}{Transport Theory and Statistical Physics}
  \bibinfo{volume}{25} (\bibinfo{year}{1996}) \bibinfo{pages}{339--352}.
\bibitem[{Ueno et~al.(2003)Ueno, Ide, and Okada}]{Ueno2003}
\bibinfo{author}{T.~Ueno}, \bibinfo{author}{T.~Ide},
  \bibinfo{author}{M.~Okada},
\newblock \bibinfo{title}{{A wavelet collocation method for evolution equations
  with energy conservation property}},
\newblock \bibinfo{journal}{Bulletin des Sciences Mathematiques}
  \bibinfo{volume}{127} (\bibinfo{year}{2003}) \bibinfo{pages}{569--583}.
\bibitem[{Quian and Weiss(1993)}]{Qian1993}
\bibinfo{author}{S.~Quian}, \bibinfo{author}{J.~Weiss},
\newblock \bibinfo{title}{{Wavelets and the numerical solution of partial
  differential equations}},
\newblock \bibinfo{journal}{Journal of Computational Physics}
  \bibinfo{volume}{106} (\bibinfo{year}{1993}) \bibinfo{pages}{155--175}.
\bibitem[{Kong et~al.(2016)Kong, Kougioumtzoglou, Spanos, and Li}]{Kong2016}
\bibinfo{author}{F.~Kong}, \bibinfo{author}{I.~A. Kougioumtzoglou},
  \bibinfo{author}{P.~D. Spanos}, \bibinfo{author}{S.~Li},
\newblock \bibinfo{title}{{Nonlinear systemresponse evolutionary power spectral
  density determination via a harmonicwavelets based galerkin technique}},
\newblock \bibinfo{journal}{International Journal for Multiscale Computational
  Engineering} \bibinfo{volume}{14} (\bibinfo{year}{2016})
  \bibinfo{pages}{255--272}.
\bibitem[{van Tuijl et~al.(2019)van Tuijl, Harnish, Matou{\v{s}}, Remmers, and
  Geers}]{Rody2019}
\bibinfo{author}{R.~van Tuijl}, \bibinfo{author}{C.~Harnish},
  \bibinfo{author}{K.~Matou{\v{s}}}, \bibinfo{author}{J.~Remmers},
  \bibinfo{author}{M.~Geers},
\newblock \bibinfo{title}{{Wavelet based reduced order models for
  microstructural analyses}},
\newblock \bibinfo{journal}{Computational Mechanics} \bibinfo{volume}{63}
  (\bibinfo{year}{2019}).
\bibitem[{Paolucci et~al.(2014{\natexlab{a}})Paolucci, Zikoski, and
  Wirasaet}]{Paolucci2014PT1}
\bibinfo{author}{S.~Paolucci}, \bibinfo{author}{Z.~J. Zikoski},
  \bibinfo{author}{D.~Wirasaet},
\newblock \bibinfo{title}{{WAMR: An adaptive wavelet method for the simulation
  of compressible reacting flow part I. accuracy and efficiency of algorithm}},
\newblock \bibinfo{journal}{Journal of Computational Physics}
  \bibinfo{volume}{272} (\bibinfo{year}{2014}{\natexlab{a}})
  \bibinfo{pages}{814--841}.
\bibitem[{Paolucci et~al.(2014{\natexlab{b}})Paolucci, Zikoski, and
  Grenga}]{Paolucci2014PT2}
\bibinfo{author}{S.~Paolucci}, \bibinfo{author}{Z.~J. Zikoski},
  \bibinfo{author}{T.~Grenga},
\newblock \bibinfo{title}{{WAMR: An adaptive wavelet method for the simulation
  of compressible reacting flow part II. the parallel algorithm}},
\newblock \bibinfo{journal}{Journal of Computational Physics}
  \bibinfo{volume}{272} (\bibinfo{year}{2014}{\natexlab{b}})
  \bibinfo{pages}{842--864}.
\bibitem[{Nejadmalayeri et~al.(2015)Nejadmalayeri, Vezolainen, Brown-Dymkoski,
  and Vasilyev}]{Nejadmalayeri2015}
\bibinfo{author}{A.~Nejadmalayeri}, \bibinfo{author}{A.~Vezolainen},
  \bibinfo{author}{E.~Brown-Dymkoski}, \bibinfo{author}{O.~V. Vasilyev},
\newblock \bibinfo{title}{{Parallel adaptive wavelet collocation method for
  PDEs}},
\newblock \bibinfo{journal}{Journal of Computational Physics}
  \bibinfo{volume}{298} (\bibinfo{year}{2015}) \bibinfo{pages}{237--253}.
\bibitem[{Dubos and Kevlahan(2013)}]{Dubos2013}
\bibinfo{author}{T.~Dubos}, \bibinfo{author}{N.~K. Kevlahan},
\newblock \bibinfo{title}{{A conservative adaptive wavelet method for the
  shallow-water equations on staggered grids}},
\newblock \bibinfo{journal}{Quarterly Journal of the Royal Meteorological
  Society} \bibinfo{volume}{139} (\bibinfo{year}{2013})
  \bibinfo{pages}{1997--2020}.
\bibitem[{Sakurai et~al.(2017)Sakurai, Yoshimatsu, Schneider, Farge, Morishita,
  and Ishihara}]{Sakurai2017}
\bibinfo{author}{T.~Sakurai}, \bibinfo{author}{K.~Yoshimatsu},
  \bibinfo{author}{K.~Schneider}, \bibinfo{author}{M.~Farge},
  \bibinfo{author}{K.~Morishita}, \bibinfo{author}{T.~Ishihara},
\newblock \bibinfo{title}{{Coherent structure extraction in turbulent channel
  flow using boundary adapted wavelets}},
\newblock \bibinfo{journal}{Journal of Turbulence} \bibinfo{volume}{18}
  (\bibinfo{year}{2017}) \bibinfo{pages}{352--372}.
\bibitem[{Fr{\"{o}}hlich and Schneider(1994)}]{Frohlich1994}
\bibinfo{author}{J.~Fr{\"{o}}hlich}, \bibinfo{author}{K.~Schneider},
\newblock \bibinfo{title}{{An Adaptive Wavelet Galerkin Algorithm for
  One-Dimensional and 2-Dimensional Flame Computations}},
\newblock \bibinfo{journal}{European Journal of Mechanics B-Fluids}
  \bibinfo{volume}{13} (\bibinfo{year}{1994}) \bibinfo{pages}{439--471}.
\bibitem[{Goedecker(2009)}]{Goedecker}
\bibinfo{author}{S.~Goedecker}, \bibinfo{title}{{Wavelets and their application
  for the solution of partial differential equations in physics}},
  \bibinfo{type}{Technical Report}, Max-Planck Institute for Solid State
  Research, \bibinfo{year}{2009}.
\bibitem[{Iqbal and Jeoti(2014)}]{Iqbal2014}
\bibinfo{author}{A.~Iqbal}, \bibinfo{author}{V.~Jeoti},
\newblock \bibinfo{title}{{An Improved Split-Step Wavelet Transform Method for
  Anomalous Radio Wave Propagation Modeling}},
\newblock \bibinfo{journal}{Radioengineering} \bibinfo{volume}{23}
  (\bibinfo{year}{2014}) \bibinfo{pages}{987--996}.
\bibitem[{Le and Caracoglia(2015)}]{Le2015}
\bibinfo{author}{T.~H. Le}, \bibinfo{author}{L.~Caracoglia},
\newblock \bibinfo{title}{{Reduced-order wavelet-Galerkin solution for the
  coupled, nonlinear stochastic response of slender buildings in transient
  winds}},
\newblock \bibinfo{journal}{Journal of Sound and Vibration}
  \bibinfo{volume}{344} (\bibinfo{year}{2015}) \bibinfo{pages}{179--208}.
\bibitem[{Lin and Zhou(2001)}]{Lin2001}
\bibinfo{author}{E.~B. Lin}, \bibinfo{author}{X.~Zhou},
\newblock \bibinfo{title}{{Connection coefficients on an interval and wavelet
  solutions of Burgers equation}},
\newblock \bibinfo{journal}{Journal of Computational and Applied Mathematics}
  \bibinfo{volume}{135} (\bibinfo{year}{2001}) \bibinfo{pages}{63--78}.
\bibitem[{Holmstr{\"{o}}m(1999)}]{Holmstrom1999}
\bibinfo{author}{M.~Holmstr{\"{o}}m},
\newblock \bibinfo{title}{{Solving hyperbolic PDEs using interpolating
  wavelets}},
\newblock \bibinfo{journal}{SIAM Journal on Scientific Computing}
  \bibinfo{volume}{21} (\bibinfo{year}{1999}) \bibinfo{pages}{405--420}.
\bibitem[{Daubechies(1992)}]{Daubechies10}
\bibinfo{author}{I.~Daubechies}, \bibinfo{title}{{Ten Lectures on Wavelets}},
  \bibinfo{publisher}{SIAM}, \bibinfo{year}{1992}.
  \DOIprefix\doi{10.1137/1.9781611970104}.
\bibitem[{Cohen et~al.(2000)Cohen, Dahmen, and Devore}]{Cohen2000_interval}
\bibinfo{author}{A.~Cohen}, \bibinfo{author}{W.~Dahmen},
  \bibinfo{author}{R.~Devore},
\newblock \bibinfo{title}{{Multiscale decompositions on bounded domains}},
\newblock \bibinfo{journal}{Transactions of the American Mathematical Society}
  \bibinfo{volume}{352} (\bibinfo{year}{2000}) \bibinfo{pages}{3651--3685}.
\bibitem[{Bacry et~al.(1992)Bacry, Mallat, and Papanicolaou}]{Bacry1992}
\bibinfo{author}{E.~Bacry}, \bibinfo{author}{S.~Mallat},
  \bibinfo{author}{G.~Papanicolaou},
\newblock \bibinfo{title}{{A Wavelet Based Space-Time Adaptive Numerical Method
  for Partial Differential Equations}},
\newblock \bibinfo{journal}{Mathematical Modelling and Numerical Analysis}
  \bibinfo{volume}{26} (\bibinfo{year}{1992}) \bibinfo{pages}{793--834}.
\bibitem[{De~Villiers et~al.(2004)De~Villiers, Goosen, and Herbst}]{interval}
\bibinfo{author}{J.~M. De~Villiers}, \bibinfo{author}{K.~M. Goosen},
  \bibinfo{author}{B.~M. Herbst},
\newblock \bibinfo{title}{{Dubuc-Deslauriers subdivision for finite sequences
  and interpolation wavelets on an interval}},
\newblock \bibinfo{journal}{SIAM Journal on Mathematical Analysis}
  \bibinfo{volume}{35} (\bibinfo{year}{2004}) \bibinfo{pages}{423--452}.
\bibitem[{Beylkin(1992)}]{Beylkin1992}
\bibinfo{author}{G.~Beylkin},
\newblock \bibinfo{title}{{On the representation of operators in bases of
  compactly supported wavelets}},
\newblock \bibinfo{journal}{SIAM Journal on Numerical Analysis}
  \bibinfo{volume}{6} (\bibinfo{year}{1992}) \bibinfo{pages}{1716--1740}.
\bibitem[{Sweldens(1998)}]{Sweldens1998}
\bibinfo{author}{W.~Sweldens},
\newblock \bibinfo{title}{{The lifting scheme: a construction of second
  generation wavelets}},
\newblock \bibinfo{journal}{SIAM Journal on Mathematical Analysis}
  \bibinfo{volume}{29} (\bibinfo{year}{1998}) \bibinfo{pages}{511--546}.
\bibitem[{Burgos et~al.(2013)Burgos, Cetale~Santos, and
  Silva}]{wavelet_FEM_beam}
\bibinfo{author}{R.~B. Burgos}, \bibinfo{author}{M.~A. Cetale~Santos},
  \bibinfo{author}{R.~R.~E. Silva},
\newblock \bibinfo{title}{{Deslauriers-Dubuc interpolating wavelet beam finite
  element}},
\newblock \bibinfo{journal}{Finite Elements in Analysis and Design}
  \bibinfo{volume}{75} (\bibinfo{year}{2013}) \bibinfo{pages}{71--77}.
\bibitem[{Fujii and Hoefer(2003)}]{maxwell}
\bibinfo{author}{M.~Fujii}, \bibinfo{author}{W.~J. Hoefer},
\newblock \bibinfo{title}{{Interpolating wavelet collocation method of time
  dependent Maxwell's equations: Characterization of electrically large optical
  waveguide discontinuities}},
\newblock \bibinfo{journal}{Journal of Computational Physics}
  \bibinfo{volume}{186} (\bibinfo{year}{2003}) \bibinfo{pages}{666--689}.
\bibitem[{Donoho(1992)}]{donoho_interpolating}
\bibinfo{author}{D.~L. Donoho}, \bibinfo{title}{{Interpolating Wavelet
  Transforms}}, \bibinfo{type}{Technical Report}, Stanford University,
  \bibinfo{year}{1992}.
\bibitem[{Harnish et~al.(2018)Harnish, Matou{\v{s}}, and Livescu}]{Harnish2018}
\bibinfo{author}{C.~Harnish}, \bibinfo{author}{K.~Matou{\v{s}}},
  \bibinfo{author}{D.~Livescu},
\newblock \bibinfo{title}{{Adaptive wavelet algorithm for solving nonlinear
  initial–boundary value problems with error control}},
\newblock \bibinfo{journal}{International Journal for Multiscale Computational
  Engineering} \bibinfo{volume}{16} (\bibinfo{year}{2018}).
\bibitem[{Harnish et~al.(2021)Harnish, Dalessandro, Matou{\v{s}}, and
  Livescu}]{Harnish2021}
\bibinfo{author}{C.~Harnish}, \bibinfo{author}{L.~Dalessandro},
  \bibinfo{author}{K.~Matou{\v{s}}}, \bibinfo{author}{D.~Livescu},
\newblock \bibinfo{title}{{A multiresolution adaptive wavelet method for
  nonlinear partial differential equations}},
\newblock \bibinfo{journal}{International Journal for Multiscale Computational
  Engineering} \bibinfo{volume}{19} (\bibinfo{year}{2021})
  \bibinfo{pages}{29--37}.
\bibitem[{Jameson(1993)}]{nasa_matrices}
\bibinfo{author}{L.~Jameson}, \bibinfo{title}{{On the Daubechies-based wavelet
  differentiation matrix}}, \bibinfo{type}{Technical Report}, ICASE 93-95,
  \bibinfo{year}{1993}.
\bibitem[{Dahmen et~al.(1999)Dahmen, Kunoth, and Urban}]{Dahmen1999}
\bibinfo{author}{W.~Dahmen}, \bibinfo{author}{A.~Kunoth},
  \bibinfo{author}{K.~Urban},
\newblock \bibinfo{title}{{Biorthogonal spline wavelets on the interval -
  stability and moment conditions}},
\newblock \bibinfo{journal}{Applied and Computational Harmonic Analysis}
  \bibinfo{volume}{6} (\bibinfo{year}{1999}) \bibinfo{pages}{132--196}.
\bibitem[{Farge(1992)}]{Farge1992}
\bibinfo{author}{M.~Farge},
\newblock \bibinfo{title}{{Wavelet transforms and their applications to
  turbulence}},
\newblock \bibinfo{journal}{Annual Review of Fluid Mechanics}
  \bibinfo{volume}{24} (\bibinfo{year}{1992}) \bibinfo{pages}{395--457}.
\bibitem[{Rioul(1992)}]{Rioul1992}
\bibinfo{author}{O.~Rioul},
\newblock \bibinfo{title}{{Simple regularity criteria for subdivision
  schemes}},
\newblock \bibinfo{journal}{SIAM Journal on Mathematical Analysis}
  \bibinfo{volume}{23} (\bibinfo{year}{1992}) \bibinfo{pages}{1544--1576}.
\bibitem[{Fehlberg(1970)}]{RKF45}
\bibinfo{author}{E.~Fehlberg},
\newblock \bibinfo{title}{{Classical fourth- and lower order Runge-Kutta
  formulas with stepsize control and their application to heat transfer
  problems}},
\newblock \bibinfo{journal}{Computing} \bibinfo{volume}{6}
  (\bibinfo{year}{1970}) \bibinfo{pages}{61--71}.
\bibitem[{Domingues et~al.(2009)Domingues, Gomes, Roussel, and
  Schneider}]{dtPercentChange}
\bibinfo{author}{M.~O. Domingues}, \bibinfo{author}{S.~M. Gomes},
  \bibinfo{author}{O.~Roussel}, \bibinfo{author}{K.~Schneider},
\newblock \bibinfo{title}{{Space-time adaptive multiresolution methods for
  hyperbolic conservation laws: Applications to compressible Euler equations}},
\newblock \bibinfo{journal}{Applied Numerical Mathematics} \bibinfo{volume}{59}
  (\bibinfo{year}{2009}) \bibinfo{pages}{2303--2321}.
\bibitem[{Pathria(1997)}]{noIntBC}
\bibinfo{author}{D.~Pathria},
\newblock \bibinfo{title}{{The correct formulation of intermediate boundary
  conditions for Runge-Kutta time integration of initial boundary value
  problems}},
\newblock \bibinfo{journal}{SIAM Journal on Scientific Computing}
  \bibinfo{volume}{18} (\bibinfo{year}{1997}) \bibinfo{pages}{1255--1266}.
\bibitem[{Hesthaven and Gottlieb(1996)}]{Hesthaven1996}
\bibinfo{author}{J.~S. Hesthaven}, \bibinfo{author}{D.~Gottlieb},
\newblock \bibinfo{title}{{A stable penalty method for the compressible
  Navier-Stokes equations: I. Open boundary conditions}},
\newblock \bibinfo{journal}{SIAM Journal on Scientific Computing}
  \bibinfo{volume}{17} (\bibinfo{year}{1996}) \bibinfo{pages}{579--612}.
\bibitem[{Schlick et~al.(2021)Schlick, Portillo-Ledesma, Blaszczyk,
  Dalessandro, Ghosh, Hackl, Harnish, Kotha, Livescu, Masud, Matou{\v{s}},
  Moyeda, Oskay, and Fish}]{HarnishSpecialIssue}
\bibinfo{author}{T.~Schlick}, \bibinfo{author}{S.~Portillo-Ledesma},
  \bibinfo{author}{M.~Blaszczyk}, \bibinfo{author}{L.~Dalessandro},
  \bibinfo{author}{S.~Ghosh}, \bibinfo{author}{K.~Hackl},
  \bibinfo{author}{C.~Harnish}, \bibinfo{author}{S.~Kotha},
  \bibinfo{author}{D.~Livescu}, \bibinfo{author}{A.~Masud},
  \bibinfo{author}{K.~Matou{\v{s}}}, \bibinfo{author}{A.~Moyeda},
  \bibinfo{author}{C.~Oskay}, \bibinfo{author}{J.~Fish},
\newblock \bibinfo{title}{{A multiscale vision -- illustrative applications
  from biology to engineering}},
\newblock \bibinfo{journal}{International Journal for Multiscale Computational
  Engineering} \bibinfo{volume}{19} (\bibinfo{year}{2021})
  \bibinfo{pages}{39--73}.
\bibitem[{Sedov(1946)}]{sedov}
\bibinfo{author}{L.~I. Sedov},
\newblock \bibinfo{title}{{Propagation of strong shock waves}},
\newblock \bibinfo{journal}{Journal of Applied Mathematics and Mechanics}
  \bibinfo{volume}{10} (\bibinfo{year}{1946}) \bibinfo{pages}{241--250}.
\bibitem[{{Ju J}(1989)}]{ju1989}
\bibinfo{author}{{Ju J}},
\newblock \bibinfo{title}{{On energy-based coupled elastoplastic damage
  theories: constitutive modeling and compuational aspects}},
\newblock \bibinfo{journal}{International Journal of Solids and Structures}
  \bibinfo{volume}{25} (\bibinfo{year}{1989}) \bibinfo{pages}{803--833}.
\bibitem[{Ahmed and Khan(2021)}]{justOmega}
\bibinfo{author}{A.~Ahmed}, \bibinfo{author}{R.~Khan},
\newblock \bibinfo{title}{{A phase field model for damage in asphalt
  concrete}},
\newblock \bibinfo{journal}{International Journal of Pavement Engineering}
  (\bibinfo{year}{2021}).
\bibitem[{Christman(1972)}]{GMelastic}
\bibinfo{author}{D.~R. Christman}, \bibinfo{title}{{Dynamic properties of poly
  (methylmethacrylate) (PMMA) (Plexiglas)}}, \bibinfo{type}{Technical Report},
  GM Technical Center, \bibinfo{address}{Warren, Michigan},
  \bibinfo{year}{1972}.
\bibitem[{Cheng et~al.(1990)Cheng, Miller, Manson, Hertzberg, and
  Sperling}]{calibrateDamage}
\bibinfo{author}{W.~M. Cheng}, \bibinfo{author}{G.~A. Miller},
  \bibinfo{author}{J.~A. Manson}, \bibinfo{author}{R.~W. Hertzberg},
  \bibinfo{author}{L.~H. Sperling},
\newblock \bibinfo{title}{{Mechanical behaviour of poly(methyl methacrylate) --
  Part 1 Tensile strength and fracture toughness}},
\newblock \bibinfo{journal}{Journal of Materials Science} \bibinfo{volume}{25}
  (\bibinfo{year}{1990}) \bibinfo{pages}{1917--1923}.
\bibitem[{Pluvinage and Capelle(2014)}]{PMMAlength}
\bibinfo{author}{G.~Pluvinage}, \bibinfo{author}{J.~Capelle},
\newblock \bibinfo{title}{{On characteristic lengths used in notch fracture
  mechanics}},
\newblock \bibinfo{journal}{International Journal of Fracture}
  \bibinfo{volume}{187} (\bibinfo{year}{2014}) \bibinfo{pages}{187--197}.
\bibitem[{Wu et~al.(2004)Wu, Ma, and Xia}]{PMMArateDep}
\bibinfo{author}{H.~Wu}, \bibinfo{author}{G.~Ma}, \bibinfo{author}{Y.~Xia},
\newblock \bibinfo{title}{{Experimental study of tensile properties of PMMA at
  intermediate strain rate}},
\newblock \bibinfo{journal}{Materials Letters} \bibinfo{volume}{58}
  (\bibinfo{year}{2004}) \bibinfo{pages}{3681--3685}.
\bibitem[{Chen et~al.(2002)Chen, Lu, and Cheng}]{PMMAsigMax}
\bibinfo{author}{W.~Chen}, \bibinfo{author}{F.~Lu}, \bibinfo{author}{M.~Cheng},
\newblock \bibinfo{title}{{Tension and compression tests of two polymers under
  quasi-static and dynamic loading}},
\newblock \bibinfo{journal}{Polymer Testing} \bibinfo{volume}{21}
  (\bibinfo{year}{2002}) \bibinfo{pages}{113--121}.
\bibitem[{Miller et~al.(1999)Miller, Freund, and Needleman}]{Miller1999}
\bibinfo{author}{O.~Miller}, \bibinfo{author}{L.~B. Freund},
  \bibinfo{author}{A.~Needleman},
\newblock \bibinfo{title}{{Energy dissipation in dynamic fracture of brittle
  materials}},
\newblock \bibinfo{journal}{Modelling and Simulation in Materials Science and
  Engineering} \bibinfo{volume}{7} (\bibinfo{year}{1999})
  \bibinfo{pages}{573--586}.
\bibitem[{Zhang et~al.(2007)Zhang, Paulino, and Celes}]{Zhang2007}
\bibinfo{author}{Z.~J. Zhang}, \bibinfo{author}{G.~H. Paulino},
  \bibinfo{author}{W.~Celes},
\newblock \bibinfo{title}{{Extrinsic cohesive modelling of dynamic fracture and
  microbranching instability in brittle materials}},
\newblock \bibinfo{journal}{International Journal for Numerical Methods in
  Engineering} \bibinfo{volume}{72} (\bibinfo{year}{2007})
  \bibinfo{pages}{893--923}.
\bibitem[{Xu and Needleman(1994)}]{Needleman1994}
\bibinfo{author}{X.-P. Xu}, \bibinfo{author}{A.~Needleman},
\newblock \bibinfo{title}{{Numerical simulations of fast crack growth in
  brittle solids}},
\newblock \bibinfo{journal}{Journal of the Mechanics and Physics of Solids}
  \bibinfo{volume}{42} (\bibinfo{year}{1994}) \bibinfo{pages}{1397--1434}.
\bibitem[{Stewart et~al.(2017)Stewart, Welch, Plale, Fox, Pierce, and
  Sterling}]{IUpti}
\bibinfo{author}{C.~Stewart}, \bibinfo{author}{V.~Welch},
  \bibinfo{author}{B.~Plale}, \bibinfo{author}{G.~Fox},
  \bibinfo{author}{M.~Pierce}, \bibinfo{author}{T.~Sterling},
  \bibinfo{title}{{Indiana University Pervasive Technology Institute}},
  \bibinfo{year}{2017}. \URLprefix \url{https://pti.iu.edu/}.
  \DOIprefix\doi{10.5967/K8G44NGB}.

\end{thebibliography}
\end{document}